\def\coralreport{1}

\documentclass[11pt]{article}

\usepackage{ifthen,amsmath,amsfonts,amssymb,amsthm,hyperref,algorithm,algpseudocode,mathtools,dpr_fec}
\newtheorem{theorem}{Theorem}
\newtheorem{lemma}{Lemma}
\newtheorem{corollary}{Corollary}
\allowdisplaybreaks

\ifthenelse{\coralreport = 1}{
\usepackage{isetechreport}
\def\reportyear{22T}
\def\reportno{001}
\def\revisionno{0}
\def\originaldate{\today}
\def\revisiondate{}
\coraltrue
\cvcrfalse

}{}

\begin{document}

\title{Worst-Case Complexity of an SQP Method for Nonlinear Equality Constrained Stochastic Optimization}






\ifthenelse{\coralreport = 1}{

\author{Frank E.~Curtis\thanks{E-mail: \texttt{frank.e.curtis@lehigh.edu}, supported by NSF Grant CCF-2008484 and ONR Grant N00014-21-1-2532}}
\author{Michael J.~O'Neill\thanks{E-mail: \texttt{moneill@lehigh.edu}, supported by the CI Fellows Program}}
\author{Daniel P.~Robinson\thanks{E-mail: \texttt{daniel.p.robinson@lehigh.edu}, supported by supported by ONR Grant N00014-21-1-2532}}
\affil{Department of Industrial and Systems Engineering, Lehigh University, USA}

\titlepage

}{

\author{Frank E.~Curtis \and Michael J.~O'Neill \and Daniel P.~Robinson}
\institute{Department of Industrial and Systems Engineering, Lehigh University, USA\\
\email{frank.e.curtis@lehigh.edu}, \email{moneill@lehigh.edu}, \email{daniel.p.robinson@lehigh.edu}
}

\date{\today}

}

\maketitle

\begin{abstract}
  A worst-case complexity bound is proved for a sequential quadratic optimization (commonly known as SQP) algorithm that has been designed for solving optimization problems involving a stochastic objective function and deterministic nonlinear equality constraints.  Barring additional terms that arise due to the adaptivity of the monotonically nonincreasing merit parameter sequence, the proved complexity bound is comparable to that known for the stochastic gradient algorithm for unconstrained nonconvex optimization.  The overall complexity bound, which accounts for the adaptivity of the merit parameter sequence, shows that a result comparable to the unconstrained setting (with additional logarithmic factors) holds with high probability.
\ifthenelse{\coralreport = 0}{
\bigskip\noindent
{\bf Keywords:} BS  
}{}
\end{abstract}

\newcommand{\flow}{f_{\text{\rm low}}}
\newcommand{\Prob}{\mathbb{P}}
\newcommand{\kappadec}{\kappa_{\text{\rm dec}}}
\newcommand{\bK}{\bar{K}}
\newcommand{\balpha}{\bar{\alpha}}
\newcommand{\btau}{\bar{\tau}}
\newcommand{\bd}{\bar{d}}
\newcommand{\bg}{\bar{g}}
\newcommand{\bxi}{\bar{\xi}}
\newcommand{\hdelta}{\hat{\delta}}

\newcommand{\xiktri}{\xi_k^{\rm trial}}
\newcommand{\ktaudec}{\kappa_{\tau_{\rm dec}}}
\newcommand{\kmax}{k_{\max}}
\newcommand{\smax}{s_{\max}}
\newcommand{\tauktri}{\tau_k^{\rm trial}}

\newcommand{\dktrue}{d_k^{\text{\rm true}}}
\newcommand{\uktrue}{u_k^{\text{\rm true}}}
\newcommand{\wktrue}{w_k^{\text{\rm true}}}
\newcommand{\yktrue}{y_k^{\text{\rm true}}}
\newcommand{\ykstartrue}{y_{k^*}^{\text{\rm true}}}
\newcommand{\tktritrue}{\tau_k^{\rm \text{trial},\text{\rm true}}}
\newcommand{\Tktritrue}{\Tcal_k^{\rm \text{trial},\text{\rm true}}}
\newcommand{\Dktrue}{D_k^{\text{\rm true}}}
\newcommand{\Yktrue}{Y_k^{\text{\rm true}}}
\newcommand{\Uktrue}{U_k^{\text{\rm true}}}
\newcommand{\Uk}{U_k}
\newcommand{\Vk}{V_k}

\newcommand{\strue}{s^{\text{\rm true}}}

\newcommand{\Lbad}{\Lcal_{\text{\rm bad}}}
\newcommand{\Lgood}{\Lcal_{\text{\rm good}}}
\newcommand{\Ebad}{E_{\text{\rm bad}}}
\newcommand{\EbadB}{E_{\text{\rm bad},B}}
\newcommand{\Cdec}{C_{\text{\rm dec}}}
\newcommand{\Idec}{\Ical_{\text{\rm dec}}}
\newcommand{\ktau}{k_{\tau}}
\newcommand{\tcalktritrue}{\Tcal_k^{\rm \text{trial},\text{\rm true}}}

\def\Ek#1{E_{k,#1}}
\def\tred#1{\textcolor{red}{#1}}
\DeclarePairedDelimiter\ceil{\lceil}{\rceil}

\section{Introduction} \label{sec:intro}

We present a worst-case complexity analysis of an algorithm for minimizing a smooth objective function subject to nonlinear equality constraints.  (Due to the nature of the algorithm, this worst-case complexity analysis holds in terms of iterations, function evaluations, and derivative evaluations.)  Problems of this type arise in various important applications throughout science and engineering, including optimal control, PDE-constrained optimization, and resource allocation \cite{DPBertsekas_1998,JTBetts_2010,FSKupfer_EWSachs_1992,TRees_HSDollar_AJWathen_2010}.  However, unlike the vast majority of the literature on equality constrained optimization, the algorithm that we consider has been designed to solve problems in which the objective function is stochastic, in the sense that it is defined by the expectation of a function that has a random variable argument.  The algorithm that we consider assumes that evaluations of the objective function and its gradient are intractable to obtain, but that it has access to (unbiased) stochastic gradient estimates.

A few algorithms have been proposed recently for solving problems of this type.  These approaches fall into two categories: penalty methods \cite{CChen_FTung_NVedula_2018,YNandwani_PAbhishek_Mausam_PSingla,SKRoy_ZMhammedi_MHarandi_2018} (which includes the class of augmented Lagrangian methods) and sequential quadratic optimization (commonly known as SQP) methods \cite{ABerahas_FECurtis_DPRobinson_BZhou_2021,SNa_MAnitescu_KMladen_2021}.  Penalty methods aim to solve the constrained optimization problem by adding a term to the objective function, weighted by a penalty parameter, that penalizes constraint violation.  Unconstrained optimization techniques are then applied to minimize the resulting penalty function, after which the penalty parameter may be modified and the minimization is performed again in an iterative manner until a solution is obtained that (approximately) satisfies the original constraints.  Methods of this type perform well in some situations, but in others they perform poorly, e.g., due to ill-conditioning and/or nonsmoothness of the subproblems. Such methods also often suffer due to their sensitivity to the particular scheme used for updating the penalty parameter.

In practice in both deterministic and stochastic optimization contexts, penalty methods are frequently outperformed by SQP methods.  Indeed, it is commonly accepted in the deterministic optimization literature that a state-of-the-art algorithm is an SQP method that chooses stepsizes based on a line search applied to a merit function.  In this deterministic setting, such an algorithm is intimately connected with applying Newton's method to the first-order primal-dual necessary conditions for optimality of the problem~\cite{RBWilson_1963}.

In this paper, we present a worst-case complexity analysis of the SQP method proposed in \cite{ABerahas_FECurtis_DPRobinson_BZhou_2021}, which can be seen as an extension of an SQP method from the deterministic to the stochastic setting.  A consequence of our analysis is that, in an idealized setting in which one knows a threshold for the merit parameter beyond which the merit function is \emph{exact} \cite{HanM79}, the number of iterations required until the method generates a point at which first-order necessary conditions for optimality hold in expectation with accuracy $\varepsilon \in (0,\infty)$ is $\Ocal(\varepsilon^{-4})$.  This is arguably the best result that one can expect given this is the same bound proved to hold for a stochastic gradient method employed to solve an unconstrained nonconvex problem~\cite{DDavis_DDrusvyatskiy_2019}.  However, \emph{our analysis does not only consider this idealized setting}; we go further and prove a worst-case complexity bound for the algorithm when the merit parameter threshold is unknown and the algorithm adaptively updates a monotonically nonincreasing merit parameter sequence.  We prove under reasonable assumptions that the aforementioned worst-case bound, with additional logarithmic factors, holds with high probability.  The high-probability aspect of this result arises purely due to the uncertainty of the behavior of the adaptive merit parameter sequence, and does not reflect any uncertainty of the behavior of the method during situations in which the merit parameter sequence remains constant.

To the best of our knowledge, ours is the first worst-case complexity result for an SQP algorithm that operates in the \emph{highly stochastic regime} (where one merely presumes that the stochastic gradient estimates have bounded variance) for solving stochastic optimization problems involving deterministic nonlinear equality constraints.  Prior to this work, the only known complexity results for stochastic constrained optimization were for algorithms for solving problems with simple constraint sets that enable projection-based methods \cite{DDavis_DDrusvyatskiy_2019,SGhadimi_GLan_HZhang_2016} and Frank-Wolfe type methods \cite{EHazan_HLuo_2016}.  (One exception is a complexity bound proved for the SQP algorithm proposed in \cite{SNa_MAnitescu_KMladen_2021}, although that result only holds for the idealized setting in which the algorithm has \emph{a priori} knowledge of a threshold for the merit function parameter.)  Our analysis focuses a great deal on the complications that arise due to the adaptivity of the merit parameter sequence, which essentially means that the algorithm in our consideration is aiming to reduce a merit function that \emph{changes} during the optimization process.  Hence, many aspects of our analysis are quite distinct from the analyses that have been presented for stochastic gradient methods in the context of unconstrained optimization or optimization over simple constraint sets, for which the tool for measuring the progress of an algorithm---namely, the objective function itself---remains the same throughout the optimization.

\subsection{Problem formulation}

The algorithm that we consider is designed to solve problems of the form
\bequation \label{eq:fdef}
  \min_{x \in \R{n}} f(x)\ \ \st\ \ c(x) = 0,\ \ \text{with}\ \ f(x) = \E[F(x,\omega)],
\eequation
where $f : \R{n} \rightarrow \R{}$, $c : \R{n} \rightarrow \R{m}$, $\omega$ is a random variable with associated probability space $(\Omega,\mathcal{F}, P)$, $F : \R{n} \times \Omega \rightarrow \R{}$, and $\E$ represents expectation with respect to $P$.  In particular, following \cite{ABerahas_FECurtis_DPRobinson_BZhou_2021}, we make the following assumption about problem \eqref{eq:fdef} and the algorithm that we analyze (stated as Algorithm~\ref{alg:ssqp} on page~\pageref{alg:ssqp}), which in any run generates a sequence of iterates $\{x_k\} \subset \R{n}$.

\begin{assumption}\label{assum:fcsmooth}
  Let $\Xcal \subseteq \R{n}$ be an open convex set that contains $\{x_k\}$ for all realizations of Algorithm \ref{alg:ssqp}.  The objective function $f : \R{n} \rightarrow \R{}$ is continuously differentiable and bounded below by $\flow \in \R{}$ over $\Xcal$ and the corresponding gradient function $\nabla f: \R{n} \rightarrow \R{n}$ is bounded and Lipschitz continuous with constant $L \in (0,\infty)$ over $\Xcal$. The constraint function $c : \R{n} \rightarrow \R{m}$ $($where $m \leq n$$)$ and the corresponding Jacobian function $J := \nabla c^\top : \R{n} \rightarrow \R{m\times n}$ are bounded over $\Xcal$, each gradient function $\nabla c_i : \R{n} \rightarrow \R{n}$ is Lipschitz continuous with constant $\gamma_i$ over $\Xcal$ for all $i \in \{1,\dots,m\}$, and the singular values of $J \equiv \nabla c^\top$ are bounded below and away from zero over $\Xcal$.
\end{assumption}

A consequence of Assumption~\ref{assum:fcsmooth} is that there exists $\kappa_g \in \R{}_{>0}$ such that
\bequation\label{eq.kappa_g}
  \|\nabla f(x_k)\| \leq \kappa_g\ \ \text{for any $k \in \N{}$ in any realization of Algorithm~\ref{alg:ssqp}}.
\eequation
Defining the Lagrangian $\ell : \R{n} \times \R{m} \rightarrow \R{}$ corresponding to \eqref{eq:fdef} by $\ell(x,y) := f(x) + c(x)^\top y$, first-order primal-dual stationarity conditions for~\eqref{eq:fdef}, which are necessary for optimality under Assumption~\ref{assum:fcsmooth}, are given by
\bequation \label{eq:optconds}
  0 =
  \left[ \begin{matrix}
    \nabla_x \ell(x,y) \\
    \nabla_y \ell(x,y)
  \end{matrix} \right]
  =
  \left[ \begin{matrix}
    \nabla f(x) + \nabla c(x) y \\
    c(x)
  \end{matrix} \right].
\eequation

\subsection{Notation}

We adopt the notation that $\|\cdot\|$ denotes the $\ell_2$-norm for vectors and the vector-induced $\ell_2$-norm for matrices.  We denote by $\mathbb{S}^{n}$ the set of $n \times n$ dimensional real symmetric matrices. The set of nonnegative integers is denoted as $\N{} := \{0,1,2,\dots,\}$. For any integer $k \in \N{}$, we use $[k]$ to denote the subset of nonnegative integers up to $k$, namely, $[k] := \{0,\dots,k\}$.  Correspondingly, to represent a set of vectors $\{v_0,\dots,v_k\}$, we define $v_{[k]} := \{v_0,\dots,v_k\}$.

Given $\phi : \R{} \to \R{}$ and $\varphi : \R{} \to [0,\infty)$, we write $\phi(\cdot) = \Ocal(\varphi(\cdot))$ to indicate that $|\phi(\cdot)| \leq c\varphi(\cdot)$ for some $c \in (0,\infty)$.  Similarly, we write $\phi(\cdot) = \widetilde\Ocal(\varphi(\cdot))$ to indicate that $|\phi(\cdot)| \leq c\varphi(\cdot)|\log^{\cbar}(\cdot)|$ for some $c \in (0,\infty)$ and $\cbar \in (0,\infty)$.  In this manner, one finds that $\Ocal(\varphi(\cdot)|\log^{\cbar}(\cdot)|) \equiv \widetilde\Ocal(\varphi(\cdot))$ for any $\cbar \in (0,\infty)$.

Algorithm~\ref{alg:ssqp} is iterative, generating in each realization a sequence $\{x_k\}$.  (See Section~\ref{subsec:process} for a complete description of the stochastic process generated by the algorithm.)  For our analysis, we also append the iteration number to other quantities corresponding to each iteration, e.g., $f_k := f(x_k)$ for all $k \in \N{}$.

When discussing stochastic quantities, we use capital letters to denote random variables and corresponding lower case letters to denote a realization of a random variable. For example, a stochastic gradient in iteration $k \in \N{}$ is denoted as $G_k$, a realization of which is written as $g_k$.

\subsection{Outline}

Section~\ref{sec:complexitypreview} provides a worst-case complexity result for the algorithm from \cite{ABerahas_FECurtis_DPRobinson_BZhou_2021} for the deterministic setting, and uses this result and further commentary to provide an overview of our main result for the stochastic setting.  Details of the algorithm for the stochastic setting are presented in Section~\ref{sec.algorithm}, followed by our main result and analysis, which are provided in Section~\ref{sec:complexity}.  Finally, we provide concluding thoughts and mention future directions in Section~\ref{sec:conclusion}.

\section{Outline of Main Results} \label{sec:complexitypreview}

Our algorithm of consideration, namely, Algorithm~3.1 in \cite{ABerahas_FECurtis_DPRobinson_BZhou_2021}, is derived from Algorithm~2.1 in \cite{ABerahas_FECurtis_DPRobinson_BZhou_2021}, which is proposed and analyzed for the deterministic setting as a precursor for the stochastic setting (of Algorithm~3.1 in \cite{ABerahas_FECurtis_DPRobinson_BZhou_2021}).  In the deterministic algorithm, the $k$th search direction $d_k \in \R{n}$ is computed by solving an optimization subproblem defined by a quadratic approximation of the objective function and an affine approximation of the constraints using derivative information at the current iterate $x_k \in \R{n}$.  This computation also results in a Lagrange multiplier vector $y_k \in \R{m}$.  The subsequent iterate is set by $x_{k+1} \gets x_k + \alpha_k d_k$, where $\alpha_k \in (0,\infty)$ is a stepsize determined by a procedure to reduce the merit function $\phi : \R{n} \times (0,\infty) \to \R{}$ defined by
\bequation \label{eq:phidef}
  \phi(x,\tau) = \tau f(x) + \|c(x)\|_1.
\eequation
In particular, based on properties of the search direction $d_k$, a value of the merit parameter $\tau_k \in (0,\tau_{k-1}]$ is set by the algorithm, after which $\alpha_k \in (0,\infty)$ is computed to ensure that $\phi(x_k,\tau_k) - \phi(x_{k+1},\tau_k)$ is sufficiently positive.

\subsection{Complexity of the Deterministic Algorithm}

To motivate our main result for the stochastic setting, it is instructive to state a worst-case complexity bound for the deterministic algorithm.  Such a result is the following; further details and a proof are provided in Appendix~\ref{app.deterministic}.

\btheorem\label{th.deterministic}
  Consider Algorithm~2.1 in \cite{ABerahas_FECurtis_DPRobinson_BZhou_2021} and suppose that Assumption~\ref{assum:fcsmooth} holds along with Assumption~2.4 from \cite{ABerahas_FECurtis_DPRobinson_BZhou_2021}.  Let $\tau_{-1} \in \R{}_{>0}$ be the initial value of the merit parameter sequence and let $\tau_{\min} \in (0,\tau_{-1}]$ be a positive lower bound for the merit parameter sequence $($the existence of which follows from Lemma~2.16 in~\cite{ABerahas_FECurtis_DPRobinson_BZhou_2021}$)$.  Then, for any $\varepsilon \in (0,1)$, there exists $(\kappa_1,\kappa_2) \in \R{}_{>0} \times \R{}_{>0}$ such that the algorithm reaches an iterate $(x_k,y_k) \in \R{n} \times \R{m}$ satisfying
  \bequation\label{eq:eps1lon}  
    \|g_k + J_k^\top y_k\| \leq \varepsilon\ \ \text{and}\ \ \sqrt{\|c_k\|_1} \leq \varepsilon
  \eequation
  in a number of iterations no more than
  \bequation\label{eq.deterministic_complexity}
    \(\frac{\tau_{-1}(f_0 - \flow) + \|c_0\|_1}{\min\{\kappa_1,\tau_{\min} \kappa_2\}}\) \varepsilon^{-2}.
  \eequation
\etheorem

Theorem~\ref{th.deterministic} is not surprising.  After all, such a complexity bound of $\Ocal(\varepsilon^{-2})$ is well-known for gradient-based algorithms in unconstrained nonconvex optimization.  Since Algorithm~2.1 in \cite{ABerahas_FECurtis_DPRobinson_BZhou_2021} and its corresponding analysis do not exploit the use of exact higher-order derivative information, this complexity bound is on the order of what could be expected for such a method.

\subsection{Preview of the Complexity of the Stochastic Algorithm}\label{sec.preview}

Moving to the stochastic setting, there are a few major technical hurdles that need to be addressed, all of which relate to the adaptivity of the merit parameter sequence.  In particular, the analysis for the deterministic setting relies heavily on the facts that $(i)$ each step of the algorithm yields a sufficient reduction in the merit function for the current value of the merit parameter, $(ii)$ each such reduction in the merit function can be tied to a first-order primal-dual stationarity measure for the current iterate, and $(iii)$ under Assumption~\ref{assum:fcsmooth}, one can be certain of the existence of a positive lower bound for the merit parameter sequence.  This lower bound for the merit parameter is referenced directly in the proof of the worst-case bound for the deterministic algorithm; in particular, it is shown (see Lemma~\ref{lem:suffdecrease} and the beginning of the proof of Theorem~\ref{th.deter_complex}) that the improvement in the merit function from any iterate that is not $\varepsilon$-stationary is at least proportional to $\min\{1,\tau_{\min}\} \varepsilon^{-2}$, even if the current value of the merit parameter is greater than $\tau_{\min}$.  Unfortunately, these properties of the steps and merit parameter sequence are not certain in the stochastic setting.  For example, as discussed in \cite{ABerahas_FECurtis_DPRobinson_BZhou_2021}, it is possible---even under Assumption~\ref{assum:fcsmooth}---for the merit parameter sequence to vanish or for it to eventually remain constant at a value that is not sufficiently small, and for there to be iterations in which the expected reduction in the merit function cannot be tied to a first-order primal-dual stationarity measure.  As a result, we have had to devise new analytical approaches that confront the fact that $\{\tau_k\}$ is a random process, the ultimate behavior of which is uncertain.

To aid the reader, we provide here an overview and commentary about our ultimate complexity bound; see Corollary~\ref{coro:complexity} on page~\pageref{coro:complexity}.  Our result is proved under Assumption~\ref{assum:fcsmooth} along with others that are introduced in the subsequent sections. For one thing, as is common in SQP methods for deterministic optimization, we assume that the subproblem defining the search direction in each iteration is defined by a matrix that is positive definite in the null space of the constraint Jacobian; see Assumption~\ref{assum:Hk} on page~\pageref{assum:Hk}.  We also assume, as is common for stochastic gradient methods, that the stochastic gradient estimates are unbiased with variance bounded by $M \in (0,\infty)$, along with some related assumptions; see Assumption~\ref{assum:eventE} on page~\pageref{assum:eventE}.  Furthermore, our analysis conditions on the occurrence of an event that we call $E$ (see \eqref{def:E}); this event captures situations in which, over a total of $\kmax+1 \in \N{}$ iterations, the merit parameter is reduced at most $\smax \in [\kmax]$ times and the merit parameter is bounded below by $\tau_{\min} \in (0,\infty)$.  Under these conditions, our main complexity result shows that, within $\kmax + 1$ iterations, it holds with probability $1-\delta \in (0,1)$ that the algorithm generates $x_{k^*} \in \R{n}$ corresponding to which there exists an associated Lagrange multiplier $\ykstartrue \in \R{m}$ such that
\bsubequations\label{eq.overview}
  \begin{align}
     &\ \E[\|\nabla f_{k^*} + J_{k^*}^\top \ykstartrue\|^2 + \|c_{k^*}\|_1 | E] \nonumber \\
    =&\ \Ocal \Bigg(\frac{\tau_{-1}(f_0 - \flow) + \|c_0\|_1 + M}{\sqrt{\kmax+1}} \label{eq.like_sg} \\
     &\qquad + \frac{(\tau_{-1} - \tau_{\min})(\smax\log(\kmax) + \log(1/\delta))}{\sqrt{\kmax+1}}\Bigg). \label{eq.new_terms}
  \end{align}
\esubequations
This form of the result is commonly called a convergence rate since it bounds the expected stationarity error from above by a function that decreases with the number of iterations performed, namely, $\kmax+1$.  This bound can be used to form a worst-case complexity result.  Specifically, the result above and Jensen's inequality imply that, within $\kmax+1$ iterations and as long as $\smax = \Ocal(\log(\kmax))$ (more on this below), it holds with probability $1-\delta$ that the algorithm requires at most $\widetilde\Ocal(\varepsilon^{-4})$ iterations to generate~$x_{k^*}$ with corresponding $\ykstartrue$ such that $\E[\|\nabla f_{k^*} + J_{k^*}^\top \ykstartrue\| |E] \leq \varepsilon$ and $\E[\sqrt{\|c_{k^*}\|_1}|E] \leq \varepsilon$.

The first three quantities on the right-hand side of the convergence rate, namely, in \eqref{eq.like_sg}, representing the initial objective function gap, initial constraint violation, and the variance of the stochastic gradient estimates, mirror the presence of similar terms that appear for comparable results for the stochastic gradient method in an unconstrained or simple-constraint-set setting.  The final term in \eqref{eq.new_terms}, on the other hand, as well as the fact that the result is stated as a high-probability result, are unique to our setting and arise due to the adaptivity of the merit parameter sequence.  If one were to have prior knowledge of $\tau_{\min}$, then one could set $\tau_{-1} = \tau_{\min}$ (and disable the update mechanism for the merit parameter in the algorithm), in which case our analysis would show that the expected stationarity error is bounded above by~\eqref{eq.like_sg} (surely, not only with high probability).

In the context of an adaptive merit parameter sequence, the particular form of our complexity result depends on the magnitude of $\smax$ relative to $\kmax+1$, i.e., the bound on the number of times that the merit parameter is decreased relative to the total number of iterations performed.  One setting in which our result is relatively straightforward is when, over all realizations of the algorithm, the differences between the stochastic gradient estimates and the true gradients are bounded deterministically, in which case $\smax$ is bounded by a value that is independent from $\kmax+1$; this follows from a deterministic lower bound on $\tau_{\min}$ \cite[Proposition 3.18]{ABerahas_FECurtis_DPRobinson_BZhou_2021} and the fact that whenever the merit parameter is decreased, it is done so by a constant factor. Beyond this setting, for another concrete example of a situation in which $\smax$ is guaranteed to be sufficiently small relative to $\kmax$, we prove in Section~\ref{subsec:subgauss} that if the distributions of the stochastic gradient estimates are sub-Gaussian, then with probability $1-\delta$ one finds that $\smax = \Ocal(\log(\log(\tfrac{\kmax}{\delta})))$, meaning that our proved convergence rate is not ruined by the term in~\eqref{eq.new_terms}.

\section{Algorithm}\label{sec.algorithm}

For ease of reference, in this section we present Algorithm~3.1 from \cite{ABerahas_FECurtis_DPRobinson_BZhou_2021}, which is designed to solve problems of the form~\eqref{eq:fdef} and is our focus for the remainder of the paper.  In the spirit of an SQP method, the algorithm computes a search direction $d_k$ and Lagrange multiplier vector $y_k$ in iteration $k \in \N{}$ by solving
\begin{equation} \label{eq:sqpsubproblem}
  \min_{d \in \R{n}}\ f_k + g_k^\top d + \thalf d^\top H_k d\ \ \st\ \ c_k + J_k d = 0,
\end{equation}
where $H_k \in \mathbb{S}^n$ is chosen independently from $g_k$, and we remind the reader that $g_k$ is a realization of the stochastic gradient $G_k$. Under Assumption~\ref{assum:fcsmooth} and the following Assumption~\ref{assum:Hk} (that we make throughout the remainder of the paper), the solution of \eqref{eq:sqpsubproblem} can be obtained from the unique solution of
\begin{equation} \label{eq:sqpsystem}
    \left[
    \begin{matrix}
    H_k & J_k^\top \\
    J_k & 0
    \end{matrix}
    \right]
    \left[
    \begin{matrix}
    d_k \\ y_k
    \end{matrix}
    \right]
    = -\left[
    \begin{matrix}
    g_k \\ c_k
    \end{matrix}
    \right].
\end{equation}

\begin{assumption} \label{assum:Hk}
  The sequence $\{\|H_k\|\}$ is bounded by $\kappa_H \in \R{}_{>0}$.  In addition, there exists $\zeta \in \R{}_{>0}$ such that, for all $k \in [\kmax]$, the matrix $H_k \in \mathbb{S}^n$ has the property that $u^\top H_k u \geq \zeta \|u\|_2^2$ for all $u \in \R{n}$ such that $J_k u = 0$.
\end{assumption}

After computation of $(d_k,y_k)$, the remainder of the $k$th iteration involves $(i)$ updating the merit parameter, $(ii)$ updating an auxiliary parameter needed for the stepsize computation, and $(iii)$ computing a positive stepsize.  These algorithmic components are designed with the aim of yielding a sufficiently positive reduction in a model of the merit function, which in turn is aimed at yielding a sufficiently positive reduction in the merit function itself.  The algorithm employs the model $q : \R{n} \times \R{}_+ \times \R{n} \times \mathbb{S}^n \times \R{n} \rightarrow \R{}$ defined by
\bequationNN
  q(x, \tau, g, H, d) = \tau (f(x) + g^\top d + \thalf \max\{d^\top H d,0\}) + \|c(x) + J(x) d\|_1,
\eequationNN
and the reduction function $\Delta q : \R{n} \times \R{}_+ \times \R{n} \times \mathbb{S}^n \times \R{n} \rightarrow \R{}$, for a given $d \in \R{n}$ satisfying $c(x) + J(x) d = 0$, defined by
\begin{equation} \label{eq:deltadef}
  \baligned
    \Delta q(x, \tau, g, H, d)
      :=&\ q(x, \tau, g, H, 0) - q(x, \tau, g, H, d) \\
       =&\ -\tau (g^\top d + \thalf \max\{d^\top H d,0\}) + \|c(x)\|_1.
  \ealigned
\end{equation}
Specifically, in order to ensure in iteration $k$ that $\tau_k \leq \tau_{k-1}$ and
\bequation \label{eq:qquadlb}
  \Delta q(x_k,\tau,g_k,H_k,d_k) \geq \thalf \tau \max\{d_k^\top H_k d_k, 0\} + \sigma \|c_k\|_1 \geq 0
\eequation
holds for all $\tau \leq \tau_k$, the algorithm sets, for user-defined $\sigma \in (0,1)$, the value
\bequation \label{eq:tautrial}
  \tauktri\gets
  \begin{cases}
    \infty & \mbox{if $g_k^\top d_k + \max\{d_k^\top H_k d_k, 0\} \leq 0$} \\
    \frac{(1-\sigma)\|c_k\|_1}{g_k^\top d_k + \max\{d_k^\top H_k d_k,0\}} & \mbox{otherwise,}
  \end{cases}
\eequation
and then sets, for some $\epsilon_\tau \in (0,1)$, the merit parameter value
\bequation \label{eq:taudef}
  \tau_k \gets \bcases \tau_{k-1} & \mbox{if $\tau_{k-1} \leq \tauktri$} \\ (1-\epsilon_\tau)\tauktri & \mbox{otherwise.} \ecases
\eequation
Then, for use in the stepsize computation (as motivated in \cite{ABerahas_FECurtis_DPRobinson_BZhou_2021}) it sets
\bequation \label{eq:xidef}
  \xiktri \gets \frac{\Delta q(x_k, \tau_k, g_k, H_k, d_k)}{\tau_k \|d_k\|^2}\ \ \text{then}\ \ \xi_k \gets \bcases \xi_{k-1} & \mbox{if }\xi_{k-1} \leq \xiktri \\ (1-\epsilon_\xi)\xiktri & \mbox{otherwise} \ecases
\eequation
for some $\epsilon_\xi \in (0,1)$, which, for one thing, ensures $\xi_k \leq \xiktri$.  The last component in the $k$th iteration is to set the stepsize, the magnitude of which is controlled by prescribed $\{\beta_k\} \subset (0,1]$, which is employed in the following projection interval that is used in the stepsize computation:
\bequationNN
    \mbox{Proj}_k(\cdot) := \mbox{Proj}\left(\cdot \; \Bigg| \left[\frac{\beta_k \xi_k \tau_k}{\tau_k L + \Gamma}, \frac{\beta_k \xi_k \tau_k}{\tau_k L + \Gamma} + \theta \beta_k^2 \right]\right),
\eequationNN
where $\mbox{Proj}(\cdot \; | \; \mathcal{I})$ represents the projection operator onto the interval $\mathcal{I} \subset \R{}$.  As in other stochastic-gradient-based methods, the convergence properties of the method depend on properties of $\{\beta_k\}$, which in many analyses is considered to be a constant or diminishing sequence.  In our analysis, we establish our result for the case of $\beta_k = \Ocal(1/\sqrt{\kmax+1})$ for all $k \in [\kmax]$.

Overall, the algorithm that we consider is stated as Algorithm~\ref{alg:ssqp}.  The only changes from Algorithm~3.1 in \cite{ABerahas_FECurtis_DPRobinson_BZhou_2021} are the fixed iteration limit ($\kmax$), a slightly stronger decrease requirement for the definition of $\tau_k$ in~\eqref{eq:taudef} when $\tau_{k-1} > \tauktri$, and the concluding step for producing the return value ($x_{k^*}$).  This method of sampling $k^*$ to produce the return value is consistent with other approaches in the literature on complexity analyses for algorithms for solving nonconvex optimization problems; see, e.g., \cite{DDavis_DDrusvyatskiy_2019}.  It amounts to uniform sampling over the iterates when constant $\{\beta_k\}$ is considered, as in our analysis.  Finally, we remark that Algorithm~\ref{alg:ssqp} presumes knowledge of Lipschitz constants for the objective and constraint gradients, although in practice one might only estimate these values using standard procedures~\cite{FECurtis_DPRobinson_2019}.

\begin{algorithm}[ht]
  \caption{Stochastic SQP Algorithm}
  \label{alg:ssqp}
  \begin{algorithmic}[1]
    \Require $x_0 \in \R{n}$; $\kmax \in \N{}$; $\tau_{-1} \in \R{}_{>0}$; $\epsilon_\tau \in (0,1)$; $\epsilon_\xi \in (0,1)$; $\sigma \in (0,1)$; $\xi_{-1} \in \R{}_{>0}$; $\{\beta_k\} \subset (0,1]$; $\theta \in \R{}_{\geq0}$; $L \in (0,\infty)$, a Lipschitz constant for $\nabla f$; $\Gamma \in [\sum_{i=1}^m \gamma_i,\infty)$, where $\gamma_i \in (0,\infty)$ is a Lipschitz constant for $\nabla c_i$ for all $i \in [m]$
    \For{\textbf{all} $k \in [\kmax]$}
	  \State Compute $(d_k,y_k)$ as the solution of \eqref{eq:sqpsystem}
	  \If{$d_k = 0$}
	    \State Set $\tauktri \gets \infty$, $\tau_k \gets \tau_{k-1}$, $\xiktri \gets \infty$, and $\xi_k \gets \xi_{k-1}$
	    \State Set $\widehat\alpha_{k,\text{init}} \gets 1$, $\widetilde\alpha_{k,\text{init}} \gets 1$, and $\alpha_k \leftarrow 1$
	  \Else\ (if $d_k \neq 0$)
	    \State Set $\tauktri$ by \eqref{eq:tautrial} and $\tau_k$ by \eqref{eq:taudef} 
	    \State Set $\xiktri$  and $\xi_k$ by \eqref{eq:xidef} 
        \State Set
        \bequationNN
          \widehat\alpha_{k,\text{init}} \gets \tfrac{\beta_k\Delta q(x_k,\tau_k,g_k,H_k,d_k)}{(\tau_k L + \Gamma)\|d_k\|_2^2}\ \ \text{and}\ \  \widetilde\alpha_{k,\text{init}} \gets \widehat\alpha_{k,\text{init}} - \tfrac{4\|c_k\|_1}{(\tau_k L +  \Gamma)\|d_k\|_2^2}
        \eequationNN
        \State Set $\widehat\alpha_k \gets \proj_k(\widehat\alpha_{k,\text{init}})$ and $\widetilde\alpha_k \gets \proj_k(\widetilde\alpha_{k,\text{init}})$, then \label{step.alpha_projection_stochastic} \label{step.alpha_stochastic}
        \bequationNN
          \alpha_k \gets
		    \bcases
		      \widehat\alpha_k & \text{if $\widehat\alpha_k < 1$} \\
		      1 & \text{if $\widetilde\alpha_k \leq 1 \leq \widehat\alpha_k$} \\
		      \widetilde\alpha_k & \text{if $\widetilde\alpha_k > 1$}
		    \ecases
        \eequationNN
      \EndIf
      \State Set $x_{k+1} \gets x_k + \alpha_k d_k$
    \EndFor
    \State Sample $k^* \in [\kmax]$, where $\Prob[k^* = k] = \frac{\beta_k}{\sum_{k=0}^{\kmax} \beta_k}$ for all $k \in [\kmax]$, then \Return $x_{k^*}$
  \end{algorithmic}
\end{algorithm}

\section{Complexity Analysis} \label{sec:complexity}

We begin our complexity analysis by describing the algorithm as a stochastic process (Section~\ref{subsec:process}), then formalizing the assumptions that we make about the stochastic gradient estimates (Section \ref{subsec:assumptions}).  We then state, in some cases in a slightly modified form, some key lemmas from \cite{ABerahas_FECurtis_DPRobinson_BZhou_2021} that are needed for our analysis (Section \ref{subsec:algprops}).  Our generic complexity result, which has been outlined in Section~\ref{sec:complexitypreview}, is then stated and proved (Section \ref{subsec:complexity}).  Consequences and extensions of our generic complexity result are then discussed for some special cases of distributions for the stochastic gradient estimates for which our required assumptions hold with high-probability (Section \ref{subsec:subgauss}).  Finally, we conclude this section by outlining a form of our generic complexity result that relaxes one of our minor simplifying assumptions (Section \ref{subsec:ratioparam}).

Similarly as for the convergence analysis in \cite{ABerahas_FECurtis_DPRobinson_BZhou_2021}, our complexity analysis makes use of orthogonal decompositions of the search directions computed by the algorithm; in particular, for all $k \in \N{}$, we express $d_k = u_k + v_k$, where $u_k \in \mbox{Null}(J_k)$ and $v_k \in \mbox{Range}(J_k^\top)$.  We note here that conditioned on the algorithm having reached $x_k$ at iteration $k$, the normal component $v_k$ is {\it deterministic}, depending only on the constraint value $c_k$ and the Jacobian $J_k$.

In addition to the quantities that are computed explicitly in Algorithm~\ref{alg:ssqp}, our analysis also refers to the quantities that would have been computed in each iteration $k \in \N{}$, conditioned on the event that the algorithm has reached~$x_k$ as the $k$th iterate, if the true gradient $\nabla f(x_k)$ is used in place of the stochastic gradient $g_k$.  These quantities are denoted by a ``true'' superscript.  For example, in iteration $k$, the true search direction and corresponding true Lagrange multiplier estimate are the solution of the linear system
\bequation \label{eq:sqpsystemdet}
    \left[
    \begin{matrix}
    H_k & J_k^\top \\
    J_k & 0
    \end{matrix}
    \right]
    \left[
    \begin{matrix}
    \dktrue \\ \yktrue
    \end{matrix}
    \right]
    = -\left[
    \begin{matrix}
    \nabla f(x_k) \\ c_k
    \end{matrix}
    \right],
\eequation
which may be decomposed as $\dktrue = \uktrue + v_k$, where $\uktrue \in \mbox{Null}(J_k)$ and $v_k \in \mbox{Range}(J_k^\top)$.  Here, we write $v_k$ (without a superscript) since the normal component of the search direction is defined in a manner that makes it independent of the objective gradient (estimate).  Similarly, the true value of the merit parameter that would have been computed is denoted
\bequationNN
    \tktritrue \gets
    \begin{cases}
        \infty\ \ \ \text{if $\nabla f(x_k)^\top \dktrue + \max\{(\dktrue)^\top H_k \dktrue, 0\} \leq 0$} & \\
        \frac{(1-\sigma)\|c_k\|_1}{\nabla f(x_k)^\top \dktrue + \max\{(\dktrue)^\top H_k \dktrue,0\}}\ \ \  \text{otherwise}. &
    \end{cases}
\eequationNN
This definition of $\tktritrue$ guarantees that, for any $\tau \leq \tktritrue$, one finds
\bequation \label{eq:qquadlbtrue}
    \Delta q(x_k,\tau,\nabla f(x_k),H_k,\dktrue) \geq \thalf \tau \max\{(\dktrue)^\top H_k \dktrue, 0\} + \sigma \|c_k\|_1.
\eequation

\subsection{Stochastic Process} \label{subsec:process}

Henceforth, for the sake of formality, we shall refer in our analysis to the stochastic process generated by Algorithm~\ref{alg:ssqp}.  Specifically, in terms of values that are computed by the algorithm itself, we have the stochastic process
\bequationNN
  \{(X_k,G_k,D_k,Y_k,\Tcal_k,\Xi_k,\Acal_k)\},
\eequationNN
where, for all $k \in \N{}$, the random variables are: the algorithm iterate $X_k$, stochastic gradient estimate $G_k$, search direction $D_k$, Lagrange multiplier estimate $Y_k$, merit parameter $\Tcal_k$, ratio parameter $\Xi_k$, and stepsize $\Acal_k$.  For all $k \in \N{}$, a realization of the corresponding element of this process is denoted $(x_k,g_k,d_k,y_k,\tau_k,\xi_k,\alpha_k)$.  Similarly, in terms of ``true'' values and step decomposition values that are not computed by the algorithm, but are defined for the sake of our analysis, we have the simultaneously generated process
\bequationNN
  \{(\Vk,\Uk,\Dktrue,\Uktrue,\Yktrue,\Tktritrue)\},
\eequationNN
where, for all $k \in \N{}$, the random variables are: the normal search direction component $\Vk$, the tangential search direction component $\Uk$, the true search direction $\Dktrue$, the true tangential search direction component $\Uktrue$, the true Lagrange multiplier estimate $\Yktrue$, and the true trial merit parameter $\Tktritrue$.  For all $k \in \N{}$, a realization of the corresponding element of this process is denoted $(v_k,u_k,\dktrue,\uktrue,\yktrue,\tktritrue)$.  Finally, for the sake of tracking the number of merit parameter updates that occur during runs of the algorithm, we define the stochastic process $\{S_k\}$, where for all $k \in \N{}$ the random variable $S_k$ represents the number of merit parameter decreases up to the end of the $k$th iteration, i.e., the number of iterations in which $\Tcal_k < \Tcal_{k-1}$.  For all $k \in \N{}$, a realization of $S_k$ is denoted $s_k$.

In any run, the behavior of Algorithm~\ref{alg:ssqp} is dictated entirely by the initial conditions and the sequence of stochastic gradient estimates that are generated.  Let $\Fcal_k$ denote the $\sigma$-algebra generated by the random variables $\{G_0,\dots,G_{k-1}\}$, a realization of which (along with all initial conditions of the algorithm, including $X_0 = x_0$) determines the realizations of
\bequationNN
  \{X_j\}_{j=1}^k\ \text{and}\ \{(D_j,Y_j,\Tcal_j,\Xi_j,\Acal_j,V_j,U_j,D_j^{\rm true},U_j^{\rm true},Y_j^{\rm true},\Tcal_j^{\rm trial, true},S_j)\}_{j=0}^{k-1}.
\eequationNN
For completeness, let $\Fcal_{0} = \sigma(x_0)$.  As a result, $\{\Fcal_k\}_{k\geq 0}$ is a filtration.  When conditioning on a specific realization of Algorithm~\ref{alg:ssqp} up to the beginning of iteration $k \in \N{}$, we condition on $G_{[k-1]} = g_{[k-1]}$, since these stochastic gradients determine $x_k$.  (Recall our notation that $g_{[k-1]}$ represents $\{g_0,\dots,g_{k-1}\}$.)  Similarly, later in our analysis when we condition on $\Fcal_{k}$, we are conditioning on all realizations of $G_{[k-1]}$ that are measurable with respect to the filtration~$\Fcal_{k}$.

\subsection{Assumptions} \label{subsec:assumptions}

Our analysis presumes certain good behavior of the sequences of merit and ratio parameters that are set adaptively by the algorithm.  Formally, given $(\kmax,\smax,\tau_{\min},\xi_{\min}) \in \N{} \times \N{} \times \R{}_{>0} \times \R{}_{>0}$, our main result characterizes the worst-case behavior of Algorithm~\ref{alg:ssqp} conditioned on the event denoted as
\bequation\label{def:E}
  E := E(\kmax,\smax,\tau_{\min},\xi_{\min}),
\eequation
which we define as the event such that, in every realization of the algorithm, the merit parameters $\{\tau_k\}_{k=0}^{\kmax}$ and ratio parameters $\{\xi_k\}_{k=0}^{\kmax}$ satisfy
\bitemize
  \item $\tau_k \geq \tau_{\min} > 0$ for all $k \in [\kmax]$,
  \item $\tktritrue \geq \tau_{\min} > 0$ for all $k \in [\kmax]$,
  \item $\xi_k = \xi_{\min} > 0$ for all $k \in [\kmax]$, and
  \item $|\{k \in [\kmax] : \tau_k < \tau_{k-1}\}| \leq \smax$.
\eitemize
Consideration of this event as a focus for proving a worst-case complexity result for Algorithm~\ref{alg:ssqp} is justifiable for the following reasons.

\bitemize
  \item The condition in $E$ that the ratio parameter sequence is constant over all iterations is not actually essential for our analysis; rather, it is made for the sake of simplicity.  Indeed, in Section \ref{subsec:ratioparam}, we present an extension of our main result to the setting in which this parameter sequence is not constant.  Observe that, as proved in \cite[Lemma~3.5]{ABerahas_FECurtis_DPRobinson_BZhou_2021}, the sequence $\{\Xi_k\}$ is bounded below by a positive real number whose value is deterministic, i.e., it is independent of the sequence of stochastic gradient estimates that are generated by the algorithm.  Hence, for the sake of simplicity, we assume for now that $\{\Xi_k\}$ is constant and leave the statement of the more complicated version of our main result to a subsection at the end of our analysis.
  \item The conditions in $E$ pertaining to the behavior of the merit parameter sequence are not necessarily minor.  That said, attention to this behavior of the algorithm is justified by arguments made in \cite{BeraCurtONeiRobi21,ABerahas_FECurtis_DPRobinson_BZhou_2021}, which under the same kinds of assumptions made in this paper argue that, in any run of the algorithm, the probability is zero that the merit parameter vanishes.  Furthermore, in Section~\ref{subsec:subgauss}, we consider a particular setting in which the distributions of the stochastic gradient estimates are sub-Gaussian over any run of the algorithm, in which case we show that the merit parameter remains bounded below with high probability, meaning that our main worst-case complexity bound---which holds with high probability due to the adaptivity of the merit parameter sequence---remains essentially unchanged in this setting when we do not presume upfront that the merit parameter sequence remains bounded above a positive real number.
  \item The condition in $E$ pertaining to the existence of $\smax$ is not actually an additional requirement beyond the existence of $\tau_{\min}$ in the event.  After all, by the construction of Algorithm~\ref{alg:ssqp}, it follows that when the merit parameter is decreased, it is decreased by at least a constant factor, from which it follows (under the existence of $\tau_{\min}$) that $\smax$ exists and satisfies
  \bequation \label{eq:smaxbound}
    \smax \leq \min\left\{\kmax+1, \left\lceil\frac{\log(\tau_{\min}/\tau_{-1})}{\log(1-\epsilon_\tau)}\right\rceil\right\}.
  \eequation
  That said, for simplicity and generality in our analysis, we define $\smax$ as a quantity that is decoupled from the above (conservative) inequality.
\eitemize

Conditioned on $E$, we assume the following about the stochastic gradient estimates.  Such an assumption, namely, that conditioned on the event that a given iterate has been reached the stochastic gradient is unbiased and has bounded variance, is common in analyses of stochastic optimization methods.  Here and throughout the remainder of the paper, we let $\Prob_k[\cdot]$ (respectively, $\E_k[\cdot]$) denote probabilty (respectively, expectation) conditioned on event $E$ and that $G_{[k-1]} = g_{[k-1]}$ for a given $g_{[k-1]}$, i.e., we define
\bequationNN
  \P_k[\cdot] := \P[\cdot | E, G_{[k-1]} = g_{[k-1]}]\ \ \text{and}\ \ \E_k[\cdot] := \E[\cdot | E, G_{[k-1]} = g_{[k-1]}].
\eequationNN

\bassumption\label{assum:eventE}
  There exists $M \in \R{}_{>0}$ such that, for all $k \in [\kmax]$ and any realization $g_{[k-1]}$ of $G_{[k-1]}$, one finds that
  \bequation\label{eq.sg}
    \E_k[G_k] = \nabla f(x_k)\ \ \text{and}\ \ \E_k[\|G_k - \nabla f(x_k)\|_2^2] \leq M.
  \eequation
  In addition, there exists $M_\tau \in \R{}_{>0}$ such that, for all $k \in [\kmax]$ and any realization $g_{[k-1]}$ of $G_{[k-1]}$, one finds that
  \bequation\label{eq.sg+}
    \baligned
      \text{either}\ \ \P_k[\nabla f(x_k)^\top(D_k - \dktrue) < 0, \Tcal_k < \tau_{k-1}] &= 0 \\
      \text{or}\ \ \E_k[\|G_k - \nabla f(x_k)\|_2 | \nabla f(x_k)^\top(D_k - \dktrue) < 0, \Tcal_k < \tau_{k-1}] &\leq M_\tau.
    \ealigned
  \eequation
\eassumption

Observe that the inequality in \eqref{eq.sg} can imply \eqref{eq.sg+}, such as when there exists $p \in (0,1]$ such that, for all $k \in [\kmax]$ and $g_{[k-1]}$, one finds that
\bequationNN
  \P_k[\nabla f(x_k)^\top(D_k - \dktrue) < 0, \Tcal_k < \tau_{k-1}] \geq p
\eequationNN
whenever this probability is nonzero.  This occurs, for example, when the objective of \eqref{eq:fdef} is a finite sum of $N$ terms and each stochastic gradient estimate is computed as a so-called mini-batch estimate through the uniform (random) selection of $b$ indices, in which case the above holds with $p = b/N$.

We make one additional assumption for our analysis, namely, the following.

\begin{assumption} \label{assum:ptau}
  There exists $p_{\tau} \in (0,1]$ such that, for all $k \in [\kmax]$ and any realization $g_{[k-1]}$ of $G_{[k-1]}$, one finds that
  \bequationNN
    \P_k[G_k^\top D_k + \max\{D_k^\top H_k D_k,0\} \geq \nabla f(x_k)^\top \dktrue + \max\{(\dktrue)^\top H_k \dktrue,0\}] \geq p_{\tau}.
  \eequationNN
\end{assumption}

Similar to \cite[Proposition~3.16]{ABerahas_FECurtis_DPRobinson_BZhou_2021}, Assumption~\ref{assum:ptau} allows us to prove that, with high probability, the number of iterations in which $\tau_k > \tktritrue$ is not too large.  In \cite[Example 3.17]{ABerahas_FECurtis_DPRobinson_BZhou_2021}, it was shown that the inequality in this assumption holds with $p_{\tau} = \frac12$ when, conditioned on having reached $x_k$, the stochastic gradient $G_k$ has a Gaussian distribution.  We show in Section~\ref{subsec:subgauss} that this result can be extended to other settings as well.

\subsection{Properties of Algorithm \ref{alg:ssqp}} \label{subsec:algprops}

In this section, we state key results from \cite{ABerahas_FECurtis_DPRobinson_BZhou_2021} that are needed for our analysis.

By \cite[Lemma 2.10]{ABerahas_FECurtis_DPRobinson_BZhou_2021}, there exists $\kappa_{uv} \in \R{}_{>0}$ such that, for all $k \in [\kmax]$, if $\|\uktrue\|^2 \geq \kappa_{uv} \|v_k\|^2$, then $\frac12 (\dktrue)^\top H_k \dktrue \geq \frac14 \zeta \|\uktrue\|^2$, where $\zeta$ is defined in Assumption \ref{assum:Hk}.  Correspondingly, let us define
\bequationNN
  \Psi_k := \begin{cases}
        \|\uktrue\|^2 + \|c_k\| & \mbox{if } \|\uktrue\|^2 \geq \kappa_{uv} \|v_k\|^2 \\
        \|c_k\| & \mbox{otherwise.}
    \end{cases}
\eequationNN

The following lemma is stated using a different norm (for $c_k$) than in \cite[Lemma 2.11]{ABerahas_FECurtis_DPRobinson_BZhou_2021}.  The result holds in the same manner due to the norm-equivalence between $\|\cdot\|$ and $\|\cdot\|_1$ in $\R{m}$.  We state the result in this manner for consistency in the measure of constraint violation stated in our final complexity bound. 

\begin{lemma}[{\cite[Lemma 2.11]{ABerahas_FECurtis_DPRobinson_BZhou_2021}}] \label{lem:kappapsi}
  Let Assumptions \ref{assum:fcsmooth} and \ref{assum:Hk} hold.  Then, there exists $\kappa_{\Psi} \in \R{}_{>0}$ such that, for all $k \in [\kmax]$, the true search direction and constraint violation satisfy $\|\dktrue\|^2 + \|c_k\|_1 \leq (\kappa_{\Psi} + 1) \Psi_k$.
\end{lemma}

\begin{lemma}[{\cite[Lemma 2.12]{ABerahas_FECurtis_DPRobinson_BZhou_2021}}] \label{lem:kappaq}
  Let Assumptions \ref{assum:fcsmooth} and \ref{assum:Hk} hold. Then, there exists $\kappa_q \in \R{}_{>0}$ such that, for all $k \in [\kmax]$ and any $\tau \leq \tktritrue$, the true reduction in the merit model satisfies $\Delta q (x_k,\tau,\nabla f(x_k),H_k,\dktrue) \geq \kappa_q \tau \Psi_k$.
\end{lemma}

\blemma[{\cite[Lemma~3.7]{ABerahas_FECurtis_DPRobinson_BZhou_2021}}]\label{lem:phidecrease}
  Let Assumptions \ref{assum:fcsmooth} and \ref{assum:Hk} hold and suppose that the sequence $\{\beta_k\}$ is chosen such that $\beta_k \xi_k \tau_k/(\tau_k L + \Gamma) \in (0,1]$ for all $k \in [\kmax]$.  Then, for all $k \in [\kmax]$, it follows that
  \bequation \label{eq:phidecrease}
    \baligned
      \phi(x_k + \alpha_k d_k, \tau_k) - \phi(x_k, \tau_k)
      &\leq -\alpha_k \Delta q(x_k,\tau_k,\nabla f(x_k),H_k,\dktrue) \\
      &\quad + \thalf \alpha_k \beta_k  \Delta q(x_k,\tau_k,g_k,H_k,d_k) \\
      &\quad + \alpha_k \tau_k \nabla f(x_k)^\top (d_k - \dktrue).
    \ealigned
  \eequation
\elemma

\blemma \label{lem:dexpectation}
  Let Assumptions \ref{assum:fcsmooth}, \ref{assum:Hk}, and \ref{assum:eventE} hold. Then, for all $k \in [\kmax]$, it follows that $\E_k[D_k] = \dktrue$, $\E_k[U_k] = \uktrue$, and $\E_k[Y_k] = \yktrue$.  Moreover, there exists $\kappa_d \in \R{}_{>0}$ such that, for all $k \in [\kmax]$, one finds that
  \bequationNN
    \baligned
      \|\dktrue\| \leq \kappa_d \|\nabla f(x_k)\| &\leq \kappa_d \kappa_g, \\
      \E_k[\|D_k - \dktrue\|] \leq \kappa_d \E_k[\|G_k - \nabla f(x_k)\|] &\leq \kappa_d \sqrt{M},\ \ \text{and} \\
      \E_k[\|D_k - \dktrue\| | \nabla f(x_k)^\top (D_k - \dktrue) < 0, \Tcal_k < \tau_{k-1}] &\leq \kappa_d M_\tau.
    \ealigned
  \eequationNN
\elemma
\bproof
  Except for the final inequality, the result follows directly from \cite[Lemma 3.8]{ABerahas_FECurtis_DPRobinson_BZhou_2021} (or the proof therein) and~\eqref{eq.kappa_g}.  As for the final inequality, observe by the arguments in \cite[Lemma 3.8]{ABerahas_FECurtis_DPRobinson_BZhou_2021} that for any realization of $G_k$ and $D_k$ one finds $\|d_k - \dktrue\| \leq \kappa_d \|g_k - \nabla f(x_k)\|$, which combined with~\eqref{eq.sg+} gives
  \begin{align*}
    &\ \E_k[\|D_k - \dktrue\| | \nabla f(x_k)^\top (D_k - \dktrue) < 0, \Tcal_k < \tau_{k-1}] \\
    \leq&\ \kappa_d \E_k[\|G_k - \nabla f(x_k)\| | \nabla f(x_k)^\top (D_k - \dktrue) < 0, \Tcal_k < \tau_{k-1}] \leq \kappa_d M_\tau,
  \end{align*}
  as desired. \qed
\eproof

\blemma[{\cite[Lemma~3.9]{ABerahas_FECurtis_DPRobinson_BZhou_2021}}]\label{lem:product_bounds}
  Let Assumptions \ref{assum:fcsmooth}, \ref{assum:Hk}, and \ref{assum:eventE} hold. Then, for all $k \in [\kmax]$, it follows that
  \bequationNN
    \baligned
      \nabla f(x_k)^\top \dktrue \geq \E_k[G_k^\top D_k] &\geq \nabla f(x_k)^\top \dktrue - \zeta^{-1}M \\
      \text{and}\ \ \E_k[D_k^\top H_k D_k] &\geq (\dktrue)^\top H_k \dktrue.
    \ealigned
  \eequationNN
\elemma

\subsection{Complexity Result} \label{subsec:complexity}

In this section, we present our main complexity results. We derive our results in largely the same manner as the global convergence result in \cite{ABerahas_FECurtis_DPRobinson_BZhou_2021}, but with two major changes that stem from the need to characterize the behavior of the algorithm in the context of an adaptive merit parameter sequence.  At a high level, the two modifications are as follows:
\begin{enumerate}
    \item We derive, in Lemma~\ref{lem:innerprod}, an upper bound for the last term in~\eqref{eq:phidecrease}, the derivation of which is complicated by the fact that, conditioned on $x_k$ being the $k$th iterate in a run of the algorithm, this term is the product of three correlated random variables: $\Acal_k$, $\Tcal_k$, and $\nabla f(x_k)^\top (D_k - \dktrue)$.  A critical aspect of our derived bound is that we isolate a term for the event when $\nabla f(x_k)^\top (D_k - \dktrue) < 0$ and $\Tcal_k < \tau_{k-1}$, since this happens to be an event that complicates subsequent aspects of our analysis.  In Lemma~\ref{lem:Eksumlsmax}, we prove a high-probability bound on the sum of the probabilities of the occurrences of this event over the entire run of the algorithm.
    \item A critical aspect of the analysis in \cite{ABerahas_FECurtis_DPRobinson_BZhou_2021} for the deterministic setting is that one can always tie the reduction in the model of the merit function to a first-order stationarity error measure (with respect to the constrained optimization problem) due to the fact that $\tktritrue \geq \tau_k$ for all $k \in \N{}$.  Unfortunately, however, this inequality is not guaranteed to hold in the stochastic setting, which is problematic for our purposes in this paper.  To account for this issue, we define an auxiliary sequence $\{\hat\Tcal_k\}$ (not generated by the algorithm) such that $\hat{\Tcal}_k := \min\{\Tcal_k, \Tktritrue\}$ for all $k \in [\kmax]$.  In Lemmas~\ref{lem:dqbound} and \ref{lem:hattaudiff}, we analyze behaviors of the algorithm with respect to this auxiliary sequence, and in Lemma~\ref{lem:Ktaulsmax} we provide a high-probability bound on the total number of iterations in which $\Tktritrue < \Tcal_k$ may occur.  (More precisely, Lemma~\ref{lem:Ktaulsmax} considers a superset of the iterations in which this bound may occur, which serves our purposes just as well.)
\end{enumerate}

The first few results in this section consider properties of algorithmic quantities conditioned on the algorithm having generated a certain sequence of stochastic gradient estimates through a given iteration.  In particular, given $k \in [\kmax]$ and $g_{[k-1]}$, values generated by the algorithm have been determined up to the beginning of iteration $k$, including $x_k$, $\dktrue$, and~$\tau_{k-1}$.  Given these quantities, let us define three events:
\bitemize
  \item $\Ek1$, the event that $\nabla f(x_k)^\top(D_k - \dktrue) \geq 0$;
  \item $\Ek2$, the event that $\nabla f(x_k)^\top(D_k - \dktrue) < 0$ and $\Tcal_{k} = \tau_{k-1}$; and
  \item $\Ek3$, the event that $\nabla f(x_k)^\top(D_k - \dktrue) < 0$ and $\Tcal_{k} < \tau_{k-1}$.
\eitemize

We now derive an upper bound on the final term in \eqref{eq:phidecrease}.  In the following lemma, we make use, for given $k \in [\kmax]$ and $g_{[k-1]}$, of the stepsize values
\bequation\label{eq.alpha_mins}
  \baligned
    \alpha_{\min,k}^{<} := \frac{\beta_k \xi_{\min} \tau_{\min}}{\tau_{\min}L + \Gamma},\ \ \alpha_{\min,k}^{=} &:= \frac{\beta_k \xi_{\min} \tau_{k-1}}{\tau_{k-1}L + \Gamma}, \\
    \text{and}\ \ \varsigma_{\max,k} &:= \alpha_{\min,k}^{=} + \theta \beta_k^2.
  \ealigned
\eequation
The first value here represents a lower bound on the smallest stepsize that may be computed in the event that $\Tcal_k < \tau_{k-1}$, while the second value is the smallest stepsize that may be computed in the event that $\Tcal_k = \tau_{k-1}$; it is easily verified that $\alpha_{\min,k}^{<} < \alpha_{\min,k}^{=}$.  Hence, $\varsigma_{\max,k}$ represents an upper bound on the largest stepsize that may be computed.

\begin{lemma} \label{lem:innerprod}
  Suppose that Assumptions \ref{assum:fcsmooth}, \ref{assum:Hk}, and \ref{assum:eventE} hold,  and let $\kappa_d \in \R{}_{>0}$ be defined by Lemma~\ref{lem:dexpectation}.  Then, for all $k \in [\kmax]$ and $g_{[k-1]}$, and with the stepsizes $(\alpha_{\min,k}^{<},\alpha_{\min,k}^{=},\varsigma_{\max,k})$ defined as in \eqref{eq.alpha_mins}, one finds that
  \begin{align*}
    &\ \E_k[\Acal_k \Tcal_k \nabla f(x_k)^\top(D_k - \dktrue)] \\
    \leq&\ (\varsigma_{\max,k}\tau_{k-1} - \alpha_{\min,k}^{<} \tau_{\min}) \kappa_g \kappa_d M_\tau \Prob_k[\Ek3] + \theta \beta_k^2 \tau_{k-1} \kappa_g \kappa_d \sqrt{M}.
  \end{align*}
\end{lemma}
\begin{proof}
  Consider arbitrary $k \in [\kmax]$ and $g_{[k-1]}$, and for ease of exposition, let us denote $\E_{k,j} = \E_k[\nabla f(x_k)^\top(D_k - \dktrue) | E_{k,j}]$ for all $j \in \{1,2,3\}$.  By the Law of Total Expectation, the fact that $0 < \tau_{\min} \leq \Tcal_k \leq \tau_{k-1}$ under $E$, and the definitions of $\alpha_{\min,k}^{<}$, $\alpha_{\min,k}^{=}$, and $\varsigma_{\max,k}$, one finds that
  \begin{align*}
    &\ \E_k[\Acal_k \Tcal_k \nabla f(x_k)^\top(D_k - \dktrue)] \\
    =&\ \E_k[\Acal_k \Tcal_k \nabla f(x_k)^\top(D_k - \dktrue)|\Ek1] \Prob_k[\Ek1] \\
     &\ + \E_k[\Acal_k \Tcal_k \nabla f(x_k)^\top(D_k - \dktrue)|\Ek2] \Prob_k[\Ek2] \\
     &\ + \E_k[\Acal_k \Tcal_k \nabla f(x_k)^\top(D_k - \dktrue)|\Ek3] \Prob_k[\Ek3] \\
    \leq&\ \varsigma_{\max,k} \tau_{k-1}  \E_{k,1} \Prob_k[\Ek1] + \alpha_{\min,k}^{=} \tau_{k-1} \E_{k,2} \Prob_k[\Ek2] + \alpha_{\min,k}^{<} \tau_{\min} \E_{k,3} \Prob_k[\Ek3].
  \end{align*}
  Using this inequality, the Law of Total Expectation, and Lemma \ref{lem:dexpectation} ($\E_k[D_k] = \dktrue$), one obtains three upper bounds by adding and subtracting like terms:
  \begin{align*}
    &\E_k[\Acal_k \Tcal_k \nabla f(x_k)^\top(D_k - \dktrue)] \\
    \leq&\ \varsigma_{\max,k} \tau_{k-1} \E_{k,1} \Prob_k[\Ek1] + \varsigma_{\max,k} \tau_{k-1} \E_{k,2} \Prob_k[\Ek2] - \theta \beta_k^2 \tau_{k-1} \E_{k,2} \Prob_k[\Ek2] \\
    &\ + \varsigma_{\max,k} \tau_{k-1} \E_{k,3} \Prob_k[\Ek3] + (\alpha_{\min,k}^{<} \tau_{\min} - \varsigma_{\max,k} \tau_{k-1}) \E_{k,3} \Prob_k[\Ek3] \\
    =&\ - \theta \beta_k^2 \tau_{k-1} \E_{k,2} \Prob_k[\Ek2] + (\alpha_{\min,k}^{<} \tau_{\min} - \varsigma_{\max,k} \tau_{k-1}) \E_{k,3} \Prob_k[\Ek3] \\ 
    \intertext{and}
    &\E_k[\Acal_k \Tcal_k \nabla f(x_k)^\top(D_k - \dktrue)] \\
    \leq&\ \alpha_{\min,k}^{=} \tau_{k-1} \E_{k,1} \Prob_k[\Ek1] + \theta \beta_k^2 \tau_{k-1} \E_{k,1} \Prob_k[\Ek1] + \alpha_{\min,k}^{=} \tau_{k-1} \E_{k,2} \Prob_k[\Ek2] \\
    &\ + \alpha_{\min,k}^{=} \tau_{k-1} \E_{k,3} \Prob_k[\Ek3] + (\alpha_{\min,k}^{<} \tau_{\min} - \alpha_{\min,k}^{=} \tau_{k-1}) \E_{k,3} \Prob_k[\Ek3] \\
    =&\ \theta \beta_k^2 \tau_{k-1} \E_{k,1} \Prob_k[\Ek1] + (\alpha_{\min,k}^{<} \tau_{\min} - \alpha_{\min,k}^{=} \tau_{k-1}) \E_{k,3} \Prob_k[\Ek3] \\ 
    \intertext{and}
    &\E_k[\Acal_k \Tcal_k \nabla f(x_k)^\top(D_k - \dktrue)] \\
    \leq&\ \alpha_{\min,k}^{<} \tau_{\min} \E_{k,1} \Prob_k[\Ek1] + (\varsigma_{\max,k} \tau_{k-1} - \alpha_{\min,k}^{<} \tau_{\min}) \E_{k,1} \Prob_k[\Ek1] \\
    &\ + \alpha_{\min,k}^{<} \tau_{\min} \E_{k,2} \Prob_k[\Ek2] + (\alpha_{\min,k}^{=} \tau_{k-1} - \alpha_{\min,k}^{<} \tau_{\min}) \E_{k,2} \Prob_k[\Ek2] \\
    &\ + \alpha_{\min,k}^{<} \tau_{\min} \E_{k,3} \Prob_k[\Ek3] \\
    =&\ (\varsigma_{\max,k} \tau_{k-1} - \alpha_{\min,k}^{<} \tau_{\min}) \E_{k,1} \Prob_k[\Ek1] 
    + (\alpha_{\min,k}^{=} \tau_{k-1} - \alpha_{\min,k}^{<} \tau_{\min}) \E_{k,2} \Prob_k[\Ek2].
  \end{align*}
  Averaging these three upper bounds and the definition of $\varsigma_{\max,k}$, one obtains
  \begin{align}
    &\E_k[\Acal_k \Tcal_k \nabla f(x_k)^\top(D_k - \dktrue)] \nonumber \\
    \leq&\ \tfrac13 ((\alpha_{\min,k}^{=} + 2\theta \beta_k^2) \tau_{k-1} - \alpha_{\min,k}^{<} \tau_{\min}) \E_{k,1} \Prob_k[\Ek1] \nonumber \\
    &\ + \tfrac13 ((\alpha_{\min,k}^{=} - \theta \beta_k^2) \tau_{k-1} - \alpha_{\min,k}^{<} \tau_{\min}) \E_{k,2} \Prob_k[\Ek2] \nonumber \\
    &\ + \tfrac13 (2\alpha_{\min,k}^{<} \tau_{\min} - (2\alpha_{\min,k}^{=} + \theta \beta_k^2) \tau_{k-1}) \E_{k,3} \Prob_k[\Ek3] \nonumber \\
    =&\ \tfrac13 ((\alpha_{\min,k}^{=} + 2\theta \beta_k^2) \tau_{k-1} - \alpha_{\min,k}^{<} \tau_{\min}) (\E_{k,1} \Prob_k[\Ek1] + \E_{k,2} \Prob_k[\Ek2]) \nonumber \\
    &\ - \theta \beta_k^2 \tau_{k-1} \E_{k,2} \Prob_k[\Ek2] \nonumber \\
    &\ - \tfrac13 ((2\alpha_{\min,k}^{=} + \theta \beta_k^2) \tau_{k-1} - 2\alpha_{\min,k}^{<} \tau_{\min}) \E_{k,3} \Prob_k[\Ek3]. \label{eq:innerprodproof1}
  \end{align}
  This bound can be rewritten as follows.  By the Law of Total expectation,
  \begin{align}
    &\ \E_{k,1} \Prob_k[\Ek1] + \E_{k,2} \Prob_k[\Ek2] \nonumber \\
    =&\ \E_k[\nabla f(x_k)^\top(D_k - \dktrue)] -\E_{k,3} \Prob_k[\Ek3] = -\E_{k,3} \Prob_k[\Ek3], \label{eq.new1}
  \end{align}
  and along with Lemma \ref{lem:dexpectation} and~\eqref{eq.kappa_g} one finds
  \begin{align}
    -\E_{k,2} \Prob_k [\Ek2] &= -\E_k[\nabla f(x_k)^\top(D_k - \dktrue)| \Ek2] \Prob_k [\Ek2] \nonumber \\
    &\leq \E_k[\|\nabla f(x_k)\| \|D_k - \dktrue\| | \Ek2] \Prob_k [\Ek2] \nonumber \\
    &= \E_k[\|\nabla f(x_k)\| \|D_k - \dktrue\|] \nonumber \\
    &\quad - \E_k[\|\nabla f(x_k)\| \|D_k - \dktrue\| | \Ek1]\Prob_k [\Ek1] \nonumber \\
    &\quad - \E_k[\|\nabla f(x_k)\| \|D_k - \dktrue\| | \Ek3] \Prob_k [\Ek3] \nonumber \\
    &\leq \E_k[\|\nabla f(x_k)\| \|D_k - \dktrue\|] \nonumber \\
    &\leq \kappa_g \E_k[\|D_k - \dktrue\|] \leq \kappa_g \kappa_d \sqrt{M}. \label{eq.new2}
  \end{align}
  In addition, Lemma \ref{lem:dexpectation} and~\eqref{eq.kappa_g} also yield that
  \begin{align}
    -\E_{k,3} \Prob_k[\Ek3]
    &= -\E_k[\nabla f(x_k)^\top(D_k - \dktrue)| \Ek3] \Prob_k[\Ek3] \nonumber \\
    &\leq \E_k[\|\nabla f(x_k)\| \|D_k - \dktrue\|| \Ek3] \Prob_k[\Ek3] \nonumber \\
    &\leq \kappa_g \kappa_d M_{\tau} \Prob_k[\Ek3]. \label{eq.new3}
  \end{align}
  Combining \eqref{eq:innerprodproof1}, \eqref{eq.new1}, \eqref{eq.new2}, and \eqref{eq.new3}, the desired result follows. \qed
\end{proof}

Now, conditioned on given $k \in [\kmax]$ and $g_{[k-1]}$, let us define
\bequation \label{eq:hattaudef}
  \hat{\Tcal}_k := \min\{\Tcal_k, \tktritrue\}.
\eequation

\begin{lemma} \label{lem:dqbound}
  Suppose that Assumptions \ref{assum:fcsmooth}, \ref{assum:Hk}, and \ref{assum:eventE} hold and let $\kappa_d \in \R{}_{>0}$ be defined by Lemma~\ref{lem:dexpectation}.  Then, for all $k \in [\kmax]$ and $g_{[k-1]}$, one finds that
  \begin{align*}
    &\ \E_k[\Delta q(x_k, \hat{\Tcal}_k, G_k, H_k, D_k)] \\
    \leq&\ \E_k[\Delta q(x_k, \hat{\Tcal}_k, \nabla f(x_k), H_k, \dktrue)] + \thalf (\tau_{k-1} + \tau_{\min}) \zeta^{-1} M \\
    &\ + (\tau_{k-1} - \tau_{\min}) (\kappa_d \sqrt{M} (2\kappa_g + \sqrt{M}) + \kappa_H \kappa_d^2 (M + \tfrac32 \kappa_g^2)).
  \end{align*}
\end{lemma}
\begin{proof}
  By the definition of $\Delta q$ in \eqref{eq:deltadef}, one has that
  \begin{align}
    &\ \E_k[\Delta q(x_k, \hat{\Tcal}_k, G_k, H_k, D_k)] \nonumber \\
    =&\ \E_k[-\hat{\Tcal}_k(G_k^\top D_k + \thalf \max\{D_k^\top H_k D_k, 0\}) + \|c_k\|_1] \nonumber \\
    =&\ \E_k[-\hat{\Tcal}_k(G_k^\top D_k - \nabla f(x_k)^\top \dktrue  + \thalf \max\{D_k^\top H_k D_k, 0\} \nonumber \\
    &\qquad - \thalf \max\{(\dktrue)^\top H_k \dktrue, 0\})] \nonumber \\
    &\ + \E_k[-\hat{\Tcal}_k(\nabla f(x_k)^\top \dktrue + \thalf \max\{(\dktrue)^\top H_k \dktrue, 0\}) + \|c_k\|_1] \nonumber \\
    =&\ \E_k[\hat{\Tcal}_k(\nabla f(x_k)^\top \dktrue - G_k^\top D_k + \thalf \max\{(\dktrue)^\top H_k \dktrue, 0\} \nonumber \\
    &\ - \thalf \max\{D_k^\top H_k D_k, 0\})] + \E_k[\Delta q(x_k, \hat{\Tcal}_k, \nabla f(x_k), H_k, \dktrue)]. \label{eq:dqbound1}
  \end{align}
  Now, for simplicity of notation, define
  \begin{align*}
    Q_k :=&\ \nabla f(x_k)^\top \dktrue - G_k^\top D_k \\
    &\ + \thalf \max\{(\dktrue)^\top H_k \dktrue, 0\} - \thalf \max\{D_k^\top H_k D_k, 0\}.
  \end{align*}
  Let $E_Q$ denote the event that $Q_k \geq 0$ occurs and let $E_Q^c$ denote the event that $Q_k < 0$ occurs.  By the Law of Total Expectation, one has that
  \begin{align*}
    \E_k[\hat{\Tcal}_k Q_k]
      &= \E_k[\hat{\Tcal}_k Q_k | E_Q] \Prob_k[E_Q] + \E_k[\hat{\Tcal}_k Q_k | E_Q^c] \Prob_k[E_Q^c] \\
      &\leq \tau_{k-1} \E_k[Q_k | E_Q] \Prob_k[E_Q] + \tau_{\min} \E_k[Q_k | E_Q^c] \Prob_k[E_Q^c].
  \end{align*}
  Therefore, by the Law of Total Probability, Lemma \ref{lem:product_bounds}, Jensen's inequality, and convexity of $\max\{\cdot,0\}$, it follows that
  \begin{align*}
    \E_k[\hat{\Tcal}_k Q_k]
    \leq&\ \tau_{k-1} \E_k[Q_k | E_Q] \Prob_k[E_Q] + \tau_{k-1} \E_k[Q_k | E_Q^c] \Prob_k[E_Q^c] \\
    &\ + (\tau_{\min} - \tau_{k-1}) \E_k[Q_k | E_Q^c] \Prob_k[E_Q^c] \\
    =&\ \tau_{k-1} \E_k[Q_k] + (\tau_{\min} - \tau_{k-1}) \E_k[Q_k | E_Q^c] \Prob_k[E_Q^c] \\
    =&\ \tau_{k-1} (\nabla f(x_k)^\top \dktrue - \E_k[G_k^\top D_k] \\
    &\ + \thalf \max\{(\dktrue)^\top H_k \dktrue, 0\} - \thalf \E_k[\max\{D_k^\top H_k D_k, 0\}]) \\
    &\ + (\tau_{\min} - \tau_{k-1}) \E_k[Q_k | E_Q^c] \Prob_k[E_Q^c] \\
    \leq&\ \tau_{k-1} \zeta^{-1} M + (\tau_{\min} - \tau_{k-1}) \E_k[Q_k | E_Q^c] \Prob_k[E_Q^c],
  \end{align*}
  and by similar reasoning one finds that
  \begin{align*}
    \E_k[\hat{\Tcal}_k Q_k]
    \leq&\ \tau_{\min} \E_k[Q_k | E_Q] \Prob_k[E_Q] + \tau_{\min} \E_k[Q_k | E_Q^c] \Prob_k[E_Q^c] \\
    &\ + (\tau_{k-1} - \tau_{\min}) \E_k[Q_k | E_Q] \Prob_k[E_Q] \\
    =&\ \tau_{\min} \E_k[Q_k] + (\tau_{k-1} - \tau_{\min}) \E_k[Q_k | E_Q] \Prob_k[E_Q] \\
    \leq&\ \tau_{\min} \zeta^{-1} M + (\tau_{k-1} - \tau_{\min}) \E_k[Q_k | E_Q] \Prob_k[E_Q].
  \end{align*}
  Averaging these two upper bounds, one finds that
  \begin{align}
    \E_k[\hat{\Tcal}_k Q_k] 
    \leq&\ \thalf (\tau_{k-1} + \tau_{\min}) \zeta^{-1} M + \thalf (\tau_{k-1} - \tau_{\min}) \E_k[Q_k | E_Q] \Prob_k[E_Q] \nonumber \\
    &\ + \thalf (\tau_{\min} - \tau_{k-1}) \E_k[Q_k | E_Q^c] \Prob_k[E_Q^c] \label{eq:dqbound2}.
  \end{align}
  Our goal now is to bound the latter two terms in \eqref{eq:dqbound2}.  Toward this end, observe that by the triangle and Cauchy-Schwarz inequalities, the proof of Lemma~\ref{lem:dexpectation}, Assumption~\ref{assum:eventE}, Jensen's inequality, and concavity of the square root over $\R{}_{\geq0}$, one finds that
  \begin{align}
    &\ \E_k[|\nabla f(x_k)^\top \dktrue - G_k^\top D_k|] \nonumber \\
    \leq&\ \E_k[|\nabla f(x_k)^\top \dktrue - G_k^\top \dktrue|] +\E_k[|G_k^\top \dktrue - G_k^\top D_k|] \nonumber \\
    \leq&\ \|\dktrue\| \E_k [\|\nabla f(x_k) - G_k\|] + \E_k[\|G_k\|\|\dktrue - D_k\|] \nonumber \\
    \leq&\ \kappa_d \kappa_g \sqrt{\E_k[\|\nabla f(x_k)-G_k\|^2]} + \kappa_d \E_k[\|G_k\| \|\nabla f(x_k) - G_k\|]  \nonumber \\
    \leq&\ \kappa_d \kappa_g \sqrt{\E_k[\|\nabla f(x_k)-G_k\|^2]}  \nonumber \\
    &\ + \kappa_d \E_k[(\|G_k - \nabla f(x_k)\| + \|\nabla f(x_k)\|) \|\nabla f(x_k) - G_k\|] \nonumber \\
    \leq&\ \kappa_d \kappa_g \sqrt{M} + \kappa_d (M + \kappa_g \sqrt{M}) = \kappa_d \sqrt{M}(2\kappa_g + \sqrt{M}). \label{eq.wednesday1}
  \end{align}
  In addition, by the Cauchy-Schwarz inequality, the proof of Lemma~\ref{lem:dexpectation}, and Assumption~\ref{assum:eventE}, and the fact that $\|a\|^2 \leq 2(\|a - b\|^2 + \|b\|^2)$ for any $(a,b) \in \R{n} \times \R{n}$, one finds
  \begin{align}
    &\ \E_k[|\thalf \max\{(\dktrue)^\top H_k \dktrue, 0\} - \thalf \max\{D_k^\top H_k D_k, 0\}|] \nonumber \\
    \leq&\ |\thalf \max\{(\dktrue)^\top H_k \dktrue, 0\}| + \E_k[|\thalf \max\{D_k^\top H_k D_k, 0\}|] \nonumber \\
    \leq&\ \thalf \|H_k\|\|\dktrue\|^2 + \thalf \|H_k\|\E_k[\|D_k\|^2] \nonumber \\
    \leq&\ \thalf \|H_k\|\|\dktrue\|^2 + \thalf \kappa_d^2 \|H_k\|\E_k[\|G_k\|^2] \nonumber \\
    \leq&\ \thalf \|H_k\|\|\dktrue\|^2 + \kappa_d^2 \|H_k\|\E_k[\|G_k - \nabla f(x_k)\|^2 + \|\nabla f(x_k)\|^2] \nonumber \\
    \leq&\ \thalf \kappa_H \kappa_d^2 \kappa_g^2 + \kappa_H \kappa_d^2 (M + \kappa_g^2) \leq \kappa_H \kappa_d^2 (M + \tfrac32 \kappa_g^2). \label{eq.wednesday2}
  \end{align}
  By the Law of Total Expectation, \eqref{eq.wednesday1}, and \eqref{eq.wednesday2}, it follows that
  \begin{align*}
    \E_k[Q_k | E_Q] \Prob_k[E_Q] &= \E_k[|Q_k| | E_Q] \Prob_k[E_Q] \\
    &= \E_k[|Q_k|] - \E_k[|Q_k| | E_Q^c] \Prob_k[E_Q^c] \\
    &\leq \kappa_d \sqrt{M} (2\kappa_g + \sqrt{M}) + \kappa_H \kappa_d^2 (M + \tfrac32 \kappa_g^2),
  \end{align*}
  and by a similar argument, one finds that
  \begin{align*}
    -\E_k[Q_k | E_Q^c] \Prob_k[E_Q^c] &= \E_k[|Q_k| | E_Q^c] \Prob_k[E_Q^c] \\
    &= \E_k[|Q_k|] - \E_k[|Q_k| | E_Q] \Prob_k[E_Q] \\
    &\leq \kappa_d \sqrt{M} (2\kappa_g + \sqrt{M}) + \kappa_H \kappa_d^2 (M + \tfrac32 \kappa_g^2).
  \end{align*}
  The conclusion follows by combining these equations, \eqref{eq:dqbound1}, and \eqref{eq:dqbound2}. \qed
\end{proof}

Our next lemma bounds differences between expected reductions in the model of the merit function that account for cases when $\hat\Tcal_k < \Tcal_k$.

\begin{lemma} \label{lem:hattaudiff}
  Let Assumptions \ref{assum:fcsmooth}, \ref{assum:Hk}, and \ref{assum:eventE} hold and let $\kappa_d \in \R{}_{>0}$ be defined by Lemma~\ref{lem:dexpectation}.  Then, for all $k \in [\kmax]$ and $g_{[k-1]}$, with $(\alpha_{\min,k}^{<},\alpha_{\min,k}^{=},\varsigma_{\max,k})$ defined as in \eqref{eq.alpha_mins} and $\hat{\Tcal}_k$ defined in \eqref{eq:hattaudef}, one finds that
  \begin{align*}
    &\ \E_k[\Acal_k \Delta q(x_k,\hat{\Tcal}_k,\nabla f(x_k),H_k,\dktrue)] - \E_k[\Acal_k \Delta q(x_k,\Tcal_k,\nabla f(x_k),H_k,\dktrue)] \\
    \leq&\ (\varsigma_{\max,k} \tau_{k-1} - \alpha_{\min,k}^{<} \tau_{\min}) \kappa_d \kappa_g^2(1 + \thalf \kappa_H \kappa_d)
  \end{align*}
  and
  \begin{align*}
    &\ \E_k[\Acal_k \Delta q(x_k, \Tcal_k,G_k,H_k,D_k)] - \E_k[\Acal_k \Delta q(x_k,\hat{\Tcal}_k,G_k,H_k,D_k)] \\
    \leq&\ (\varsigma_{\max,k} \tau_{k-1} - \alpha_{\min,k}^{<} \tau_{\min}) \kappa_d(2 + \kappa_H \kappa_d) (M + \kappa_g^2).
  \end{align*}
\end{lemma}
\begin{proof}
  Under the stated assumptions, it follows from the stated lemma and definitions, along with the definition of $\Delta q$ in \eqref{eq:deltadef}, that
  \begin{align*}
    &\ \E_k[\Acal_k \Delta q(x_k,\hat{\Tcal}_k,\nabla f(x_k),H_k,\dktrue)] - \E_k[\Acal_k \Delta q(x_k,\Tcal_k,\nabla f(x_k),H_k,\dktrue)] \\
    =&\ \E_k[\Acal_k (\Tcal_k - \hat{\Tcal}_k) (\nabla f(x_k)^\top \dktrue + \thalf \max\{(\dktrue)^\top H_k \dktrue,0\})] \\
    \leq&\ (\varsigma_{\max,k} \tau_{k-1} - \alpha_{\min,k}^{<} \tau_{\min}) |\nabla f(x_k)^\top \dktrue + \thalf \max\{(\dktrue)^\top H_k \dktrue,0\}| \\
    \leq&\ (\varsigma_{\max,k} \tau_{k-1} - \alpha_{\min,k}^{<} \tau_{\min}) (\kappa_d \kappa_g^2 + \thalf\kappa_H \kappa_d^2 \kappa_g^2),
  \end{align*}
  and, along with $\|a\|^2 \leq 2(\|a - b\|^2 + \|b\|^2)$ for any $(a,b) \in \R{n} \times \R{n}$, one finds
  \begin{align*}
    &\ \E_k[\Acal_k \Delta q(x_k,\Tcal_k,G_k,H_k,D_k)] - \E_k[\Acal_k \Delta q(x_k,\hat{\Tcal}_k,G_k,H_k,D_k)] \\
    =&\ \E_k[\Acal_k (\hat{\Tcal}_k - \Tcal_k) (G_k^\top D_k + \thalf \max\{D_k^\top H_k D_k,0\})] \\
    \leq&\ (\varsigma_{\max,k} \tau_{k-1} - \alpha_{\min,k}^{<} \tau_{\min}) \E_k[|G_k^\top D_k + \thalf \max\{D_k^\top H_k D_k,0\}|] \\
    \leq&\ (\varsigma_{\max,k} \tau_{k-1} - \alpha_{\min,k}^{<} \tau_{\min}) (\kappa_d + \thalf \kappa_d^2 \kappa_H) \E_k[\|G_k\|^2] \\
    \leq&\ (\varsigma_{\max,k} \tau_{k-1} - \alpha_{\min,k}^{<} \tau_{\min}) (\kappa_d + \thalf \kappa_d^2 \kappa_H) (2(M + \kappa_g^2)),
  \end{align*}
  which together are the desired conclusions. \qed
\end{proof}

Our next two lemmas are critical elements of our analysis.  For both lemmas, we define, for any $s \in \N{}$ and $\delta \in (0,1)$, the quantities
\bequation \label{eq:hatdelta}
    \hdelta := \frac{\delta}{\sum_{j=0}^{\max\{\smax-1,0\}} \binom{\kmax}{j}}
\eequation
and
\bequation \label{eq:lsd}
  \ell(s,\hdelta) := s + \log(1/\hdelta) + \sqrt{\log(1/\hdelta)^2 + 2 s \log(1/\hdelta)}.
\eequation
The first of these lemmas bounds, with high probability, the sum of the probabilities of the occurrences of event $E_{k,3}$ over the run of the algorithm.

\begin{lemma} \label{lem:Eksumlsmax}
  Suppose Assumptions \ref{assum:fcsmooth}, \ref{assum:Hk}, and \ref{assum:eventE} hold. Then, for any $\delta \in (0,1)$,
  \bequation \label{eq:Ek3lsmax}
    \Prob\left[\sum_{k=0}^{\kmax} \Prob [E_{k,3} | E, \mathcal{F}_k] \leq \ell(\smax,\hdelta) + 1\Bigg | E \right] \geq 1-\delta.
  \eequation
\end{lemma}
\begin{proof}
  This result is proved in Appendix~\ref{app:probbounds}. \qed
\end{proof}

For our next lemma, let us define the random index set
\bequation \label{eq:Ktau}
  \Kcal_{\tau} := \{k \in [\kmax] : \Tktritrue < \Tcal_{k-1}\}.
\eequation
By the manner in which $\{\Tcal_k\}$, $\{\Tktritrue\}$, and $\{\hat{\Tcal}_k\}$ are defined, this set is always a superset of the iterations in which $\hat{\Tcal}_k < \Tcal_k$; hence, by bounding the cardinality of \eqref{eq:Ktau}, one bounds the cardinality of the set of iterations in which $\hat{\Tcal}_k < \Tcal_k$, which is needed for our main theorem.  The reason that we consider the set $\Kcal_\tau$ in \eqref{eq:Ktau} is the fact that, conditioned on the algorithm having generated a particular set of stochastic gradients up to the beginning of iteration $k$, whether or not $\tktritrue < \tau_{k-1}$ holds has already been determined; in other words, the occurrence of this inequality does not depend on $G_k$.  This means that  $G_k$ is conditionally independent of the event $\Tktritrue < \Tcal_{k-1}$ given $g_{[k-1]}$ and $E$, which is a fact that is exploited in the proof of the lemma.

\begin{lemma} \label{lem:Ktaulsmax}
  Suppose Assumptions \ref{assum:fcsmooth}, \ref{assum:Hk}, \ref{assum:eventE}, and \ref{assum:ptau} hold and let $\Kcal_\tau$ be defined as in \eqref{eq:Ktau}.  Then, for any $\delta \in (0,1)$, it follows that
  \bequation \label{eq:Ktaulsmax}
    \Prob\left[|\mathcal{K}_{\tau}| \leq \left\lceil \frac{\ell(\smax,\hdelta)+1}{p_{\tau}}\right\rceil \Bigg| E\right] \geq 1-\delta.
  \eequation
\end{lemma}
\begin{proof}
  This result is proved in Appendix~\ref{app:probbounds}.
\end{proof}

We are now prepared to prove a convergence rate result.

\begin{theorem} \label{thm:complexity}
  Suppose Assumptions \ref{assum:fcsmooth}, \ref{assum:Hk}, \ref{assum:eventE}, and \ref{assum:ptau} hold, let $\smax \geq 1$, let $\kappa_d \in \R{}_{>0}$ be defined by Lemma~\ref{lem:dexpectation}, define
  \bequationNN
    A_{\min} :=\ \frac{\xi_{\min} \tau_{\min}}{\tau_{\min}L + \Gamma}\ \ \text{and}\ \ A_{\max} := \frac{\xi_{-1} \tau_{-1}}{\tau_{-1}L + \Gamma},
  \eequationNN
  suppose that $\beta_k = \beta$ for all $k \in [\kmax]$ where
  \bequation\label{def.beta}
    \beta := \frac{\gamma}{\sqrt{\kmax+1}}\ \ \text{for some}\ \ \gamma \in \(0,\frac{A_{\min}}{A_{\max} + \theta}\right],
  \eequation
  define
  \bequationNN
    \baligned
      \Mbar :=&\ \tfrac14 (A_{\max} + \theta\beta) (\tau_{-1} + \tau_{\min}) \zeta^{-1} M \\
      &\ + \thalf (A_{\max} + \theta\beta) (\tau_{-1} - \tau_{\min}) (\kappa_d \sqrt{M} (2\kappa_g + \sqrt{M}) + \kappa_H \kappa_d^2 (M + \tfrac32 \kappa_g^2)) \\
      &\ + \theta \tau_{-1} \kappa_g \kappa_d \sqrt{M} \\
      \kappa_{E_3} :=&\ ((A_{\max} + \theta \beta)\tau_{-1} - A_{\min} \tau_{\min}) \kappa_g \kappa_d M_\tau, \\
      \kappa_{\Delta q,1} :=&\ ((A_{\max} + \theta \beta) \tau_{-1} - A_{\min} \tau_{\min}) \kappa_d \kappa_g^2(1 + \thalf \kappa_H \kappa_d)\ \ \text{and} \\
      \kappa_{\Delta q,2} :=&\ ((A_{\max} + \theta \beta) \tau_{-1} - A_{\min} \tau_{\min}) \kappa_d (1 + \thalf \kappa_H \kappa_d) (M + \kappa_g^2),
    \ealigned
  \eequationNN
  and, for all $k \in [\kmax]$ and $g_{[k-1]}$, let $\hat{\Tcal}_k$ be defined as in \eqref{eq:hattaudef}.  Then, for any $\delta \in (0,1)$, it follows with $K^*$ having a discrete uniform distribution over $[\kmax]$ and $\hat\delta$ and $\ell$ defined as in \eqref{eq:hatdelta} and \eqref{eq:lsd} that, with probability at least $1-\delta$,
  \begin{align}
    &\ \E[\Delta q(X_{K^*},\hat{\Tcal}_{K^*},\nabla f(X_{K^*}),H_{K^*},D_{K^*}^{\rm true}) | E] \nonumber \\
    \leq &\ 2\(\frac{\tau_{-1} (f_0 - f_{\min}) + \|c_0\|_1 + \Mbar \gamma^2 + \kappa_{E_3} \gamma (\ell(\smax,\hdelta/2)+1)/\sqrt{\kmax+1}}{A_{\min} \gamma \sqrt{\kmax+1}}\) \nonumber \\
    &\ + \frac{2(\kappa_{\Delta q,1} \gamma + \kappa_{\Delta q,2} \gamma^2/\sqrt{\kmax+1})}{A_{\min} \gamma (\kmax+1)} \left\lceil\frac{\ell(\smax,\hdelta/2)+1}{p_{\tau}}\right\rceil. \label{eq:complexity}
  \end{align}
\end{theorem}
\begin{proof}
  First, consider arbitrary $k \in [\kmax]$.  By Lemmas~\ref{lem:phidecrease} and \ref{lem:innerprod}, one has that
  \begin{align}
    &\ \E_k [\phi(x_k + \Acal_k D_k, \Tcal_k)] - \E_k [\phi(x_k, \Tcal_k)] \nonumber \\
    \leq&\ \E_k[-\Acal_k \Delta q(x_k,\Tcal_k,\nabla f(x_k),H_k,\dktrue)
    + \thalf \Acal_k \beta  \Delta q(x_k,\Tcal_k,G_k,H_k,D_k) \nonumber \\
    &\ + \Acal_k \Tcal_k \nabla f(x_k)^\top (D_k - \dktrue)] \nonumber \\
    \leq&\ \E_k[-\Acal_k \Delta q(x_k,\Tcal_k,\nabla f(x_k),H_k,\dktrue)
    + \thalf \Acal_k \beta  \Delta q(x_k,\Tcal_k,G_k,H_k,D_k)] \nonumber \\
    &\ + (\varsigma_{\max,k}\tau_{k-1} - \alpha_{\min,k}^{<} \tau_{\min}) \kappa_g \kappa_d M_\tau \Prob_k[\Ek3] + \theta \beta^2 \tau_{k-1} \kappa_g \kappa_d \sqrt{M}. \label{eq.danny_boy}
  \end{align}
  Our next aim is to prove that, roughly speaking, one in fact finds that
  \begin{multline}\label{eq.danny_elton_john}
    \E_k [\phi(x_k + \Acal_k D_k, \Tcal_k)] - \E_k [\phi(x_k, \Tcal_k)] \\
    \leq -\thalf A_{\min} \beta \E_k[\Delta q(x_k,\hat{\Tcal}_k,\nabla f(x_k),H_k,\dktrue)] + \text{``noise.''}
  \end{multline}
  Such a bound does not follow directly from \eqref{eq.danny_boy} since the first term on the right-hand side in \eqref{eq.danny_boy} involves a model reduction with respect to $\Tcal_k$ (which cannot be tied to a stationarity measure), whereas the first term on the right-hand side of the bound in \eqref{eq.danny_elton_john} involves a model reduction with respect to $\hat\Tcal_k$ (which can be tied to a stationarity measure).  Toward the aim of proving a bound of the form in \eqref{eq.danny_elton_john}, first observe that it follows with Lemma~\ref{lem:dqbound} that
  \begin{align}
    &\ \E_k[-\Acal_k \Delta q(x_k,\hat{\Tcal}_k,\nabla f(x_k),H_k,\dktrue) + \thalf \Acal_k \beta  \Delta q(x_k,\hat{\Tcal}_k,G_k,H_k,D_k)] \nonumber \\
    &\ + (\varsigma_{\max,k} \tau_{k-1} - \alpha_{\min,k}^{<} \tau_{\min}) \kappa_g \kappa_d M_\tau \Prob_k[\Ek3] + \theta \beta^2 \tau_{k-1} \kappa_g \kappa_d \sqrt{M} \nonumber \\
    \leq&\ -A_{\min} \beta \E_k[\Delta q(x_k,\hat{\Tcal}_k,\nabla f(x_k),H_k,\dktrue)] \nonumber \\
    &\ + \thalf (A_{\max}\beta + \theta\beta^2) \beta \E_k[\Delta q(x_k,\hat{\Tcal}_k,G_k,H_k,D_k)] \nonumber \\
    &\ + ((A_{\max}\beta + \theta\beta^2) \tau_{k-1} - A_{\min} \beta \tau_{\min}) \kappa_g \kappa_d M_\tau \Prob_k[\Ek3] + \theta \beta^2 \tau_{k-1} \kappa_g \kappa_d \sqrt{M} \nonumber \\
    \leq&\ -(A_{\min} - \thalf (A_{\max} + \theta\beta) \beta) \beta \E_k[\Delta q(x_k,\hat{\Tcal}_k,\nabla f(x_k),H_k,\dktrue)] \nonumber \\
    &\ + \tfrac14 (A_{\max} + \theta\beta) \beta^2 (\tau_{k-1} + \tau_{\min}) \zeta^{-1} M \nonumber \\
    &\ + \thalf (A_{\max} + \theta\beta) \beta^2 (\tau_{k-1} - \tau_{\min}) (\kappa_d \sqrt{M} (2\kappa_g + \sqrt{M}) + \kappa_H \kappa_d^2 (M + \tfrac32 \kappa_g^2)) \nonumber \\
    &\ + ((A_{\max} + \theta\beta) \tau_{k-1} - A_{\min} \tau_{\min}) \beta \kappa_g \kappa_d M_\tau \Prob_k[\Ek3] + \theta \beta^2 \tau_{k-1} \kappa_g \kappa_d \sqrt{M} \nonumber \\
    \leq&\ - \thalf A_{\min} \beta \E_k[\Delta q(x_k,\hat{\Tcal}_k,\nabla f(x_k),H_k,\dktrue)] \nonumber \\
    &\ + \kappa_{E_3} \beta \Prob_k[\Ek3] + \Mbar \beta^2, \label{eq.danny's_song}
  \end{align}
  where the final inequality follows due to the fact that $(A_{\max}+\theta\beta)\beta \leq A_{\min}$ holds by the definitions of $\beta$, $\gamma$, $A_{\min}$, and $A_{\max}$.
  
  Let us now combine \eqref{eq.danny_boy} and \eqref{eq.danny's_song} to prove a bound of the form in \eqref{eq.danny_elton_john} by considering two complementary events.  In particular, let $E_{k,\tau}$ be the event that $\Tktritrue < \tau_{k-1}$ and let $E_{k,\tau}^c$ be the complementary event that $\Tktritrue \geq \tau_{k-1}$.  Observe that whether $E_{k,\tau}$ or $E_{k,\tau}^c$ occurs is determined by the condition that $G_{[k-1]} = g_{k-1}$.  Hence, the bound in Lemma~\ref{lem:dqbound}---and in Lemma~\ref{lem:hattaudiff} as well, which is used below---holds even if one conditions on the occurence of $E_{k,\tau}$ or of $E_{k,\tau}^c$.  Consequently, the bounds in \eqref{eq.danny_boy} and \eqref{eq.danny's_song} hold even if one also conditions on the occurrence of $E_{k,\tau}$ or of $E_{k,\tau}^c$.  Let us now consider $E_{k,\tau}^c$ and $E_{k,\tau}$ in turn.  Conditioning on $E_{k,\tau}^c$, one finds from \eqref{eq.danny_boy}, \eqref{eq.danny's_song}, and the fact that $\Tktritrue \geq \tau_{k-1} \geq \Tcal_k = \hat\Tcal_k$ (by \eqref{eq:hattaudef}) in $E_{k,\tau}^c$ that
  \begin{align*}
    &\ \E_k [\phi(x_k + \alpha_k D_k, \Tcal_k) | E_{k,\tau}^c] - \E_k [\phi(x_k, \Tcal_k) | E_{k,\tau}^c] \\
    \leq&\ \E_k[-\Acal_k \Delta q(x_k,\Tcal_k,\nabla f(x_k),H_k,\dktrue) + \thalf \Acal_k \beta  \Delta q(x_k,\Tcal_k,G_k,H_k,D_k) | E_{k,\tau}^c] \\
    &\ + (\varsigma_{\max,k} \tau_{k-1} - \alpha_{\min,k}^{<} \tau_{\min}) \kappa_g \kappa_d M_\tau \Prob_k[\Ek3 | E_{k,\tau}^c] + \theta \beta^2 \tau_{k-1} \kappa_g \kappa_d \sqrt{M} \\
    =&\ \E_k[-\Acal_k \Delta q(x_k,\hat{\Tcal}_k,\nabla f(x_k),H_k,\dktrue) + \thalf \Acal_k \beta  \Delta q(x_k,\hat{\Tcal}_k,G_k,H_k,D_k) | E_{k,\tau}^c] \\
    &\ + (\varsigma_{\max,k} \tau_{k-1} - \alpha_{\min,k}^{<} \tau_{\min}) \kappa_g \kappa_d M_\tau \Prob_k[\Ek3 | E_{k,\tau}^c] + \theta \beta^2 \tau_{k-1} \kappa_g \kappa_d \sqrt{M} \\
    \leq&\ - \thalf A_{\min} \beta \E_k[\Delta q(x_k,\hat{\Tcal}_k,\nabla f(x_k),H_k,\dktrue) | E_{k,\tau}^c] \\
    &\ + \kappa_{E_3} \beta \Prob_k[\Ek3 | E_{k,\tau}^c] + \Mbar \beta^2.
  \end{align*}
  On the other hand, one finds from \eqref{eq.danny_boy}, \eqref{eq.danny's_song}, and Lemma~\ref{lem:hattaudiff} that
  \begin{align*}
    &\ \E_k [\phi(x_k + \Acal_k D_k, \Tcal_k) | E_{k,\tau}] - \E_k [\phi(x_k, \Tcal_k) | E_{k,\tau}] \\
    \leq&\ \E_k[-\Acal_k \Delta q(x_k,\Tcal_k,\nabla f(x_k),H_k,\dktrue) + \thalf \Acal_k \beta  \Delta q(x_k,\Tcal_k,G_k,H_k,D_k) | E_{k,\tau}] \\
    &\ + \E_k[-\Acal_k \Delta q(x_k,\hat{\Tcal}_k,\nabla f(x_k),H_k,\dktrue) + \thalf \Acal_k \beta  \Delta q(x_k,\hat{\Tcal}_k,G_k,H_k,D_k) | E_{k,\tau}] \\
    &\ - \E_k[-\Acal_k \Delta q(x_k,\hat{\Tcal}_k,\nabla f(x_k),H_k,\dktrue) + \thalf \Acal_k \beta  \Delta q(x_k,\hat{\Tcal}_k,G_k,H_k,D_k) | E_{k,\tau}] \\
    &\ + (\varsigma_{\max,k}\tau_{k-1} - \alpha_{\min,k}^{<} \tau_{\min}) \kappa_g \kappa_d M_\tau \Prob_k[\Ek3 | E_{k,\tau}] + \theta \beta^2 \tau_{k-1} \kappa_g \kappa_d \sqrt{M} \\ 
    \leq&\ \E_k[-\Acal_k \Delta q(x_k,\Tcal_k,\nabla f(x_k),H_k,\dktrue) + \thalf \Acal_k \beta  \Delta q(x_k,\Tcal_k,G_k,H_k,D_k) | E_{k,\tau}] \\
    &\ - \thalf A_{\min} \beta \E_k[\Delta q(x_k,\hat{\Tcal}_k,\nabla f(x_k),H_k,\dktrue) | E_{k,\tau}] \\
    &\ - \E_k[-\Acal_k \Delta q(x_k,\hat{\Tcal}_k,\nabla f(x_k),H_k,\dktrue) + \thalf \Acal_k \beta  \Delta q(x_k,\hat{\Tcal}_k,G_k,H_k,D_k) | E_{k,\tau}] \\
    &\ + \kappa_{E_3} \beta \Prob_k[\Ek3 | E_{k,\tau}] + \Mbar \beta^2 \\ 
    \leq&\ -\thalf A_{\min} \beta \E_k[\Delta q(x_k,\hat{\Tcal}_k,\nabla f(x_k),H_k,\dktrue) | E_{k,\tau}] \\
    &\ + (\varsigma_{\max,k} \tau_{k-1} - \alpha_{\min,k}^{<} \tau_{\min}) \kappa_d \kappa_g^2(1 + \thalf \kappa_H \kappa_d) \\
    &\ + \thalf (\varsigma_{\max,k} \tau_{k-1} - \alpha_{\min,k}^{<} \tau_{\min}) \kappa_d(2 + \kappa_H \kappa_d) (M + \kappa_g^2) \beta \\
    &\ + \kappa_{E_3} \beta \Prob_k[E_{k,3} | E_{k,\tau}] + \Mbar \beta^2 \\
    \leq&\ -\thalf A_{\min} \beta \E_k[\Delta q(x_k,\hat{\Tcal}_k,\nabla f(x_k),H_k,\dktrue) | E_{k,\tau}] \\
    &\ + \kappa_{E_3} \beta \Prob_k[E_{k,3} | E_{k,\tau}] + \Mbar \beta^2 + \kappa_{\Delta q,1} \beta + \kappa_{\Delta q,2} \beta^2.
  \end{align*}
  Hence, by the laws of total probability and expectation, one finds that
  \begin{align*}
    &\ \E_k [\phi(x_k + \Acal_k D_k, \Tcal_k)] - \E_k[\phi(x_k, \Tcal_k)] \\
    =&\ (\E_k [\phi(x_k + \Acal_k D_k, \Tcal_k)|E_{k,\tau}] - \E_k[\phi(x_k, \Tcal_k)| E_{k,\tau}]) \Prob_k[E_{k,\tau}] \\
    &\ + (\E_k [\phi(x_k + \Acal_k D_k, \Tcal_k)|E_{k,\tau}^c] - \E_k[\phi(x_k, \Tcal_k)| E_{k,\tau}^c]) \Prob_k[E_{k,\tau}^c] \\
    \leq&\ -\thalf A_{\min} \beta \E_k[\Delta q(x_k,\hat{\Tcal}_k,\nabla f(x_k),H_k,\dktrue)] \\
    &\ + \kappa_{E_3} \beta \Prob_k[E_{k,3}] + \Mbar \beta^2 + (\kappa_{\Delta q,1} \beta + \kappa_{\Delta q,2} \beta^2) \Prob_k[E_{k,\tau}].
  \end{align*}
  Summing this inequality for all $k \in [\kmax]$ yields
  \begin{align*}
    &\ \sum_{k=0}^{\kmax} (\E_k [\phi(x_k + \Acal_k D_k, \Tcal_k)] - \E_k[\phi(x_k, \Tcal_k)]) \\
    \leq&\ \sum_{k=0}^{\kmax} (-\thalf A_{\min} \beta \E_k[\Delta q(x_k,\hat{\Tcal}_k,\nabla f(x_k),H_k,\dktrue)]) \\
    &\ + \sum_{k=0}^{\kmax} (\kappa_{E_3} \beta \Prob_k[E_{k,3}] + (\kappa_{\Delta q,1} \beta + \kappa_{\Delta q,2} \beta^2) \Prob_k[E_{k,\tau}]) + (\kmax + 1) \Mbar \beta^2.
  \end{align*}
  
  Let us now turn to analyzing the overall behavior of the algorithm through iterations $k \in [\kmax]$.  Let $f_{G_{[\kmax]}}$ denote the probability density function of $G_{[\kmax]}$ and observe that, for all $k \in [\kmax]$, one finds $\Prob_k[E_{k,\tau}] \in \{0,1\}$ since, when conditioned on $E$ and $\Fcal_k$, the event $E_{k,\tau}$ is independent of $G_k$.  Therefore, by the bound above, the law of total expectation, Lemma \ref{lem:Eksumlsmax}, Lemma \ref{lem:Ktaulsmax}, and the union bound, it follows that, with probability at least $1-\delta$, one finds
  \begin{align}
    &\ \sum_{k=0}^{\kmax} (\E[\phi(X_k + \Acal_k D_k, \Tcal_k) | E, \Fcal_{k+1}] - \E[\phi(X_k,\Tcal_k) | E, \Fcal_{k+1}]) \nonumber \\
    =&\ \int_{g_{[\kmax]} \in \Fcal_{\kmax+1}} \sum_{k=0}^{\kmax} \E[\phi(x_k + \Acal_k D_k, \Tcal_k) - \phi(x_k,\Tcal_k) | E, G_{[\kmax]} = g_{[\kmax]}] \nonumber \\
    &\hspace{1.3in} \cdot f_{G_{[\kmax]}}(g_{[\kmax]}) \textrm{d}g_{[\kmax]} \nonumber \\
    \leq&\ \int_{g_{[\kmax]} \in \Fcal_{\kmax+1}} \Bigg(-\thalf A_{\min} \beta \sum_{k=0}^{\kmax} \E[\Delta q(x_k,\hat{\Tcal}_k,\nabla f(x_k),H_k,\dktrue) | E, G_{[\kmax]} = g_{[\kmax]}] \nonumber \\
    &\hspace{1.3in} + \kappa_{E_3} \beta \sum_{k=0}^{\kmax} \Prob[\Ek3 | E, G_{[\kmax]} = g_{[\kmax]}] \nonumber \\
    &\hspace{1.3in} + (\kappa_{\Delta q,1} \beta + \kappa_{\Delta q,2} \beta^2) \sum_{k=0}^{\kmax} \Prob[E_{k,\tau} | E, G_{[\kmax]} = g_{[\kmax]}] \Bigg) \nonumber \\
    &\hspace{1.3in} \cdot f_{G_{[\kmax]}}(g_{[\kmax]}) \textrm{d}g_{[\kmax]} \nonumber \\
    &\ + (\kmax+1) \Mbar \beta^2 \nonumber \\
    =&\ -\thalf A_{\min} \beta \sum_{k=0}^{\kmax} \E[\Delta q(X_k,\hat{\Tcal}_k,\nabla f(X_k),H_k,\Dktrue) | E, \Fcal_k] \nonumber \\
    &\ + \kappa_{E_3} \beta \sum_{k=0}^{\kmax} \Prob[\Ek3 | E, \Fcal_k] + (\kappa_{\Delta q,1} \beta + \kappa_{\Delta q,2} \beta^2) \E[|\Kcal_\tau| | E] + (\kmax + 1) \Mbar \beta^2 \nonumber \\
    \leq&\ -\thalf A_{\min} \beta \sum_{k=0}^{\kmax} \E[\Delta q(X_k,\hat{\Tcal}_k,\nabla f(X_k),H_k,\Dktrue) | E, \Fcal_k] \nonumber \\
    &\ + \kappa_{E_3} \beta (\ell(\smax,\hdelta/2)+1) \nonumber \\
    &\ + (\kappa_{\Delta q,1} \beta + \kappa_{\Delta q,2} \beta^2) \left\lceil\frac{\ell(\smax,\hdelta/2)+1}{p_{\tau}}\right\rceil + (\kmax+1) \Mbar \beta^2. \label{eq:complexityproof1}
  \end{align}
  The left-hand side of this inequality satisfies
  \begin{align}
    &\ \sum_{k=0}^{\kmax} (\E [\phi(X_k + \Acal_k D_k, \Tcal_k) | E, \Fcal_{k+1}] - \E[\phi(X_k,\Tcal_k) | E, \Fcal_{k+1}]) \nonumber \\
    =&\ \sum_{k=0}^{\kmax} \Bigg( \E[\Tcal_k (f(X_k + \Acal_k D_k) - f_{\min}) + \|c(X_k + \Acal_k D_k)\|_1 | E,\Fcal_{k+1}] \nonumber \\
    &\hspace{0.45in} - \E[\Tcal_k (f(X_k) - f_{\min}) + \|c(X_k)\|_1 | E,\Fcal_{k+1}] \Bigg). \label{eq.newnumb}
  \end{align}
  By the merit parameter updating strategy and the definition of the filtration (which guarantees that $\Tcal_{k-1}$ and $X_k$ are fully determined by $G_{[k-1]}$, and are therefore conditionally independent of $\Fcal_{k+1}$) that, for all $k \in [\kmax]$, one finds
  \bequationNN
    \baligned
      &\ -\E[\Tcal_k (f(X_k) - f_{\min}) + \|c(X_k)\|_1 | E, \Fcal_{k+1}] \\
      \geq&\ - \E[\Tcal_{k-1} (f(X_k) - f_{\min}) + \|c(X_k)\|_1 | E, \Fcal_{k+1}] \\
      =&\ - \E[\Tcal_{k-1} (f(X_k) - f_{\min}) + \|c(X_k)\|_1 | E, \Fcal_k].
    \ealigned
  \eequationNN
  Thus, from \eqref{eq.newnumb}, it follows that
  \begin{align*}
    &\ \sum_{k=0}^{\kmax} (\E[\phi(X_k + \Acal_k D_k, \Tcal_k) | E, \Fcal_{k+1}] - \E[\phi(X_k,\Tcal_k) | E, \Fcal_{k+1}]) \\
    &\geq \E[\Tcal_{\kmax} (f(X_{\kmax+1}) - f_{\min}) + \|c(X_{\kmax+1})\|_1 | E,\Fcal_{\kmax+1}] - \tau_{-1} (f_0 - f_{\min}) - \|c_0\|_1 \\
    &\geq - \tau_{-1} (f_0 - f_{\min}) - \|c_0\|_1.
  \end{align*}
  Combining this with \eqref{eq:complexityproof1}, one obtains that
  \begin{align*}
    &\ \frac{\beta}{\sum_{k=0}^{\kmax} \beta} \sum_{k=0}^{\kmax} \E[\Delta q(X_k,\hat{\Tcal}_k,\nabla f(X_k),H_k,\Dktrue) | E, \Fcal_k] \\
    \leq &\ 2\(\frac{\tau_{-1} (f_0 - f_{\min}) + \|c_0\|_1 + (\kmax+1) \Mbar \beta^2 + \kappa_{E_3} \beta (\ell(\smax,\hdelta/2)+1)}{A_{\min} {\sum_{k=0}^{\kmax} \beta}}\) \nonumber \\
    &\ + \frac{2(\kappa_{\Delta q,1} \beta + \kappa_{\Delta q,2} \beta^2)}{A_{\min} {\sum_{k=0}^{\kmax} \beta}} \left\lceil\frac{\ell(\smax,\hdelta/2)+1}{p_{\tau}}\right\rceil.
  \end{align*}
  Hence, by the definitions of $K^*$ and $\beta$, the desired conclusion follows. \qed
\end{proof}

The following corollary translates the result of the preceding theorem to a result pertaining to a stationary measure of \eqref{eq:fdef}; recall \eqref{eq:optconds}.

\begin{corollary} \label{coro:complexity}
  Under the assumptions, conditions, and definitions of Theorem~\ref{thm:complexity}, it holds with probability at least $1-\delta \in (0,1)$ that
  \begin{align*}
    &\ \E\left[\frac{\|\nabla f(X_{K^*}) + J_{K^*}^\top Y_{K^*}^{true}\|^2}{\kappa_H^2} + \|c(X_{K^*})\|_1 \Bigg| E\right] \\
    \leq&\ 2 (\kappa_\Psi + 1) \(\frac{\tau_{-1} (f_0 - f_{\min}) + \|c_0\|_1 + \Mbar \gamma^2}{\kappa_q \tau_{\min} A_{\min} \gamma \sqrt{\kmax+1}}\) \\
    &\ + 2(\kappa_\Psi + 1) \(\frac{\kappa_{E_3} \gamma (\ell(\smax,\hdelta/2)+1)}{\kappa_q \tau_{\min} A_{\min} \gamma \sqrt{\kmax+1}}\) \\
    &\ + (\kappa_\Psi + 1) \(\frac{2(\kappa_{\Delta q,1} \gamma + \kappa_{\Delta q,2} \gamma^2/\sqrt{\kmax+1})}{\kappa_q \tau_{\min} A_{\min} \gamma (\kmax+1)}\) \left\lceil\frac{\ell(\smax,\hdelta/2)+1}{p_{\tau}}\right\rceil.
  \end{align*}
  Hence, the complexity bound described in Section~\ref{sec.preview} $($see \eqref{eq.overview}$)$ holds.
\end{corollary}
\begin{proof}
  The result follows by Lemma \ref{lem:kappapsi}, Lemma \ref{lem:kappaq}, Theorem \ref{thm:complexity}, and \eqref{eq:sqpsystemdet}, which implies that $\|\nabla f(X_{K^*}) + J_{K^*}^\top Y_{K^*}^{true}\| = \|H_{k^*} D^{true}_{K^*}\| \leq \kappa_H \|D^{true}_{K^*}\|$.  The worst-case complexity bound described in Section~\ref{sec.preview} follows by combining this result with the definitions of $\kappa_{E_3}$, $\kappa_{\Delta q,1}$, $\kappa_{\Delta q,2}$, and Lemma \ref{lem:Orderlsd} in Appendix~\ref{app:probbounds}.
  \qed
\end{proof}

This result, as well as that of Theorem \ref{thm:complexity}, is proven under the assumption that $\smax \geq 1$. When $\smax = 0$, this result simplifies to a \textit{deterministic} complexity bound with the terms dependent on $\smax$ and $\delta$ ommitted. Under the condition $\smax = 0$, the proof follows by noting that $\Prob_k[\Ek3] = \Prob_k[E_{k,\tau}] = 0$ for all $k \in [\kmax]$ (where $E_{k,\tau}$ is defined in the proof of Theorem \ref{thm:complexity}) along with a similar argument to the proof of Theorem \ref{thm:complexity}.

Again, we remark that this result, when viewed in terms of the squared norm of the gradient of the Lagrangian, matches the worst-case complexity of the stochastic gradient method for the unconstrained setting~\cite{DDavis_DDrusvyatskiy_2019}.

\subsection{Complexity Result for Symmetric Sub-Gaussian Distributions} \label{subsec:subgauss}

In this section, we show that Assumption~\ref{assum:eventE}, Assumption~\ref{assum:ptau}, and the event $E$ in \eqref{def:E} occur with high probability when each stochastic gradient is unbiased and has a symmetric, sub-Gaussian distribution and the ratio parameter sequence remains constant.  For these purposes, we make the following assumption.

\begin{assumption} \label{assum:symmsubgauss}
  There exists $M \in \R{}_{>0}$ such that, for all $k \in [\kmax]$ and any realization $g_{[k-1]}$ of $G_{[k-1]}$, one finds that
  \begin{align}
    \E[G_k | G_{[k-1]} = g_{[k-1]}] &= \nabla f(x_k) \nonumber \\ \text{and}\ \ 
    \E[\exp(\|G_k-\nabla f(x_k)\|^2/M) | G_{[k-1]} = g_{[k-1]}] &\leq \exp(1), \label{eq:subgauss}
  \end{align}
  and the random vectors $G_k - \nabla f(x_k)$ and $\nabla f(x_k) - G_k$ have equal distributions. Finally, for all $k \in [\kmax]$, the ratio parameter $\Xi_k$ satisfies $\Xi_k = \xi_{\min}$.
\end{assumption}

Our first lemma shows that, under Assumption~\ref{assum:symmsubgauss}, Assumption~\ref{assum:ptau} holds.

\blemma \label{lem:subgptau}
  Under Assumptions~\ref{assum:fcsmooth}, \ref{assum:Hk}, and \ref{assum:symmsubgauss}, it follows for all $k \in [\kmax]$ and any realization $g_{[k-1]}$ of $G_{[k-1]}$ that
  \bequationNN
    \Prob[G_k^\top D_k + \max\{D_k^\top H_k D_k,0\} \geq \nabla f(x_k)^\top \dktrue + \max\{(\dktrue)^\top H_k \dktrue,0\} | g_{[k-1]}] \geq \thalf.
  \eequationNN
\elemma
\bproof
  Consider arbitrary $k \in [\kmax]$.  Let $Z_k$ be a basis for the null space of~$J_k$, which under Assumption \ref{assum:fcsmooth} is a matrix in $\R{n\times(n-m)}$.  Then, let $W_k \in \R{n-m}$ be a random vector such that $U_k = Z_k W_k$, and let $\wktrue \in \R{n-m}$ be such that $\uktrue = Z_k \wktrue$.  By \eqref{eq:sqpsystem}, $Z_k W_k = -Z_k (Z_k^\top H_k Z_k)^{-1} Z_k^\top (G_k + H_k v_k)$, so that
  \bequationNN
    \baligned
      &\ G_k^\top D_k + D_k^\top H_k D_k \\
     =&\ v_k^\top H_k^{1/2}(I - H_k^{1/2} Z_k(Z_k^\top H_k Z_k)^{-1} Z_k^\top H_k^{1/2})(H_k^{-1/2} G_k + H_k^{1/2} v_k)
    \ealigned
  \eequationNN
  and similarly
  \bequationNN
    \baligned
      &\ \nabla f(x_k)^\top \dktrue + (\dktrue)^\top H_k \dktrue \\
     =&\ v_k^\top H_k^{1/2}(I - H_k^{1/2} Z_k(Z_k^\top H_k Z_k)^{-1} Z_k^\top H_k^{1/2})(H_k^{-1/2} \nabla f(x_k) + H_k^{1/2} v_k).
    \ealigned
  \eequationNN
  Hence, the random variables
  \bequationNN
    \baligned
      &\ G_k^\top D_k + \max\{D_k^\top H_k D_k,0\} - \nabla f(x_k)^\top \dktrue - \max\{(\dktrue)^\top H_k \dktrue,0\} \\
     =&\ v_k^\top H_k^{1/2}(I - H_k^{1/2} Z_k(Z_k^\top H_k Z_k)^{-1} Z_k^\top H_k^{1/2})(H_k^{-1/2} (G_k -\nabla f(x_k)))
    \ealigned
  \eequationNN
  and
  \bequationNN
    \baligned
      &\ \nabla f(x_k)^\top \dktrue + \max\{(\dktrue)^\top H_k \dktrue,0\} - G_k^\top D_k - \max\{D_k^\top H_k D_k,0\} \\
     =&\ v_k^\top H_k^{1/2}(I - H_k^{1/2} Z_k(Z_k^\top H_k Z_k)^{-1} Z_k^\top H_k^{1/2})(H_k^{-1/2} (\nabla f(x_k) - G_k))
    \ealigned
  \eequationNN
  are equivalent in distribution by Assumption \ref{assum:symmsubgauss}. Therefore,
  \bequationNN
    \baligned
      &\ \Prob[G_k^\top D_k + \max\{D_k^\top H_k D_k,0\} - \nabla f(x_k)^\top \dktrue - \max\{(\dktrue)^\top H_k \dktrue,0\} \geq 0 | g_{[k-1]}] \\
     =&\ \Prob[\nabla f(x_k)^\top \dktrue + \max\{(\dktrue)^\top H_k \dktrue,0\} - G_k^\top D_k - \max\{D_k^\top H_k D_k,0\} \geq 0 | g_{[k-1]}]
    \ealigned
  \eequationNN
  and
  \bequationNN
    \baligned
      1 &= \Prob[G_k^\top D_k + \max\{D_k^\top H_k D_k,0\} - \nabla f(x_k)^\top \dktrue - \max\{(\dktrue)^\top H_k \dktrue,0\} \geq 0 | g_{[k-1]}] \\
        &+ \Prob[\nabla f(x_k)^\top \dktrue + \max\{(\dktrue)^\top H_k \dktrue,0\} - G_k^\top D_k - \max\{D_k^\top H_k D_k,0\} \geq 0 | g_{[k-1]}] \\
        &- \Prob[\nabla f(x_k)^\top \dktrue + \max\{(\dktrue)^\top H_k \dktrue,0\} - G_k^\top D_k - \max\{D_k^\top H_k D_k,0\} = 0 | g_{[k-1]}],
    \ealigned
  \eequationNN
  which combined leads to the desired conclusion. \qed
\eproof

Next, we state a result based on well-known properties of sub-Gaussian random variables. This lemma follows in the same manner as \cite[Lemma~5]{LiOrab20}.
\blemma \label{lem:subgboundedrv}
  Suppose Assumption \ref{assum:symmsubgauss} holds. Then, for any $\delta \in (0,1)$,
  \bequationNN
    \Prob\left[\underset{k \in [\kmax]}{\max} \|G_k - \nabla f(X_k)\| \leq \sqrt{M \left(1 + \log \left(\frac{\kmax+1}{\delta}\right)\right)}\right] \geq 1 - \delta.
  \eequationNN
\elemma

We conclude this subsection by showing that, under Assumption~\ref{assum:symmsubgauss}, both Assumption~\ref{assum:eventE} and event $E$ occur with high probability.

\blemma \label{lem:subgeventE}
  Suppose that Assumptions \ref{assum:fcsmooth}, \ref{assum:Hk}, and \ref{assum:symmsubgauss} hold, let $\kappa_v$ be defined as in \cite[Lemma 2.9]{ABerahas_FECurtis_DPRobinson_BZhou_2021}, let $\kappa_c$ be an upper bound for $\|c_k\|_2$ for all $k \in [\kmax]$ $($the existence of which follows under Assumption~\ref{assum:fcsmooth}$)$, and define
  \begin{align*}
    \kappa_{\tau_{\min}} :=&\ \kappa_v \Bigg(\kappa_g + \sqrt{M \left(1 + \log\left(\frac{\kmax+1}{\delta}\right)\right)} \\
    &\ + \frac{\kappa_H}{\zeta} \left(\sqrt{M \left(1 + \log\left(\frac{\kmax+1}{\delta}\right)\right)} + \kappa_g + \zeta + \kappa_H \kappa_v \kappa_c\right)\Bigg).
  \end{align*}
  Then, for any $\delta \in (0,1)$, it follows with probability at least $1-\delta$ that the conditions in Assumption~\ref{assum:eventE} and event $E$ hold with
  \begin{align*}
    M_\tau &= \sqrt{M \left(1 + \log\left( \frac{\kmax+1}{\delta}\right)\right)},\ \ \tau_{\min} = \frac{(1-\sigma)(1-\epsilon_\tau)}{\kappa_{\tau_{\min}}}, \\
    \text{and}\ \ \smax &= \min\left\{\kmax + 1, \left\lceil\frac{\log \left(\frac{\tau_{-1}\kappa_{\tau_{\min}}}{(1-\sigma)(1-\epsilon_\tau)}\right)}{\log\left(\frac{1}{1 - \epsilon_\tau}\right)}\right\rceil\right\}.
  \end{align*}
\elemma
\bproof
  By Lemma~\ref{lem:subgboundedrv}, the event considered in that lemma holds with probability at least $1-\delta$.  Hence, for the purposes of this proof, suppose that event holds.  By Jensen's inequality, convexity of $\exp(\cdot)$, and \eqref{eq:subgauss}, it follows that
  \bequationNN
    \E[\|G_k - \nabla f(x_k)\|^2| G_{[k-1]} = g_{[k-1]}] \leq M.
  \eequationNN
  In addition, it follows from the event in Lemma~\ref{lem:subgboundedrv} that \eqref{eq.sg+} holds with $M_\tau$ as stated in the lemma.  This accounts for Assumption~\ref{assum:eventE}.  Now consider event~$E$.  First, it follows from the arguments of \cite[Lemma 2.15 and 2.16]{ABerahas_FECurtis_DPRobinson_BZhou_2021} combined with the event in Lemma~\ref{lem:subgboundedrv} that $\Tcal_k \geq \tau_{\min}$ and $\Tktritrue \geq \tau_{\min}$ for all $k \in [\kmax]$ for $\tau_{\min}$ as stated in the lemma.  Second, it follows from the stated value of $\tau_{\min}$ and~\eqref{eq:smaxbound} that $|\{k \in [\kmax] : \Tcal_k < \Tcal_{k-1}\}| \leq \smax$ for $\smax$ as stated in the lemma.  Finally, the desired behavior of $\{\Xi_k\}$ follows from Assumption~\ref{assum:symmsubgauss}. \qed
\eproof

\subsection{Adaptive Ratio Parameter} \label{subsec:ratioparam}

In this section, we state a convergence rate result, which can be translated to a worst-case complexity result, that relaxes the definition of the event $E$ considered in prior sections.  In particular, we remove the assumption that $\Xi_k = \xi_{\min} \in (0,\infty)$ for all $k \in [\kmax]$.  Importantly, it has been proved in \cite[Lemma 3.5]{ABerahas_FECurtis_DPRobinson_BZhou_2021} that, under our remaining assumptions, there still exists $\xi_{\min} \in (0,\infty)$ such that $\Xi_k \geq \xi_{\min}$ for all $k \in [\kmax]$.  Therefore, by the manner in which the ratio parameter sequence is set, it follows that there exists a maximum number of $k \in [\kmax]$ such that $\Xi_k < \Xi_{k-1}$.  Denoting this limit as $r_{\max} \in \N{}$, it follows (for the same reasons as the bound for $\smax$ in \eqref{eq:smaxbound}) that
\bequationNN
  r_{\max} \leq \min\left\{\kmax+1, \left\lceil \frac{\log(\xi_{\min}/\xi_{-1})}{\log(1-\epsilon_\xi)} \right\rceil\right\}.
\eequationNN
To formalize our new assumption, we define
\bequationNN
    E_{\xi} := E(\kmax,\smax,r_{\max},\tau_{\min},\xi_{\min})
\eequationNN
as the event such that, in every realization of the algorithm, the merit parameters $\{\tau_k\}^{\kmax}_{k=0}$ and ratio parameters $\{\xi_k\}^{\kmax}_{k=0}$ satisfy
\bitemize
  \item $\tau_k \geq \tau_{\min} > 0$ for all $k \in [\kmax]$,
  \item $\tktritrue \geq \tau_{\min} > 0$ for all $k \in [\kmax]$,
  \item $\xi_k \geq \xi_{\min} > 0$ for all $k \in [\kmax]$,
  \item $|\{k \in [\kmax] : \tau_k < \tau_{k-1}\}| \leq \smax$, and
  \item $|\{k \in [\kmax] : \xi_k < \xi_{k-1}\}| \leq r_{\max}$.
\eitemize
Since $\{\Xi_k\}$ is bounded below deterministically, this event is identical to the event~$E$ defined in~\eqref{def:E}, except that one may have $\xi_0 > \xi_{\min}$.

For the purposes of this section, redefining
\bequationNN
    \Prob_k[\cdot] := \Prob[\cdot | E_{\xi}, G_{[k-1]} = g_{[k-1]}] \ \ \text{and} \ \ \E_k[\cdot] := \Prob[\cdot | E_{\xi}, G_{[k-1]} = g_{[k-1]}],
\eequationNN
our analysis of this case considers the following replacement of Assumption~\ref{assum:eventE}.

\bassumption\label{assum:eventExi}
  There exists $M \in \R{}_{>0}$ such that, for all $k \in [\kmax]$ and any realization $g_{[k-1]}$ of $G_{[k-1]}$, one finds that
  \bequationNN
    \E_k[G_k] = \nabla f(x_k) \quad \text{and} \quad
     \E_k[\|G_k - \nabla f(x_k)\|_2^2] \leq M.
  \eequationNN
  In addition, there exist $M_1 \in \R{}_{>0}$, $M_2 \in \R{}_{>0}$, and $M_3 \in \R{}_{>0}$ such that, for all $k \in [\kmax]$ and any realization $g_{[k-1]}$ of $G_{[k-1]}$, one finds that
  \bequationNN
    \baligned
      \text{either}\ \ \P_k[\nabla f(x_k)^\top(D_k - \dktrue)) < 0, \Tcal_k < \tau_{k-1}, \Xi_k = \xi_{k-1}] &= 0 \\
      \text{or}\ \ \E_k[\|G_k - \nabla f(x_k)\| |\nabla f(x_k)^\top(D_k - \dktrue)) < 0, \Tcal_k < \tau_{k-1}, \Xi_k = \xi_{k-1}] &\leq \kappa_1; \\
      \text{either}\ \ \P_k[\nabla f(x_k)^\top(D_k - \dktrue)) < 0, \Tcal_k = \tau_{k-1}, \Xi_k < \xi_{k-1}] &= 0 \\
      \text{or}\ \ \E_k[\|G_k - \nabla f(x_k)\| |\nabla f(x_k)^\top(D_k - \dktrue)) < 0, \Tcal_k = \tau_{k-1}, \Xi_k < \xi_{k-1}] &\leq \kappa_2; \\
      \text{and either}\ \ \P_k[\nabla f(x_k)^\top(D_k - \dktrue)) < 0, \Tcal_k < \tau_{k-1}, \Xi_k < \xi_{k-1}] &= 0 \\
      \text{or}\ \ \E_k[\|G_k - \nabla f(x_k)\| |\nabla f(x_k)^\top(D_k - \dktrue)) < 0, \Tcal_k < \tau_{k-1}, \Xi_k < \xi_{k-1}] &\leq \kappa_3.
    \ealigned
  \eequationNN
\eassumption

We claim that this assumption holds with high probability under Assumption \ref{assum:symmsubgauss} (without the assumption that $\Xi_k = \xi_{\min}$ for all $k \in [\kmax]$), a result that can be derived by modification of the techniques used in Section~\ref{subsec:subgauss}.

The complexity analysis for this case follows by essentially the same arguments as those used to derive a complexity result under Assumption \ref{assum:eventE}. A slight modification of Lemma \ref{lem:innerprod} is needed to include the three events related to the sign of $\nabla f(x_k)^\top(D_k-\dktrue)$ that appear in Assumption \ref{assum:eventExi} (as opposed to the one in Assumption \ref{assum:eventE}), which yields a result in terms of the probabilities of these three events.  Then, a slightly modified Lemma \ref{lem:Eksumlsmax} and the union bound can be applied two additional times to derive a complexity result. Since the analysis is a straightforward, but tedious extension of the results in Section~\ref{subsec:complexity}, we simply state the extension of \eqref{eq.overview} to this case without proof.
\begin{theorem}
  Suppose that Assumptions~\ref{assum:fcsmooth}, \ref{assum:Hk}, \ref{assum:ptau}, and \ref{assum:eventExi} hold and consider arbitrary $\delta \in (0,1)$.  Then, within $\kmax+1$ iterations, it holds with probability at least $1-\delta$ that the algorithm generates $x_{k^*}\in \R{n}$ corresponding to which there exists an associated Lagrange multiplier $y_{k^*}^{\text{true}} \in \R{m}$ such that
  \begin{align*}
    &\ \E[\|\nabla f_{k^*} + J_{k^*}^\top y_{k^*}^{true}\|^2 + \|c_{k^*}\| | E_{\xi}] \\
    =&\ \mathcal{O}\Bigg(\frac{\tau_{-1}(f_0 - \flow) + \|c_0\|_1 + M}{\sqrt{\kmax+1}} \\
     &\qquad + \frac{(\tau_{-1}\xi_{-1} - \tau_{\min}\xi_{\min})((\smax + r_{\max})\log(\kmax) + \log(1/\delta))}{\sqrt{\kmax+1}}\Bigg).
  \end{align*}
\end{theorem}

\section{Conclusion} \label{sec:conclusion}

We proved a worst-case complexity bound (in terms of iterations, function evaluations, and (stochastic) derivative evaluations) for the stochastic sequential quadratic optimization method for solving optimization problems involving a stochastic objective function and deterministic equality constraints proposed in \cite{ABerahas_FECurtis_DPRobinson_BZhou_2021}. While key to the practical performance of the algorithm, the adaptivity of the merit parameter introduced a number of theoretical challenges to overcome.  Under mostly standard assumptions, we proved that, with high probability, a measure of primal-dual stationarity decays at a rate of $k^{-4}$ (ignoring log factors), which translates into a worst-case complexity bound on par with the stochastic gradient method in the unconstrained setting.

While our analytical approach has been developed for an SQP method that uses the merit function in \eqref{eq:phidef}, it may be applicable to a wide variety of algorithmic frameworks for constrained stochastic optimization.  For example, our approach may be modified to apply to methods that  adaptively update critical parameters at each iteration, such as adaptive penalty methods \cite{RHByrd_GLopez-Calva_JNocedal_2012,RHByrd_JNocedal_RAWaltz,MMongeau_ASartenaer_1995}, adaptive augmented Lagrangian methods \cite{FECurtis_HJiang_DPRobinson_2015}, adaptive barrier methods \cite{JNocedal_AWachter_RAWaltz_2009}, and penalty-interior point methods \cite{FECurtis_2012}.  In addition, many constrained optimization algorithms generate (often unconstrained) subproblems defined by an auxilliary parameter sequence that is updated dynamically based off of the solution to the previous subproblem.  Algorithms of this type include penalty methods, augmented Lagrangian methods, and interior point methods \cite{JNocedal_SJWright_2006}.  In cases when the objective is stochastic, this auxilliary sequence would also be a random process, in which case analyzing the behavior of such a process would be paramount to proving a complexity result for such a method.  We believe that the techniques that we have devised for this paper are broadly applicable and foundational for performing complexity analyses of deterministically constrained stochastic optimization methods.

\bibliographystyle{plain}
\bibliography{references}

\section{Proof of Theorem~\ref{th.deterministic} (Deterministic Algorithm Complexity)}\label{app.deterministic}

In this appendix, we prove Theorem~\ref{th.deterministic}, which states a worst-case complexity bound for Algorithm 2.1 of \cite{ABerahas_FECurtis_DPRobinson_BZhou_2021}.  We refer to quantities defined and employed in the analysis in \cite{ABerahas_FECurtis_DPRobinson_BZhou_2021}.  In particular, in this appendix, for all $k \in \N{}$, we suppose that $g_k = \nabla f(x_k)$ and $d_k = u_k + v_k$ with $u_k \in \Null(J_k)$ and $v_k \in \Range(J_k^\top)$ is the search direction computed by solving the SQP subproblem with $g_k = \nabla f(x_k)$.  As seen in \cite{ABerahas_FECurtis_DPRobinson_BZhou_2021}, the convergence properties of Algorithm~2.1 in that paper are driven by reductions in a model of the merit function in each iteration.  Our first lemma proves a useful lower bound for such a reduction.

\begin{lemma}\label{lem:suffdecrease}
  Define $(\kappa_{uv},\kappa_H,\kappa_v,\tau_{\min},\zeta,\sigma) \in (0,\infty)^5 \times (0,1)$ as in \cite{ABerahas_FECurtis_DPRobinson_BZhou_2021} and let
  \begin{equation} \label{eq:kappahat}
    \hat{\kappa} := \min\left\{1,\frac{1}{(1+\kappa_{uv})\kappa_v\kappa_H^2}\right\}\ \ \text{and}\ \ \tilde\kappa := \tfrac14 \zeta \kappa_{uv} \kappa_v \hat{\kappa}.
  \end{equation}
  Then, for any $\varepsilon \in (0,1)$, if $\|g_k + J_k^\top y_k\| > \varepsilon$ and/or $\sqrt{\|c_k\|_1} > \varepsilon$, then
  \begin{equation} \label{eq:suffdecrease}
    \Delta q(x_k, \tau_k, g_k, H_k, d_k) \geq \min\left\{\sigma \hat{\kappa},
    \tau_{\min} \tilde\kappa \right\} \varepsilon^2.
  \end{equation}
\end{lemma}
\begin{proof}
  Consider arbitrary $(\varepsilon,k) \in (0,1) \times \N{}$ such that $\|g_k + J_k^\top y_k\| > \varepsilon$ and/or $\sqrt{\|c_k\|_1} > \varepsilon$.  Let us consider two cases.  First, suppose that $\|c_k\|_1 > \hat{\kappa} \varepsilon^2$. Then, by \cite[equation (2.9)]{ABerahas_FECurtis_DPRobinson_BZhou_2021},
  \bequationNN
    \Delta q(x_k, \tau_k, g_k, H_k, d_k) \geq \tfrac12 \tau_k \max\{d_k^\top H_k d_k 0\} + \sigma \|c_k\|_1 \geq \sigma \|c_k\|_1 \geq \sigma \hat{\kappa} \varepsilon^2,
  \eequationNN
  which implies \eqref{eq:suffdecrease}, as desired.  Second, suppose that $\|c_k\|_1 \leq \hat{\kappa} \varepsilon^2 \leq \varepsilon^2$, which by the definition of $(\varepsilon,k)$ implies that $\|g_k + J_k^\top y_k\| > \varepsilon$.  It follows from this fact that $\|d_k\| > \varepsilon/\kappa_H$; indeed, if $\|d_k\| \leq \varepsilon/\kappa_H$, then by  \cite[equation (2.6) and Assumption~2.4]{ABerahas_FECurtis_DPRobinson_BZhou_2021} one would find
  \bequationNN
    \|g_k + J_k^\top y_k\| = \|H_k d_k\| \leq \kappa_H \|d_k\| \leq \varepsilon,
  \eequationNN
  which is a contradiction.  Hence, $\|d_k\| > \varepsilon/\kappa_H$, and by \cite[Lemma~2.9]{ABerahas_FECurtis_DPRobinson_BZhou_2021}, it follows that $\|v_k\|^2 \leq \kappa_v \|c_k\| \leq \kappa_v \|c_k\|_1 \leq \kappa_v \hat{\kappa} \varepsilon^2$, which combined shows that
  \bequationNN
    \varepsilon^2/\kappa_H^2 < \|d_k\|^2 = \|u_k\|^2 + \|v_k\|^2 \leq \|u_k\|^2 + \kappa_v \hat{\kappa} \varepsilon^2.
  \eequationNN
  From this fact and the definition of $\hat{\kappa}$, it follows that
  \bequationNN
    \|u_k\|^2 > \frac{\varepsilon^2}{\kappa_H^2} - \kappa_v \hat{\kappa} \varepsilon^2 \geq \frac{\varepsilon^2}{\kappa_H^2}\(1 - \frac{1}{(1+\kappa_{uv})}\) = \frac{\kappa_{uv}\varepsilon^2}{(1+\kappa_{uv})\kappa_H^2} \geq \kappa_{uv} \kappa_v \hat{\kappa} \varepsilon^2 \geq \kappa_{uv} \|v_k\|^2,
  \eequationNN
  which along with \cite[Lemma~2.10]{ABerahas_FECurtis_DPRobinson_BZhou_2021} implies $d_k^\top H_k d_k \geq \thalf \zeta \|u_k\|^2 \geq \thalf \zeta \kappa_{uv} \kappa_v \hat{\kappa} \varepsilon^2$.  Thus,
  \begin{align*}
    \Delta q(x_k, \tau_k, g_k, H_k, d_k) &\geq \thalf \tau_k \max\{d_k^\top H_k d_k 0\} + \sigma \|c_k\|_1 \geq \tfrac14 \tau_{\min} \zeta \kappa_{uv} \kappa_v \hat{\kappa} \varepsilon^2,
  \end{align*}
  which implies \eqref{eq:suffdecrease}, as desired. \qed
\end{proof}

We now prove Theorem~\ref{th.deterministic}, further details of which are provided in the statement below.

\begin{theorem}\label{th.deter_complex}
  Define $(\tau_{-1},\flow,\alpha_{\min},\tau_{\min},\eta,\sigma) \in (0,\infty)^4 \times (0,1)^2$ as in \cite{ABerahas_FECurtis_DPRobinson_BZhou_2021} and $(\hat\kappa,\tilde\kappa) \in (0,1] \times (0,\infty)$ as in \eqref{eq:kappahat}.  Then, for any $\varepsilon \in (0,1)$, Theorem~\ref{th.deterministic} holds with \eqref{eq.deterministic_complexity} given by
  \bequationNN
    \overline{K}_\varepsilon := \(\frac{\tau_{-1}(f_0 - \flow) + \|c_0\|_1}{\eta \alpha_{\min} \min\{\sigma \hat{\kappa}, \tau_{\min} \tilde\kappa \}}\) \varepsilon^{-2}.
  \eequationNN
\end{theorem}
\begin{proof}
  To derive a contradiction, suppose \eqref{eq:eps1lon}
does not hold for all $k \in \{0,\dots,\overline{K}_\varepsilon\}$.  Then, along with Lemma~\ref{lem:suffdecrease} and~\cite[equation~(2.10) and Lemma~2.17]{ABerahas_FECurtis_DPRobinson_BZhou_2021}, it follows for all such $k$ that
  \bequationNN
    \phi (x_k + \alpha_k d_k, \tau_k) - \phi (x_k, \tau_k) \leq - \eta \alpha_k \Delta q(x_k, \tau_k, g_k, H_k, d_k) \leq - \eta \alpha_{\min} \min\{\sigma\hat\kappa, \tau_{\min} \tilde\kappa\} \varepsilon^2.
  \eequationNN
  By the definition of $\phi$, this means for all such $k$ that
  \bequationNN
    \tau_k f_{k+1} + \|c_{k+1}\|_1 \leq \tau_k f_k + \|c_k\|_1 - \eta \alpha_{\min} \min\{\sigma\hat\kappa, \tau_{\min} \tilde\kappa\} \varepsilon^2.
  \eequationNN
  Summing this inequality for all $k \in \{0,\dots,\overline{K}_\varepsilon\}$, one can deduce that
  \bequationNN
    \|c_{\overline{K}_\varepsilon+1}\|_1 - \|c_0\|_1 + \tau_{\overline{K}_\varepsilon}f_{\overline{K}_\varepsilon+1} - \tau_0 f_0 + \sum_{k=1}^{\overline{K}_\varepsilon} f_k (\tau_{k-1} - \tau_k) \leq - (\overline{K}_\varepsilon+1) \eta \alpha_{\min} \min\{\sigma\hat\kappa, \tau_{\min} \tilde\kappa\} \varepsilon^2.
  \eequationNN
  Since $\{\tau_k\}$ is monotonically nonincreasing, $\|c_{\overline{K}_\varepsilon+1}\|_1 \geq 0$, and $f_k \geq \flow$ for all $k \in \N{}$,
  \bequationNN
    - \|c_0\|_1 + \tau_{\overline{K}_\varepsilon}\flow - \tau_0 f_0 + \flow (\tau_{0} - \tau_{\overline{K}_\varepsilon}) \leq - (\overline{K}_\varepsilon+1) \eta \alpha_{\min} \min\{\sigma\hat\kappa, \tau_{\min} \tilde\kappa\} \varepsilon^2.
  \eequationNN
  Rearranging this inequality, one arrives at the conclusion that
  \bequationNN
    \overline{K}_\varepsilon+1 \leq \(\frac{\tau_0 (f_0 - \flow) + \|c_0\|_1}{\eta \alpha_{\min} \min\{\sigma\hat\kappa, \tau_{\min} \tilde\kappa\}}\) \varepsilon^{-2} \leq \(\frac{\tau_{-1} (f_0 - \flow) + \|c_0\|_1}{\eta \alpha_{\min} \min\{\sigma\hat\kappa, \tau_{\min} \tilde\kappa\}}\) \varepsilon^{-2} \equiv \overline{K}_\varepsilon,
  \eequationNN
  which is a contradiction.  Therefore, one arrives at the desired conclusion that Algorithm 2.1 yields an iterate satisfying \eqref{eq:eps1lon} in at most $\overline{K}_\varepsilon$ iterations. \qed
\end{proof}

\section{Proofs of Lemmas \ref{lem:Eksumlsmax} and \ref{lem:Ktaulsmax}} \label{app:probbounds}

In this appendix, we prove Lemmas~\ref{lem:Eksumlsmax} and \ref{lem:Ktaulsmax}.  Toward this end, we prove for any $\delta \in (0,1)$ with $\hat{\delta}$ as defined in \eqref{eq:hatdelta} and $\ell(\smax,\hat{\delta})$ as defined in \eqref{eq:lsd}, one finds
\bequation \label{eq:sumptaudecbound}
    \Prob\left[\sum_{i=0}^{\kmax} \Prob[\Tcal_i < \Tcal_{i-1} | E, G_{[i-1]}] \leq \ell(s_{\max},\hat{\delta})+1 \Bigg| E \right] \geq 1-\delta.
\eequation
We build to this result, ultimately proved as Lemma~\ref{coro:probbounded}, with a series of preliminary lemmas.

As our first preliminary result, we state a particular form of Chernoff's bound in the following lemma, which will prove instrumental in deriving \eqref{eq:sumptaudecbound}.

\blemma\label{lem:chernoff}
  For any $k \in \N{}$, let $\{Y_0,\dots,Y_k\}$ be independent Bernoulli random variables.  Then, for any $s \in \N{}$ and $\bar{\delta} \in (0,1)$, it follows that
  \bequation \label{eq:chernoffeq1}
    \mu := \sum_{j=0}^k \P[Y_j = 1] \geq \ell(s,\bar{\delta})\ \ \implies \ \ \P\left[\sum_{j=0}^k Y_j \leq s\right] \leq \bar{\delta}.
  \eequation
\elemma
\bproof
  Suppose that $\mu \geq \ell(s,\bar\delta)$.  By the multiplicative form of Chernoff's bound, it follows for $\rho := 1 - s/\mu$ (which is in the interval $(0,1)$ by (\ref{eq:chernoffeq1})) that
  \bequationNN
    \P\left[\sum_{j=0}^{k} Y_j \leq s \right] \leq e^{-\thalf \mu \rho^2} = e^{- \thalf \mu (1-s/\mu)^2}.
  \eequationNN
  Hence, to prove the result, all that remains is to show that $e^{-\thalf \mu (1-s/\mu)^2} \leq \bar{\delta}$, i.e., that $-\thalf \mu(1-s/\mu)^2 \leq \log(\bar{\delta})$.  Using $\log(\bar{\delta}) = -\log(1/\bar{\delta})$, this inequality is equivalent to
  \bequationNN
    0 \leq \thalf \mu(1-s/\mu)^2 - \log(1/\bar{\delta}) = \tfrac{1}{2\mu} (\mu - s)^2 - \log(1/\bar{\delta}),
  \eequationNN
  which holds if and only if $\mu^2 - 2\mu (s+\log(1/\bar{\delta})) + s^2 \geq 0$.  Viewing the left-hand side of this inequality as a convex quadratic function in $\mu$, one finds that the inequality holds as long as $\mu$ is greater than or equal to the positive root of the quadratic, i.e.,
  \bequationNN
    s + \log(1/\bar{\delta}) + \sqrt{(s+\log(1/\bar{\delta}))^2 - s^2} = s + \log(1/\bar{\delta}) + \sqrt{\log(1/\bar{\delta})^2 + 2s\log(1/\bar{\delta})}.
  \eequationNN
  This holds since $\mu \geq \ell(s,\bar\delta)$; hence, the result is proved. \qed
\eproof

Now, we turn our attention to proving \eqref{eq:sumptaudecbound}. For any realization of a run of the algorithm up to iteration $k \in [\kmax]$, let $w_k$ denote the number of times that the merit parameter has been decreased up to the beginning of iteration $k$ and let $\pbar_k$ denote the probability that the merit parameter is decreased during iteration $k$.  The \emph{signature} of a realization up to iteration $k \in \N{}$ is $(\pbar_0,\dots,\pbar_k,w_0,\dots,w_k)$, which encodes all of the pertinent information regarding the behavior of the merit parameter sequence up to the start of iteration $k$.

One could imagine using all possible signatures to define a tree whereby each node contains a subset of all realizations of the algorithm. To construct such tree, one could first consider the root node, which could be denoted by $\tilde{N}(\pbar_0,w_0)$, where $\pbar_0$ is uniquely defined by the starting conditions of our algorithm and $w_0 = 0$. All realizations of our algorithm follow the same initialization, so $\pbar_0$ and $w_0$ would be in the signature of every realization. Now, one could define a node $\tilde{N}(\pbar_{[k]},w_{[k]})$ at depth $k \in [\kmax]$ (where the root node has a depth of $0$) in the tree as the set of all realizations of our algorithm for which the signature of the realization up to iteration $k$ is $(\pbar_0,\dots,\pbar_k,w_0,\dots,w_k)$.  One could then define the edges in the tree by connecting nodes at adjacent levels, where node $\tilde{N}(\pbar_{[k]},w_{[k]})$ is connected to node $\tilde{N}(\bar{p}_{[k]},\pbar_{k+1},w_{[k]},w_{k+1})$ for any $\pbar_{k+1} \in [0,1]$ and $w_{k+1} \in \{w_k,w_k+1,\dots\}$.

Unfortunately, the construction described in the previous paragraph may lead to nodes in the tree representing realizations with probability zero occurrence. In order to remedy this, we instead construct a tree where the nodes contain all realizations whose probability signatures fall within specified intervals. To define such intervals, consider arbitrary $B \in \N{} \setminus \{0\}$ and let us restrict the sequence of values $p_{[k]}$ used to define our nodes as those with
\begin{equation}
  p_{[k]} = (p_0, \dots, p_k) \in \left\{0, \tfrac{1}{B}, \dots, \tfrac{B-1}{B}\right\}^{k+1}.
\end{equation}
For $p \in \{0, 1/B, \dots, (B-1)/B\}$, these define the open probability intervals $\iota(p)$ given by
\begin{equation*}
    \iota(p) = \begin{cases}
        \left[p, p + \tfrac{1}{B}\right) & \text{if} \;p \in \left\{0, \tfrac{1}{B}, \dots, \tfrac{B-2}{B}\right\}, \\
        \left[\tfrac{B-1}{B}, 1\right] & \text{if} \;p = \tfrac{B-1}{B}.
    \end{cases}
\end{equation*}

Now, we can construct our tree as follows. As before, first consider the root node, which we denote by $N(p_0,w_0)$, where $p_0 \in \{0, 1/B, \dots, (B-1)/B\}$ is uniquely defined by the starting conditions of our algorithm so that $\Prob[\Tcal_0 < \tau_{-1} | E] \in \iota(p_0)$ and $w_0 = 0$. All realizations of our algorithm follow the same initialization, so with $\bar{p}_0 = \Prob[\Tcal_0 < \tau_{-1} | E]$ one finds that $\bar{p}_0 \in \iota(p_0)$ and $w_0$ are in the signature of every realization.  We define a node $N(p_{[k]},w_{[k]})$ at depth $k \in [\kmax]$ as the set of all realizations for which the signature of the realization at iteration $k$ exactly matches $w_{[k]}$ and has probabilities that fall within the intervals defined by $p_{[k]}$; i.e., a realization with signature $(\pbar_{[k]},w_{[k]})$ is a member of $N(p_{[k]},w_{[k]})$ if and only if, for all $j \in [k]$, one finds that $\pbar_j \in \iota(p_j)$.  The edges in the tree connect nodes in adjacent levels, where $N(p_{[k]},w_{[k]})$ is connected to $N(p_{[k]},p_{k+1},w_{[k]},w_{k+1})$ for any $p_{k+1} \in \{0,1/B,\dots,(B-1)/B\}$ and $w_{k+1} \in \{w_k,w_k+1,\dots\}$.

Notationally, since the behavior of a realization of the algorithm up to iteration $k \in \N{}$ is completely determined by the initial conditions and the realization of $G_{[k-1]}$, we say that a realization described by $G_{[k-1]}$ belongs in node $N(p_{[k]},w_{[k]})$ by writing that
\bequationNN
  G_{[k-1]} \in N(p_{[k]},w_{[k]}).
\eequationNN
The initial condition, denoted for consistency as $G_{[-1]} \in N(p_0,w_0)$, occurs with probability one.  Based on the description above, the nodes of our tree satisfy: For any node at a depth of $k \geq 2$, the event $G_{[k-1]} \in N(p_{[k]},w_{[k]})$ occurs if and only if
\bequation\label{eq.defdef}
  \baligned
    \Prob[\Tcal_k < \Tcal_{k-1} | E, G_{[k-1]}] &\in \iota(p_k), \\
    S_{k-1} := \sum_{i=0}^{k-1} \Ical[\Tcal_i < \Tcal_{i-1}] &= w_k, \\
    \text{and}\ \ G_{[k-2]} &\in N(p_{[k-1]},w_{[k-1]}).
  \ealigned
\eequation

Let us now define certain important sets of nodes in the tree.  First, let
\bequationNN
  \Lgood := \left\{N(p_{[k]},w_{[k]}) : \left(\sum_{i=0}^k p_i \leq \ell(\smax,\hat{\delta}) + 1\right) \land (w_{k} = s_{\max} \lor k = \kmax) \right\}
\eequationNN
be the set of nodes at which the sum of the elements of $p_{[k]}$ is sufficiently small and either $w_k$ has reached $\smax$ or $k$ has reached $\kmax$.  Second, let
\bequationNN
  \Lbad := \left\{N(p_{[k]},w_{[k]}) : \sum_{i=0}^{k} p_i > \ell(s_{\max},\hat{\delta})+1 \right\}
\eequationNN
be the nodes in the complement of $\Lgood$ at which the sum of the elements of $p_{[k]}$ has exceeded the threshold $\ell(s_{\max},\hat\delta)+1$.  Going forward, we restrict attention to the tree defined by the root node and all paths from the root node that terminate at a node contained in $\Lgood \cup \Lbad$.  It is clear from this restriction and the definitions of $\Lgood$ and $\Lbad$ that this tree is finite with the elements of $\Lgood \cup \Lbad$ being leaves.

Let us now define relationships between nodes.  The parent of a node is defined as
\bequationNN
  P(N(p_{[k]},w_{[k]})) = N(p_{[k-1]},w_{[k-1]}).
\eequationNN
On the other hand, the children of node $N(p_{[k]},w_{[k]})$ are defined as
\bequationNN
  C(N(p_{[k]},w_{[k]})) = \bcases \{N(p_{[k]},p_{k+1},w_{[k]},w_{k+1})\} & \text{if $N(p_{[k]},w_{[k]}) \not\in \Lgood \cup \Lbad$} \\ \emptyset & \text{otherwise.} \ecases
\eequationNN
This ensures that paths down the tree terminate at nodes in $\Lgood \cup \Lbad$, making these nodes the leaves of the tree.  For convenience in the remainder of our discussions, let $C(\emptyset) = \emptyset$.

We define the height of node $N(p_{[k]},w_{[k]})$ as the length of the longest path from $N(p_{[k]},w_{[k]})$ to a leaf node, i.e., the height is denoted as
\bequationNN
  h(N(p_{[k]},w_{[k]})) := \left(\min \{j \in \N{} \setminus \{0\} : C^j(N(p_{[k]},w_{[k]})) = \emptyset\}\right)-1,
\eequationNN
where $C^j(N(p_{[k]},w_{[k]}))$ is shorthand for applying the mapping $C(\cdot)$ consecutively $j$ times. From this definition, $h(N(p_{[k]},w_{[k]})) = 0$ for all $N(p_{[k]},w_{[k]}) \in \Lgood \cup \Lbad$.

Next, let us define two more sets of nodes that will be useful later.  Let $C_{dec}(N(p_{[k]},w_{[k]}))$ denote the set of children of $N(p_{[k]},w_{[k]})$ such that the merit parameter decreases and let $\Cdec^c(N(p_{[k]},w_{[k]}))$ denote set of children of $N(p_{[k]},w_{[k]})$ such that it does not decrease, so
\begin{align}
  \Cdec(N(p_{[k]},w_{[k]})) := \{
  & N(p_{[k]},p_{k+1},w_{[k]},w_{k+1}) : \nonumber \\
  & (N(p_{[k]},p_{k+1},w_{[k]},w_{k+1}) \in C(N(p_{[k]},w_{[k]}))) \nonumber \\
  & \land (w_{k+1} = w_{k} + 1)\} \label{eq:Cdec}
\end{align}
and
\begin{align}
  \Cdec^c(N(p_{[k]},w_{[k]})) := \{
  & N(p_{[k]},p_{k+1},w_{[k]},w_{k+1}) : \nonumber \\
  & (N(p_{[k]},p_{k+1},w_{[k]},w_{k+1}) \in C(N(p_{[k]},w_{[k]}))) \nonumber \\
  & \land (w_{k+1} = w_{k})\}. \label{eq:CdecC}
\end{align}

Finally, let us define the event $\EbadB$ as the event that for some $j \in [\kmax]$ one finds
\bequation\label{eq:EbadB}
  \left(\sum_{i=0}^{j} \Prob[\Tcal_i < \Tcal_{i-1} | E, G_{[i-1]}] > \ell(s_{\max},\hat{\delta}) +  \tfrac{\kmax+1}{B} + 1\right).
\eequation
With respect to our goal of proving \eqref{eq:sumptaudecbound}, the event $\EbadB$ is of interest since it is the event that the given probabilities accumulated up to iteration $j \in [\kmax]$ (and beyond) exceed the threshold found in~\eqref{eq:sumptaudecbound} plus a factor that is inversely proportional to $B$.

Let us now prove some properties of the leaf nodes.

\blemma \label{lem:leafnodes}
  For any $k \in [\kmax]$ and $(p_{[k]},w_{[k]})$ with $N(p_{[k]},w_{[k]}) \in \Lgood$, one finds
  \bequationNN
    \Prob[G_{[k-1]} \in N(p_{[k]},w_{[k]}) \land \EbadB | E] = 0.
  \eequationNN
  On the other hand, for all $k \in [\kmax]$ and $(p_{[k]},w_{[k]})$ with $N(p_{[k]},w_{[k]}) \in \Lbad$, one finds
  \begin{align*}
    &\Prob[G_{[k-1]} \in N(p_{[k]},w_{[k]}) \land \EbadB | E] \\
    &\leq \hat{\delta} \prod_{i=1}^k \Prob\left[\Prob[\Tcal_i < \Tcal_{i-1} | E, G_{[i-1]}] \in \iota(p_i) \big| E, S_{i-1} = w_{i}, G_{[i-2]} \in N(p_{[i-1]},w_{[i-1]})\right].
  \end{align*}
\elemma
\bproof
  Consider an arbitrary index $k \in [\kmax]$ and an arbitrary pair $(p_{[k]},w_{[k]})$ such that $N(p_{[k]},w_{[k]}) \in \Lgood$.  By the definition of $\Lgood$, it follows that
  \begin{equation}\label{eq.garble}
    \sum_{i=0}^k p_i \leq \ell(\smax, \hat{\delta}) + 1.
  \end{equation}
  Since the maximum depth of a node is $\kmax$, it follows from \eqref{eq.garble} that
  \begin{align*}
    &\Prob\left[\sum_{i=0}^k \Prob[\Tcal_i < \Tcal_{i-1} | E, G_{[i-1]}] > \ell(\smax,\hat{\delta}) + \tfrac{\kmax+1}{B} + 1 \Big| E, G_{[k-1]} \in N(p_{[k]},w_{[k]})\right] \\
    &\leq \Prob\left[\sum_{i=0}^k \left(p_i + \tfrac{1}{B}\right) > \ell(\smax,\hat{\delta}) + \tfrac{\kmax+1}{B} + 1 \Big| E, G_{[k-1]} \in N(p_{[k]},w_{[k]})\right] \\
    &\leq \Prob\left[\ell(\smax,\hat{\delta}) + \tfrac{k+1}{B} + 1 > \ell(\smax,\hat{\delta}) + \tfrac{\kmax+1}{B} + 1 \Big| E, G_{[k-1]} \in N(p_{[k]},w_{[k]})\right] = 0.
  \end{align*}
  Therefore, for any $j \in \{1,\dots,k\}$, one finds from conditional probability that
  \begin{align*}
    &\Prob\left[G_{[j-1]} \in N(p_{[j]},w_{[j]}) \land \text{\eqref{eq:EbadB} holds} | E\right] \\
    &=\Prob\left[\sum_{i=0}^{j} \Prob[\Tcal_i < \Tcal_{i-1} | E, G_{[i-1]}] > \ell(\smax,\hat{\delta}) + \tfrac{\kmax+1}{B} + 1 \Big| E, G_{[j-1]} \in N(p_{[j]},w_{[j]})\right] \\
    &\quad \cdot \Prob\left[G_{[j-1]} \in N(p_{[j]},w_{[j]}) | E\right] = 0.
  \end{align*}
  In addition, \eqref{eq:EbadB} cannot hold for $j=0$ since $\ell(\smax,\hat{\delta})+1 > 1$.  Hence, along with the conclusion above, it follows that $\EbadB$ does not occur in any realization whose signature up to iteration $j \in \{1,\dots,k\}$ falls into a node along any path from the root note to $N(p_{[k]},w_{[k]})$.  Now, by the definition of $\Lgood$, at least one of $w_{k} = \smax$ or $k = \kmax$ holds.  Let us consider each case in turn.  If $k = \kmax$, then it follows by the preceding arguments that 
  \begin{equation*}
    \Prob\left[\sum_{i=0}^{\kmax} \Prob[\Tcal_i < \Tcal_{i-1} | E, G_{[i-1]}] \leq \ell(\smax,\hat{\delta}) + \tfrac{\kmax+1}{B} + 1 \Big| E, G_{[k-1]} \in N(p_{[k]},w_{[k]})\right] = 1.
  \end{equation*}
  Otherwise, if $w_{k} = \smax$, then it follows by Assumption~\ref{assum:eventE} that $\Prob[\Tcal_i < \Tcal_{i-1} | E, G_{[i-1]}] = 0$ for all $i \in \{k,\dots,\kmax\}$, and therefore the equation above again follows.  Overall, it follows that $\Prob[G_{[k-1]} \in N(p_{[k]},w_{[k-1]}) \land \EbadB | E] = 0$, as desired.

  Now consider arbitrary $k \in \N{}$ and $(p_{[k]},w_{[k]})$ with $N(p_{[k]},w_{[k]}) \in \Lbad$.  One finds
  \begin{align*}
    &\ \Prob[G_{[k-1]} \in N(p_{[k]},w_{[k]}) \land \EbadB | E] \\
    =&\ \Prob[\EbadB | E, G_{[k-1]} \in N(p_{[k]},w_{[k]})] \cdot \Prob[G_{[k-1]} \in N(p_{[k]},w_{[k]}) | E] \\
    \leq&\ \Prob[G_{[k-1]} \in N(p_{[k]},w_{[k]}) | E].
  \end{align*}
  Hence, using the initial condition that $G_{[-1]} \in N(p_0,w_0)$, it follows that
  \begin{align}
    &\ \Prob[G_{[k-1]} \in N(p_{[k]},w_{[k]}) \land \EbadB | E] \nonumber \\
    \leq&\ \Prob[G_{[k-1]} \in N(p_{[k]},w_{[k]}) | E] = \Prob\left[\text{\eqref{eq.defdef} holds} \big| E\right] \nonumber \\
    =&\ \Prob\left[\Prob[\Tcal_{k} < \Tcal_{k-1} | E, G_{[k-1]}] \in \iota(p_k) \big| E, S_{k-1} = w_{k}, G_{[k-2]} \in N(p_{[k-1]},w_{[k-1]})\right] \nonumber \\
    &\ \cdot \Prob\left[S_{k-1} = w_{k} \land G_{[k-2]} \in N(p_{[k-1]},w_{[k-1]}) \big| E\right] \nonumber \\
    =&\ \Prob\left[\Prob[\Tcal_{k} < \Tcal_{k-1} | E, G_{[k-1]}] \in \iota(p_k) \big| E, S_{k-1} = w_{k}, G_{[k-2]} \in N(p_{[k-1]},w_{[k-1]})\right] \nonumber \\
    &\ \cdot\Prob\left[S_{k-1} = w_{k} \big| E, G_{[k-2]} \in N(p_{[k-1]},w_{[k-1]})\right]
    \Prob\left[G_{[k-2]} \in N(p_{[k-1]},w_{[k-1]}) \big| E\right] \nonumber \\
    =&\ \Prob[G_{-1} \in N(p_0,w_0)] \nonumber \\
    &\ \cdot \prod_{i=1}^k \Big(\Prob\left[\Prob[\Tcal_{i} < \Tcal_{i-1} | E, G_{[i-1]}] \in \iota(p_i) \big| E, S_{i-1} = w_{i}, G_{[i-2]} \in N(p_{[i-1]},w_{[i-1]})\right] \nonumber \\
    &\hspace{28pt} \cdot \Prob\left[S_{i-1} = w_{i} \big| E, G_{[i-2]} \in N(p_{[i-1]},w_{[i-1]})\right]\Big) \nonumber \\
    =&\ \prod_{i=1}^k \Big(\Prob\left[\Prob[\Tcal_{i} < \Tcal_{i-1} | E, G_{[i-1]}] \in \iota(p_i) \big| E, S_{i-1} = w_{i}, G_{[i-2]} \in N(p_{[i-1]},w_{[i-1]})\right] \nonumber \\
    &\ \cdot \Prob\left[S_{i-1} = w_{i} \big| E, G_{[i-2]} \in N(p_{[i-1]},w_{[i-1]})\right]\Big). \label{eq:smax1proof1}
  \end{align}
  Our goal is to bound \eqref{eq:smax1proof1}.  Toward this end, define
  \bequationNN
    \Idec := \{i \in \{1,\dots,k\} : w_{i} = w_{i-1} + 1\}\ \ \text{and}\ \ \Idec^c := \{i \in \{1,\dots,k\} : w_{i} = w_{i-1}\},
  \eequationNN
  which by the definition of $w_{[k]}$ form a partition of $\{1,\dots,k\}$. For any $i \in \Idec$,
  \begin{align*}
    &\ \Prob[S_{i-1} = w_{i} | E, G_{[i-2]} \in N(p_{[i-1]},w_{[i-1]})] \\
    =&\ \Prob[\Tcal_{i-1} < \Tcal_{i-2} |E,  G_{[i-2]} \in N(p_{[i-1]},w_{[i-2]})] \leq p_{i-1} + \tfrac{1}{B}.
  \end{align*}
  On the other hand, for any $i \in \Idec^c$,
  \begin{align*}
    &\ \Prob[S_{i-1} = w_{i} | E, G_{[i-2]} \in N(p_{[i-1]},w_{[i-1]})] \\
    =&\ \Prob[\Tcal_{i-1} = \Tcal_{i-2} |E,  G_{[i-2]} \in N(p_{[i-1]},w_{[i-1]})] \\
    =&\ 1 - \Prob[\Tcal_{i-1} < \Tcal_{i-2} |E,  G_{[i-2]} \in N(p_{[i-1]},w_{[i-1]})] \leq 1-p_{i-1}.
  \end{align*}
  Thus, it follows that the latter term in \eqref{eq:smax1proof1} satisfies
  \bequationNN
    \prod_{i=1}^{k} \Prob[S_{i-1} = w_{i} | E, G_{[i-2]} \in N(p_{[i-1]},w_{[i-1]})] \leq \(\prod_{i \in \Idec} (p_{i-1} + \tfrac{1}{B})\) \(\prod_{i \in \Idec^c} (1-p_{i-1})\).
  \eequationNN
  Now let us bound this term.  By the definition of $\Lbad$, one finds that
  \bequation\label{eq.badbadbad}
    \sum_{i=0}^k p_i > \ell(\smax,\hat{\delta}) + 1 \implies \sum_{i=0}^{k-1} p_i > \ell(\smax,\hat{\delta}).
  \eequation
  In addition, by the definition of $\smax$, it follows that $w_k \leq \smax$ for all $k \in [\kmax]$, from which it follows that $|\Idec| \leq \smax$.  Now, let $\{Z_0,\dots,Z_{k-1}\}$ be independent Bernoulli random variables such that, for all $i \in \{0,\dots,k-1$\}, one has
  \begin{equation*}
    \Prob[Z_i=1] = \begin{cases} p_i + \frac{1}{B} & \text{if} \; i+1 \in \Idec \\
    p_i & \text{if} \; i+1 \in \Idec^c.
    \end{cases}
  \end{equation*}
  By \eqref{eq.badbadbad}, it follows from the definition of these random variables that $\sum_{i=0}^{k-1} \Prob[Z_i=1] \geq \ell(\smax, \hat{\delta})$. Then, it follows by Lemma \ref{lem:chernoff} and the preceding argument that
  \begin{align*}
    &\prod_{i \in \Idec} \left(p_{i-1} + \frac{1}{B}\right) \prod_{i \in \Idec^c} (1-p_{i-1}) \\
    &= \Prob\left[(\text{$Z_{i-1} = 1$ for all $i \in \Idec$}) \land (\text{$Z_{i-1} = 0$ for all $i \in \Idec^c$}) \right] \\
    &= \Prob\left[\(\sum_{i=0}^{k-1} Z_i \leq \smax\) \land (\text{$Z_{i-1} = 1$ for all $i \in \Idec$}) \land (\text{$Z_{i-1} = 0$ for all $i \in \Idec^c$}) \right] \\
    &\leq \Prob\left[\sum_{i=0}^{k-1} Z_i \leq \smax\right] \leq \hat{\delta}.
  \end{align*}
  Combining this with \eqref{eq:smax1proof1}, the desired conclusion follows. \qed
\eproof

Next, we present a lemma about nodes in the sets defined in \eqref{eq:Cdec} and \eqref{eq:CdecC}.  The lemma essentially states that a certain probability of interest, defined as the product of probabilities along a path to a child node, can be reduced to a product of probabilities to the child's parent node by partitioning the childen into those at which a merit parameter decrease has occurred and children at which a merit parameter decrease has not occurred.

\blemma \label{lem:cdec}
  For all $k \in [\kmax]$ and $(p_{[k]},w_{[k]})$, one finds that
  \begin{align*}
    &\sum_{\{(p_{k+1},w_{k+1}) : N(p_{[k+1]},w_{[k+1]}) \in \Cdec(N(p_{[k]},w_{[k]})) \}} \\
    &\prod_{i=1}^{k+1} \Prob\left[\Prob[\Tcal_{i} < \Tcal_{i-1} | E, G_{[i-1]}] \in \iota(p_i) \big| E, S_{i-1} = w_{i}, G_{[i-2]} \in N(p_{[i-1]},w_{[i-1]})\right] \\
    &\quad= \prod_{i=1}^k \Prob\left[\Prob[\Tcal_{i} < \Tcal_{i-1} | E, G_{[i-1]}] \in \iota(p_i) \big| E, S_{i-1} = w_{i}, G_{[i-2]} \in N(p_{[i-1]},w_{[i-1]})\right]
  \end{align*}
  and, similarly, one finds that
  \begin{align*}
    &\sum_{\{(p_{k+1},w_{k+1}) : N(p_{[k+1]},w_{[k+1]}) \in \Cdec^c(N(p_{[k]},w_{[k]})) \}} \\
    &\prod_{i=1}^{k+1} \Prob\left[\Prob[\Tcal_{i} < \Tcal_{i-1} | E, G_{[i-1]}] \in \iota(p_i) \big| E, S_{i-1} = w_{i}, G_{[i-2]} \in N(p_{[i-1]},w_{[i-1]})\right] \\
    &\quad= \prod_{i=1}^k \Prob\left[\Prob[\Tcal_{i} < \Tcal_{i-1} | E, G_{[i-1]}] \in \iota(p_i) \big| E, S_{i-1} = w_{i}, G_{[i-2]} \in N(p_{[i-1]},w_{[i-1]})\right],
  \end{align*}
  where by the definitions of $\Cdec$ and $\Cdec^c$ it follows that the value of $w_{k+1}$ in the sum in the former equation is one greater than the value of $w_{k+1}$ in the sum in the latter equation.
\elemma
\bproof
  One finds that
  \begin{align*}
    &\sum_{\{(p_{[k+1]},w_{[k+1]}) : N(p_{[k+1]},w_{[k+1]}) \in \Cdec(N(p_{[k]},w_{[k]}))\}} \\
    &\prod_{i=1}^{k+1} \Prob\left[\Prob[\Tcal_{i} < \Tcal_{i-1} | E, G_{[i-1]}] \in \iota(p_i) \big| E, S_{i-1} = w_{i}, G_{[i-2]} \in N(p_{[i-1]},w_{[i-1]})\right] \\
    &= \prod_{i=1}^{k} \Prob\left[\Prob[\Tcal_{i} < \Tcal_{i-1} | E, G_{[i-1]}] \in \iota(p_i) \big| E, S_{i-1} = w_{i}, G_{[i-2]} \in N(p_{[i-1]},w_{[i-1]})\right] \\
    &\qquad \cdot \sum_{\{(p_{[k+1]},w_{[k+1]}) : N(p_{[k+1]},w_{[k+1]}) \in \Cdec(N(p_{[k]},w_{[k]}))\}} \\
    &\qquad\qquad \Prob\left[\Prob[\Tcal_{k+1} < \Tcal_{k} | E, G_{[k]}] \in \iota(p_{k+1}) \big| E, S_{k} = w_{k+1}, G_{[k-1]} \in N(p_{[k]},w_{[k]})\right].
  \end{align*}
  The desired conclusion follows since, by the definition of $\Cdec(N(p_{[k]},w_{[k]}))$, all elements in the latter sum have $S_k = w_{k+1} = w_{k}+1$, meaning that the sum on the right-hand side is the sum of all conditional probabilities with the same conditions, and hence the sum is 1.
  
  The proof of the second desired conclusion follows in the same manner with $\Cdec^c$ in place of $\Cdec$ and $S_k=w_{k+1}=w_{k}$ in place of $S_k=w_{k}=w_{k-1}+1$. \qed
\eproof

Next, we derive a result for certain nodes containing realizations with $w_{k} = \smax - 1$.

\blemma \label{lem:smaxw1}
  For any $k \in [\kmax]$ and $(p_{[k]},w_{[k]})$ such that $w_{k} = s_{\max} - 1$ and $N(p_{[k]},w_{[k]}) \not\in \Lgood$, it follows that
  \begin{align}
    &\Prob[G_{[k-1]} \in N(p_{[k]},w_{[k]}) \land \EbadB | E] \nonumber \\
    &\leq \hat{\delta} \prod_{i=1}^k \Prob\left[\Prob[\Tcal_{i} < \Tcal_{i-1} | E, G_{[i-1]}] \in \iota(p_i) \big| E, S_{i-1} = w_{i}, G_{[i-2]} \in N(p_{[i-1]},w_{[i-1]})\right].  \label{eq:recursivebound1}
  \end{align}
\elemma
\bproof
  By the supposition that $N(p_{[k]},w_{[k]}) \not\in \Lgood$, it follows that any $(p_{[k]},w_{[k]})$ with $h(N(p_{[k]},w_{[k]})) = 0$ has $N(p_{[k]},w_{[k]}) \in \Lbad$, in which case the desired conclusion follows from Lemma~\ref{lem:leafnodes}.  With this base case being established, we now prove the result by induction.  Suppose that the result holds for all $(p_{[k]},w_{[k]})$ such that $w_{k} = s_{\max} - 1$, $N(p_{[k]},w_{[k]}) \not\in \Lgood$, and $h(N(p_{[k]},w_{[k]})) \leq j$ for some $j \in \N{}$.  Our goal is to show that the same statement holds with $j$ replaced by $j+1$.  For this purpose, consider arbitrary $(p_{[k]},w_{[k]})$ such that $w_{k} = s_{\max} - 1$, $N(p_{[k]},w_{[k]}) \not\in \Lgood$, and $h(N(p_{[k]},w_{[k]})) = j + 1$.  Observe that by the definition of the child operators $C$, $\Cdec$, and $\Cdec^c$, it follows that
  \begin{align*}
    &\Prob[G_{[k-1]} \in N(p_{[k]},w_{[k]}) \land \EbadB | E] \\
    &= \sum_{\{(p_{k+1},w_{k+1}) : N(p_{[k+1]},w_{[k+1]}) \in C(N(p_{[k]},w_{[k]})) \}} \Prob[G_{[k]} \in N(p_{[k+1]},w_{[k+1]}) \land \EbadB | E] \\
    &= \sum_{\{(p_{k+1},w_{k+1}) : N(p_{[k+1]},w_{[k+1]}) \in \Cdec(N(p_{[k]},w_{[k]})) \}} \Prob[G_{[k]} \in N(p_{[k+1]},w_{[k+1]}) \land \EbadB | E] \\
    &+ \sum_{\{(p_{k+1},w_{k+1}) : N(p_{[k+1]},w_{[k+1]}) \in \Cdec^c(N(p_{[k]},w_{[k]})) \}} \Prob[G_{[k]} \in N(p_{[k+1]},w_{[k+1]}) \land \EbadB | E].
  \end{align*}
  Since $w_{k} = \smax-1$, it follows from the definition of $\Cdec$ that for any $(p_{k+1},w_{k+1})$ with $N(p_{[k+1]},w_{[k+1]}) \in \Cdec(N(p_{[k]},w_{[k]}))$, one finds that $w_{k+1} = w_{k}+1 = \smax$. By the definition of $\smax$, this implies that $\Prob[\Tcal_{k+1} < \Tcal_k | E, G_{[k]}] = 0$, so $p_{k+1} = 0$. In addition, since $N(p_{[k]},w_{[k]}) \not\in \Lbad$ since $C(N(p_{[k]},w_{[k]})) \neq \emptyset$, it follows that $\sum_{i=0}^{k+1} p_{k+1} \leq \ell(\smax,\hat{\delta})+1$, meaning $N(p_{[k+1]},w_{[k+1]}) \in \Lgood$.  Consequently, from above and Lemma \ref{lem:leafnodes}, one finds
  \begin{align*}
    &\Prob[G_{[k-1]} \in N(p_{[k]},w_{[k]}) \land \EbadB | E] \\
    &= \sum_{\{(p_{k+1},w_{k+1}) : N(p_{[k+1]},w_{[k+1]}) \in \Cdec^c(N(p_{[k]},w_{[k]})) \}} \Prob[G_{[k]} \in N(p_{[k+1]},w_{[k+1]}) \land \EbadB | E].
  \end{align*}
  Since $h(N(p_{[k]},w_{[k]}) = j+1$, it follows that $h(N(p_{[k+1]},w_{[k+1]})) \leq j$ for any $(p_{[k+1]},w_{[k+1]})$ with $h(N(p_{[k+1]},w_{[k+1]})) \in \Cdec^c(N(p_{[k]},w_{[k]}))$.  Therefore, by the induction hypothesis and the result of Lemma~\ref{lem:cdec}, it follows that
  \begin{align*}
    &\Prob[G_{[k-1]} \in N(p_{[k]},w_{[k]}) \land \EbadB | E] \\
    &= \sum_{\{(p_{k+1},w_{k+1}) : N(p_{[k+1]},w_{[k+1]}) \in \Cdec^c(N(p_{[k]},w_{[k]})) \}} \\
    &\quad \hat{\delta} \prod_{i=1}^{k+1} \Prob\left[\Prob[\Tcal_{i} < \Tcal_{i-1} | E, G_{[i-1]}] \in \iota(p_i) \big| E, S_{i-1} = w_{i}, G_{[i-2]} \in N(p_{[i-1]},w_{[i-1]})\right] \\
    &\quad\leq \hat{\delta} \prod_{i=1}^k \Prob\left[\Prob[\Tcal_{i} < \Tcal_{i-1} | E, G_{[i-1]}] \in \iota(p_i) \big| E, S_{i-1} = w_{i}, G_{[i-2]} \in N(p_{[i-1]},w_{[i-1]})\right],
  \end{align*}
  which completes the proof. \qed
\eproof

Using the preceding lemma as a base case, we now perform induction on the difference $\smax - w_{k}$ to prove a similar result for arbitrary $\smax$.

\blemma \label{lem:nodeprobfull}
  For any $k \in [\kmax]$ and $(p_{[k]},w_{[k]})$ with $N(p_{[k]},w_{[k]}) \not\in \Lgood$, it follows that
  \begin{align*}
    &\Prob[G_{[k-1]} \in N(p_{[k]},w_{[k]}) \land \EbadB | E] \\
    &\leq \hat{\delta} \cdot \sum_{j=0}^{\min\{\smax-w_{k}-1,h(N(p_{[k]},w_{[k]}))\}} \( \begin{matrix} h(N(p_{[k]},w_{[k]})) \\ j \end{matrix} \) \\
    &\qquad \cdot \prod_{i=1}^k \Prob\left[\Prob[\Tcal_{i} < \Tcal_{i-1} | E, G_{[i-1]}] \in \iota(p_i) \big| E, S_{i-1} = w_{i}, G_{[i-2]} \in N(p_{[i-1]},w_{[i-1]})\right].
  \end{align*}
\elemma
\bproof
  For all $(p_{[k]},w_{[k]})$ such that $N(p_{[k]},w_{[k]}) \not\in \Lgood$ and $h(N(p_{[k]},w_{[k]})) = 0$, it follows that $N(p_{[k]},w_{[k]}) \in \Lbad$.  The result holds in this case according to Lemma~\ref{lem:leafnodes} since one finds that $\sum_{j=0}^{\min\{\smax-w_{k}-1,h(N(p_{[k]},w_{[k]}))\}} \binom{h(N(p_{[k]},w_{[k]}))}{j} = \binom{0}{0} = 1$.  On the other hand, for all $(p_{[k]},w_{[k]})$ such that $N(p_{[k]},w_{[k]}) \not\in \Lgood$ and $\smax - w_{k} = 1$, the result follows from Lemma~\ref{lem:smaxw1}.  Hence, to prove the remainder of the result by induction, one may assume that it holds for all $(p_{[k]},w_{[k]})$ such that $N(p_{[k]},w_{[k]}) \not\in \Lgood$, $h(N(p_{[k]},w_{[k]})) \leq t$ for some $t \in \N{}$, and $\smax - w_{k} = r$ for some $r \in \N{} \setminus \{0\}$, and show that it holds for all $(p_{[k]},w_{[k]})$ such that $N(p_{[k]},w_{[k]}) \not\in \Lgood$, $h(N(p_{[k]},w_{[k]})) = t + 1$, and $\smax - w_{k} = r$.
  
  Consider arbitrary $(p_{[k]},w_{[k]})$ such that $N(p_{[k]},w_{[k]}) \not\in \Lgood$, $h(N(p_{[k]},w_{[k]})) = t + 1$, and $\smax - w_{k} = r$.  By the definitions of $C$, $\Cdec$, and $\Cdec^c$, it follows that
  \begin{align*}
    &\Prob[G_{[k-1]} \in N(p_{[k]},w_{[k]}) \land \EbadB | E] \\
    &= \sum_{\{(p_{[k+1]},w_{[k+1]}) : N(p_{[k+1]},w_{[k+1]}) \in C(N(p_{[k]},w_{[k]})) \}} \Prob[G_{[k]} \in N(p_{[k+1]},w_{[k+1]}) \land \EbadB | E] \\
    &= \sum_{\{(p_{[k+1]},w_{[k+1]}) : N(p_{[k+1]},w_{[k+1]}) \in \Cdec(N(p_{[k]},w_{[k]})) \}} \Prob[G_{[k]} \in N(p_{[k+1]},w_{[k+1]}) \land \EbadB | E] \\
    &+ \sum_{\{(p_{[k+1]},w_{[k+1]}) : N(p_{[k+1]},w_{[k+1]}) \in \Cdec^c(N(p_{[k]},w_{[k]})) \}} \Prob[G_{[k]} \in N(p_{[k+1]},w_{[k+1]}) \land \EbadB | E].
  \end{align*}
  Further by the definition of $\Cdec$, it follows that $w_{k+1} = w_k + 1$ (thus $\smax - w_{k+1} = r-1$) for all terms in the former sum on the right-hand side, whereas by the definition of $\Cdec^c$ it follows that $w_{k+1} = w_k$ (thus $\smax - w_{k+1} = r$) for all terms in the latter sum on the right-hand side.  Moreover, from $h(N(p_{[k]},w_{[k]})) = t+1$, it follows that $h(N(p_{[k+1]},w_{[k+1]})) \leq t$ for all terms on the right-hand side.  Therefore, by the induction hypothesis, it follows that
  \begin{align*}
    &\Prob[G_{[k-1]} \in N(p_{[k]},w_{[k]}) \land \EbadB | E] \\
    &\leq \sum_{\{(p_{[k+1]},w_{[k+1]}) : N(p_{[k+1]},w_{[k+1]}) \in \Cdec(N(p_{[k]},w_{[k]})) \}} \hat{\delta} \sum_{j=0}^{\min\{r-2,t\}} \binom{t}{j} \\
    &\quad \cdot \prod_{i=1}^{k+1} \Prob\left[\Prob[\Tcal_{i} < \Tcal_{i-1} | E, G_{[i-1]}] \in \iota(p_{i}) \big| E, S_{i-1} = w_i, G_{[i-2]} \in N(p_{[i-1]},w_{[i-1]})\right] \\
    &+ \sum_{\{(p_{[k+1]},w_{[k+1]}) : N(p_{[k+1]},w_{[k+1]}) \in \Cdec^c(N(p_{[k]},w_{[k]})) \}} \hat{\delta} \sum_{j=0}^{\min\{r-1,t\}} \binom{t}{j} \\
    &\quad \cdot \prod_{i=1}^{k+1} \Prob\left[\Prob[\Tcal_{i} < \Tcal_{i-1} | E, G_{[i-1]}] \in \iota(p_{i}) \big| E, S_{i-1} = w_i, G_{[i-2]} \in N(p_{[i-1]},w_{[i-2]})\right],
  \end{align*}
  which by Lemma \ref{lem:cdec} implies that
  \begin{align}
    &\Prob[G_{[k-1]} \in N(p_{[k]},w_{[k]}) \land \EbadB | E] \nonumber \\
    &\leq \hat{\delta} \left(\sum_{j=0}^{\min\{r-2,t\}} \binom{t}{j} + \sum_{j=0}^{\min\{r-1,t\}} \binom{t}{j}\right) \label{eq:nodeprobfull1} \\
    &\quad \cdot \prod_{i=1}^{k} \Prob\left[\Prob[\Tcal_{i} < \Tcal_{i-1} | E, G_{[i-1]}] \in \iota(p_i) \big| E, S_{i-1} = w_i, G_{[i-2]} \in N(p_{[i-1]},w_{[i-1]})\right]. \nonumber
  \end{align}
  To complete the proof, we need only consider two cases on the relationship between $t$ and~$r$.  First, if $t \leq r-2$, then Pascal's rule implies that
  \begin{align*}
    \sum_{j=0}^{\min\{r-2,t\}} \binom{t}{j} + \sum_{j=0}^{\min\{r-1,t\}} \binom{t}{j} &= 2\sum_{j=0}^{t} \binom{t}{j} \\
    &= \binom{t}{t} + \binom{t}{0} + \sum_{j=1}^t \left(\binom{t}{j} + \binom{t}{j-1}\right) \\
    &= \binom{t+1}{t+1} + \binom{t+1}{0} + \sum_{j=1}^t \binom{t+1}{j} \\
    &= \sum_{j=0}^{t+1} \binom{t+1}{j} = \sum_{j=0}^{h(N_{p_{[k]},w_{[k]}})}  \binom{h(N_{p_{[k]},w_{[k]}})}{j}.
  \end{align*}
  Since $t \leq r-2$, it follows that $h(N_{p_{[k]},w_{[k]}}) = t+1 \leq r-1 = \smax-w_{k}-1$, which combined with \eqref{eq:nodeprobfull1} proves the result in this case.  Second, if $t > r-2$, then similarly
  \begin{align*}
    \sum_{j=0}^{\min\{r-2,t\}} \binom{t}{j} + \sum_{j=0}^{\min\{r-1,t\}} \binom{t}{j} &= \sum_{j=0}^{r-2} \binom{t}{j} + \sum_{j=0}^{r-1} \binom{t}{j} \\
    &= \binom{t}{0} + \sum_{j=1}^{r-1} \left(\binom{t}{j} + \binom{t}{j-1}\right) \\
    &= \binom{t+1}{0} + \sum_{j=1}^{r-1} \binom{t+1}{j} \\
    &= \sum_{j=0}^{r-1} \binom{t+1}{j} = \sum_{j=0}^{\smax-w_{k}-1} \binom{h(N_{p_{[k]},w_{[k]}})}{j}.
  \end{align*}
  Since $t > r-2$, $h(N_{p_{[k]},w_{[k]}}) = t+1 > r-1 = \smax-w_{k-1}-1$, which combined with \eqref{eq:nodeprobfull1} proves the result for this case as well. \qed
\eproof

We now prove our first main result of this section.

\begin{theorem} \label{coro:probbounded}
  For any $\delta \in (0,1)$ with $\hat{\delta}$ as defined in \eqref{eq:hatdelta} and $\ell(\smax,\hat{\delta})$ as defined in \eqref{eq:lsd}, one finds that \eqref{eq:sumptaudecbound} holds.
\end{theorem}
\bproof
  First, consider the case where $\smax = 0$. Then, by the definition of $\smax$,
  \bequationNN
    \Prob[\Tcal_k < \Tcal_{k-1}|E,G_{[k-1]}] = 0,
  \eequationNN
  for all $k = [\kmax]$, so the result holds trivially.
  
  Now, let $\smax \in \N{}_{>0}$. By construction of our tree and the definitions of $\Lgood$ and $\Lbad$, one finds that $h(N(p_0,w_0)) \leq \kmax$. In addition, by the definition of $\smax$, $\smax - 1 < \kmax$, so $\min\{\smax-w_0-1,h(N(p_0,w_0))\} = \smax-1 \geq 0$.  Consider arbitrary $B \in \N{} \setminus \{0\}$ (see \eqref{eq:EbadB}).  By Lemma \ref{lem:nodeprobfull} and \eqref{eq:hatdelta},
  \bequationNN
    \Prob[\EbadB | E] = \Prob[G_{[-1]} \in N(p_0,w_0) \land \EbadB | E] \leq \hat{\delta} \sum_{j=0}^{\min\{\smax-1,\kmax\}} \binom{\kmax}{j} = \delta.
  \eequationNN
  Therefore, by the definition of $\EbadB$ (see \eqref{eq:EbadB}), it follows that
  \bequationNN
    \Prob\left[\sum_{i=0}^{\kmax} \Prob[\Tcal_i < \Tcal_{i-1} | E, G_{[i-1]}] \leq \ell(s_{\max},\hat{\delta}) + \tfrac{\kmax+1}{B}+1 \Bigg| E \right] \geq 1-\delta.
  \eequationNN
  Now, let us define the event $E_{\text{good},B}$ for $B \in \N{} \setminus \{0\}$ as the event that
  \bequationNN
    \sum_{i=0}^{\kmax} \Prob[\Tcal_i < \Tcal_{i-1} | E, G_{[i-1]}] \leq \ell(s_{\max},\hat{\delta}) + \tfrac{\kmax+1}{B}+1,
  \eequationNN
  One sees that $E_{\text{good},B} \supseteq E_{\text{good},B+1}$ for all such $B$.  Therefore, by the properties of a decreasing sequence of events (see, for example \cite[Section 1.5]{DStirzaker_2003}), it follows that
  \begin{align*}
    &\Prob\left[\sum_{i=0}^{\kmax} \Prob[\Tcal_i < \Tcal_{i-1} | E, G_{[i-1]}] \leq \ell(s_{\max},\hat{\delta}) +1 \Bigg| E \right] \\
    &= \Prob\left[\lim_{B \rightarrow \infty} \right(\sum_{i=0}^{\kmax} \Prob[\Tcal_i < \Tcal_{i-1} | E, G_{[i-1]}] \leq \ell(s_{\max},\hat{\delta}) + \frac{\kmax+1}{B}+1\left) \Bigg| E \right] \\
    &= \lim_{B \rightarrow \infty} \Prob\left[\sum_{i=0}^{\kmax} \Prob[\Tcal_i < \Tcal_{i-1} | E, G_{[i-1]}] \leq \ell(s_{\max},\hat{\delta}) + \frac{\kmax+1}{B}+1 \Bigg| E \right] \geq 1-\delta,
  \end{align*}
  as desired.  \qed
\eproof

Now, we are prepared to prove Lemma \ref{lem:Eksumlsmax}.

\begin{proof}[Lemma~\ref{lem:Eksumlsmax}]
  Observe that, for any $k \in [\kmax]$, by the defintion of $\Ek3$, the event $\Tcal_k < \Tcal_{k-1}$ must occur whenever $\Ek3$ occurs. Therefore, for any $k \in [\kmax]$, one finds
  \bequationNN
    \Prob[\Tcal_k < \Tcal_{k-1} | E, \Fcal_k] \geq \Prob[\Ek3 | E, \Fcal_k].
  \eequationNN
  The result then follows directly from Theorem~\ref{coro:probbounded}. \qed
\end{proof}

Now, we turn our attention to Lemma \ref{lem:Ktaulsmax}.  Let a realization of the random index set $\Kcal_{\tau}$ defined in \eqref{eq:Ktau} be denoted by $\ktau$.  Our next lemma shows an important property about any iteration $k \in [\kmax]$ in which $k \in \ktau$ for a given realization $k_\tau$.

\blemma \label{lem:Etautrial}
  For any $k \in [\kmax]$ and $g_{[k-1]}$ such that $k \in k_\tau$ for some realization $k_\tau$ of the random index set $\Kcal_{\tau}$, one finds that $\Prob[\Tcal_k < \tau_{k-1} | E, g_{[k-1]}, k \in \ktau] \geq p_{\tau}$.
\elemma
\bproof
  In any iteration during which $\tktritrue < \tau_{k-1}$, it follows that $\tktritrue < \infty$, so
    \bequationNN
        \tktritrue = \frac{(1-\sigma)\|c_k\|_1}{\nabla f(x_k)^\top \dktrue + \max\{(\dktrue)^\top H_k \dktrue,0\}}
    \eequationNN
    and thus
    \bequationNN
        (1-\sigma)\|c_k\|_1 < (\nabla f(x_k)^\top \dktrue + \max\{(\dktrue)^\top H_k \dktrue,0\}) \tau_{k-1}.
    \eequationNN
    By the definition of $\tau_k$, if
    \bequationNN
        g_k^\top d_k + \max\{d_k^\top H_k d_k,0\} \geq \nabla f(x_k)^\top \dktrue + \max\{(\dktrue)^\top H_k \dktrue,0\}
    \eequationNN
    in an iteration such that $\tktritrue < \tau_{k-1}$, then
    \bequationNN
        (1-\sigma)\|c_k\|_1 < (g_k^\top d_k + \max\{d_k^\top H_k d_k,0\}) \tau_{k-1},
    \eequationNN
    meaning that $\tau_k < \tau_{k-1}$.  Noting that the event $k \in \ktau$ is conditionally independent of $G_k$ given $E$ and $g_{[k-1]}$, it follows from Assumption~\ref{assum:ptau} that
    \begin{align*}
        &\Prob[\Tcal_k < \tau_{k-1} | E, g_{[k-1]}, k \in \ktau] \\
        &\geq \Prob[G_k^\top D_k + \max\{D_k^\top H_k D_k,0\} \geq \nabla f(x_k)^\top \dktrue + \max\{(\dktrue)^\top H_k \dktrue,0\} | E, g_{[k-1]}, k \in \ktau] \\
        &= \Prob[G_k^\top D_k + \max\{D_k^\top H_k D_k,0\} \geq \nabla f(x_k)^\top \dktrue + \max\{(\dktrue)^\top H_k \dktrue,0\} | E, g_{[k-1]}] \geq p_{\tau},
    \end{align*}
    as desired. \qed
\eproof

The previous lemma guarantees that in any iteration in which $\tktritrue < \tau_{k-1}$, the probability is at least $p_{\tau}$ that the merit parameter decreases. By scheme for setting $\tau_k$,
\bequation \label{eq:cardinality1}
    \Prob[\tktritrue < \tau_k | E, g_{[k-1]}, \tktritrue \geq \tau_{k-1}] = 0,
\eequation
so one must have $\tktritrue < \tau_{k-1}$ in any iteration when $\hat{\tau}_k < \tau_k$.  Thus, we can obtain a bound on the number of iterations at which $\tktritrue < \tau_k$ by bounding the number of iterations at which $\tktritrue < \tau_{k-1}$. Now we prove a result relating $|\Kcal_{\tau}|$ to the probabilities of decreasing the merit paremeter over all iterations.

\blemma \label{lem:probktaudec}
  One finds that
  \bequationNN
    \sum_{k=0}^{\kmax} \Prob[\Tcal_k < \Tcal_{k-1} | E, \Fcal_k] \geq |\Kcal_{\tau}| p_{\tau}.
  \eequationNN
\elemma
\bproof
  Consider arbitrary $g_{[\kmax]}$ and the corresponding realization $k_\tau$ of the random index set $\Kcal_\tau$.  By Lemma \ref{lem:Etautrial}, it follows for all $k \in [\kmax]$ that
  \begin{align*}
    \Prob[\Tcal_k < \tau_{k-1} | E, g_{[k-1]}] 
    &\geq \Prob[\Tcal_k < \tau_{k-1} \land k \in \ktau | E, g_{[k-1]}] \\
    &= \Prob[\Tcal_k < \tau_{k-1} | E, g_{[k-1]}, k \in \ktau] \cdot \Prob[k \in \ktau | E, g_{[k-1]}] \\
    &= \Prob[\Tcal_k < \tau_{k-1} | E, g_{[k-1]}, k \in \ktau] \cdot \Ical[k \in \ktau] \\
    &\geq \Ical[k \in \ktau] p_{\tau},
  \end{align*}
  where $\Ical[k \in \ktau]$ is the indicator function for the event $k \in \ktau$ and the second equality follows due to the fact that the event $k \in \ktau$ is deterministically known when conditioned on $g_{[k-1]}$. Summing this inequality over $k \in [\kmax]$, one finds that
    \bequationNN
        \sum_{k=0}^{\kmax} \Prob[\Tcal_k < \tau_{k-1} | E, g_{[k-1]}] \geq \sum_{k=0}^{\kmax} \Ical[k \in \ktau] p_{\tau} = |\ktau| p_{\tau}.
    \eequationNN
    Letting $f_{G_{[\kmax]}}$ denote the probability density function of $G_{[\kmax]}$, the fact that the bound above holds \emph{deterministically} for any realization $g_{[\kmax]}$ that
    \begin{align*}
        &\Prob\left[\sum_{k=0}^{\kmax} \Prob[\Tcal_k < \Tcal_{k-1} | E, \Fcal_{k-1}] \geq |\Kcal_{\tau}| p_{\tau} \ \Bigg| E\right] \\
        &=\int_{g_{[\kmax]} \in \Fcal_{\kmax+1}} \Prob\left[\sum_{k=0}^{\kmax} \Prob[\Tcal_k < \tau_{k-1} | E, g_{[k-1]}] \geq |\ktau| p_{\tau} \ \Bigg| E, g_{[\kmax]}\right] f_{G_{[\kmax]}}(g_{[k-1]}) \textrm{d} g_{[\kmax]} \\
        &= \int_{g_{[\kmax]} \in \Fcal_{\kmax+1}} 1 \cdot f_{G_{[\kmax]}}(g_{[k-1]}) \textrm{d} g_{[\kmax]} = 1.
    \end{align*}
    Therefore, the desired result holds as well. \qed
\eproof

We now claim that Lemma \ref{lem:Ktaulsmax} follows.

\begin{proof}[Lemma \ref{lem:Ktaulsmax}]
    The proof follows by combining Theorem \ref{coro:probbounded} and Lemma \ref{lem:probktaudec}. \qed
\end{proof}

We conclude this appendix by showing that the order notation result in \eqref{eq.like_sg} and \eqref{eq.new_terms} holds, as required in the proof of Corollary~\ref{coro:complexity}.

\blemma \label{lem:Orderlsd}
    Let $\delta \in (0,1)$, $\hdelta$ be defined in \eqref{eq:hatdelta}, $\smax \in \N{}_{>0}$ and $\ell(\smax,\hdelta)$ be defined in \eqref{eq:lsd}. Then,
    \bequationNN
        \ell(\smax,\hdelta) = \Ocal\left(\smax\log(\kmax) + \log(1/\delta)\right).
    \eequationNN
\elemma 
\bproof
  Since $\smax \in \N{}_{>0}$, it follows
    \begin{align*}
        \sum_{j=0}^{\max\{\smax-1,0\}} \binom{\kmax}{j}
        &= \sum_{j=0}^{\smax-1} \frac{(\kmax)!}{j!(\kmax-j)!} \\
        &\leq \sum_{j=0}^{\smax-1} \frac{(\kmax)!}{(\kmax-j)!} \\
        &= 1+\sum_{j=1}^{\smax-1} \prod_{i=\kmax+1-j}^{\kmax} i \\
        &\leq 1+\sum_{j=1}^{\smax-1} (\kmax)^j \\
        &\leq 1+(\smax-2) (\kmax)^{\smax-1} \\
        &\leq (\smax-1) (\kmax)^{\smax-1}.
    \end{align*}
    Then, by the definitions of $\ell(\smax,\hdelta)$ and $\hdelta$, it follows that
    \begin{align*}
        \ell(\smax,\hdelta) &= \Ocal\left(\smax + \log(1/\hdelta)\right) \\
        &= \Ocal\left(\smax + \log(\smax-1) + (\smax-1)\log(\kmax) + \log(1/\delta)\right) \\
        &= \Ocal\left(\smax\log(\kmax) + \log(1/\delta)\right),
    \end{align*}
    as desired. \qed
\eproof

\end{document}